\documentclass[final]{elsarticle}
\usepackage{xurl}
\usepackage{lineno,hyperref,xcolor}
\usepackage{amsmath} 
\usepackage{amsfonts} 
\usepackage{mathtools} 
\usepackage{amssymb} 
\usepackage{siunitx} 
\usepackage[table]{xcolor} 
\usepackage{graphicx}

\usepackage{subcaption} 
\usepackage{tikz}
\usetikzlibrary{math}
\usetikzlibrary{shapes.geometric}
\usetikzlibrary{decorations.pathreplacing}
\usetikzlibrary{arrows.meta}

\newcommand{\R}{\mathbb{R}}
\newcommand{\C}{\mathbb{C}}
\newcommand{\N}{\mathbb{N}}
\newcommand{\md}{\mathrm{d}}

\newcommand{\mi}{\mathrm{i}}

\DeclareMathOperator{\diag}{diag}
\DeclareMathOperator{\rank}{rank}

\usepackage{amsthm}
\theoremstyle{definition}
\newtheorem{definition}{Definition}[section]

\theoremstyle{plain}

\newtheorem{theorem}[definition]{Theorem}

\modulolinenumbers[1]

\begin{document}

\begin{frontmatter}

\title{Investigating parameter dependence for counterexamples to Payne’s nodal line conjecture using HPC}

\author[first]{Sina Mattfeldt\corref{mycorrespondingauthor}}
\cortext[mycorrespondingauthor]{Corresponding author}
\ead{s.mattfeldt@fz-juelich.de}
\author[first,second,third]{Andreas Kleefeld}

\address[first]{Forschungszentrum J\"ulich GmbH,
J\"ulich Supercomputing Centre, Wilhelm-Johnen-Str., 52425 J\"ulich, Germany.}

\address[second]{FH Aachen - University of Applied Sciences, Faculty of Medical Engineering and Technomathematics, Heinrich-Mu\ss{}mann-Str. 1, 52428 J\"ulich, Germany.}

\address[third]{FH Aachen - University of Applied Sciences, Institute for Data-driven Technologies, Heinrich-Mu\ss{}mann-Str. 1, 52428 J\"ulich, Germany.}

\begin{abstract}
The behavior of the nodal line of the second eigenfunction of the negative Laplacian has been widely studied. Dahne could show numerically for one domain with six holes in two dimensions that it is possible to reject Payne's nodal line conjecture. We investigate the behavior of the nodal line if it does not touch the boundary, specifically the general shape, its prediction, its dependence on the parameters, and the numbers of holes. In this paper, we could address these issues numerically with a very efficient parallelized algorithm showing the rejection for various domains with six and five holes.
\end{abstract}


\end{frontmatter}

\clearpage

\section{Introduction}
\label{sec:Introduction}
\subsection{Related work}
\label{sec:RelatedWork}
The structure of eigenfunctions can be used to describe many phenomena in physics \cite{Hoffmann_NodalLine}. In this paper, the second eigenfunction of the negative Laplacian with homogeneous Dirichlet boundary conditions is considered. For instance, eigenfunctions of the negative Laplacian can be used to describe the vibration of plates and membranes. Nodal lines are of special interest. They are the set of points in a bounded planar domain $\Omega$ where no vibration exists. Mathematically described, they are the zero sets of the $i$-th negative Laplacian eigenfunction $u_i$, given by $\mathcal{N}(u_i) = \overline{\{x \in \Omega : u_i(x) = 0 \}}$.\cite{Zelditch2017_Eigenfunctions}

The eigenfunctions are defined by the Dirichlet problem
\begin{equation}
    \begin{aligned}
    -\Delta u_i &= \lambda_i u_i, &&\qquad i = 1, 2, \dots \qquad &\mbox{in } \Omega\\
    u_i &= 0, &&\qquad i = 1, 2, \dots \qquad &\mbox{on } \Gamma
    \end{aligned}
    \label{LaplaceEigenvalueProblem}
\end{equation}
in a bounded planar domain $\Omega$ with the eigenvalues $0 < \lambda_1 < \lambda_2 \leq \lambda_3 \leq \dots$ of the Dirichlet Laplacian and the boundary $\Gamma$. Here, the first eigenvalue is a simple eigenvalue and there are no constraints on the boundary $\Gamma$ like smoothness. \cite{Dahne_Counterexample} In generally, it is not possible to say something about the multiplicity of the second eigenvalue, but with a few restrictions it is. If a smooth and bounded domain $\Omega \in \R$ is given, then the multiplicity is at most two \cite{Liqun1995_MultiplicityEigenvalue}.

This article focuses on the second eigenfunction. 
In 1967, Payne wrote the ``nodal line conjecture'' \cite{Payne1967_PaynesConjecture}. It says that on a bounded domain in $\mathbb{R}^2$ the nodal line of the second Dirichlet eigenfunction $u_2$ must touch the boundary \cite{Dahne_Counterexample}. Several proofs exist that Payne's conjecture is true under special circumstances.  In 1973, Payne proved that if $\Omega$ is convex in $x$ and symmetric about the $y$-axis, $u_2$ can not have an interior closed nodal curve \cite{Payne1973_TwoConjectures}. Here, $x,y \in \R$ are the two directions of the Cartesian coordinate system. $\Omega$ is convex in $x$, if every line parallel to the $x$-axis which intersects $\Omega$, cuts $\Gamma$ in at most two points. This result was proved by Damascelli \cite{Damascelli2000_NodalSet} in 2000 for $\mathbb{R}^n$. In 1987, Lin \cite{Lin1987_SecondEigenfunctions} was able to prove that the conjecture is true if the given smooth domain is symmetric under a rotation with angle $\frac{2\pi p}{q}$, where $p$ and $q$ are positive integers.
In 1990, Pütter \cite{Putter1990_NodalLine} could show that Payne's conjecture is fulfilled if the domain is simply-connected in $\mathbb{R}^2$ of class $C^{2,\alpha}$, symmetric about the axis $\{x_2=0\}$ and convex with respect to a point lying on this axis, but outside of the given domain $\Omega$.
Without assuming smoothness of the domain, Jerison \cite{Jerison1991_ConvexPlanarDomain} was able to show, in 1991, for a convex domain that is narrow enough, i.e., the diameter divided by the maximum inner radius of the domain is greater than an absolute constant $C$, Payne's conjecture is fulfilled.
One year later, Melas \cite{Melas1992_NodalLine} could prove that if the domain is convex and the boundary is smooth, which means $C^{\infty}$, than the conjecture is true.
Based on this proof, Alessandrini \cite{Alessandrini1994_NodalLine} were able to show this result for general convex domains without this smoothness assumptions in 1994.
In 2007 Freitas and Krej\v{c}i\v{r}ík \cite{Freitas2007LocationNodalSetCurved} proved that the conjecture is true even in non-convex domains and higher dimensions $\mathbb{R}^n$, with $n\geq 2$, if the region $\Omega$ is a sufficiently thin curved tube in $\mathbb{R}^n$.
Based on the result for convex domains in 2013, Yang and Guo \cite{Yang2013NodalLine} constructed a concave domain via an inverse polar coordinate transformation $\rho^{-1}(r,\theta) = (x,y) = (r\cos(\theta), r\sin(\theta))$. This transformation is based on a strictly convex domain in $\theta$, meaning that when varying $\theta$ while holding $r$ constant yields in strictly convex domains. Additionally, the domain must be symmetric with respect to the $r$-axes. In these concave domains the nodal line touches the boundary at exactly two points.
Most of the proofs are based on Courant's nodal domain theorem, in two dimensions also known as Courant's nodal line theorem, which states that the nodal set $u_n$ has at most $n$ nodal domains. Furthermore, eigenfunctions $u_n$ with $n \geq 2$ have at least two nodal domains due to the possibility of orthogonalizing eigenfunctions belonging to the same eigenvalue \cite{Pleijel1956_CourantTherorem}. Combined this means that the second eigenfunction $u_2$ has always two nodal domains.

On the other hand there are some counterexamples to Payne's conjecture. In 2005, Freitas and Krej\v{c}i\v{r}ík found out that there exist simply-connected, unbounded domains for which the nodal line does not touch the boundary. These domains have to be convex and symmetric with respect to two orthogonal directions. Even if the nodal line never touches the boundary, the nodal line is not closed in this counterexample. \cite{Freitas2005UnboundedNodalLine}

\begin{figure}[h!]
    \centering
    \begin{subfigure}[t]{0.28\linewidth}
        \tikzmath{\n = 5;
        \eps = 180/(\n+10);}
        \resizebox{\linewidth}{!}{
        \begin{tikzpicture}
            \definecolor{pastelgruen}{HTML}{acf6c8}
            \fill[pastelgruen] (0,0) circle (3cm);
            \foreach \angle in {0,...,\n}{
                \tikzmath{\start = 360*\angle/\n+\eps;
                \end = 360*(\angle+1)/\n-\eps;}
                \draw[thick, black, domain= \start: \end] plot({cos(\x)}, {sin(\x)});}
                \draw[thick, black](0,0) circle (3cm);
            \draw[thick, blue] (0,0) --  ({1/sqrt(2)},{1/sqrt(2)})  node[midway, above, sloped]{$R_1$};
            \draw[thick, blue] (0,0) --  (3,0)  node[midway, above, sloped]{$R_2$};
            \draw[thick,blue,decorate,decoration={brace,amplitude=12pt}] ({cos(360*(\n-1)/\n+\eps)},{sin(360*(\n-1)/\n+\eps)}) -- ({cos(360*(\n-1)/\n-\eps)},{sin(360*(\n-1)/\n-\eps)}) node[midway, below,yshift=-12pt, xshift= 6pt]{$\epsilon$};
        \end{tikzpicture}}
        \caption{Domain $D_{N, \epsilon}$ used in \cite{Hoffmann_NodalLine}, with $R_1 = 1$, $R_2 = 3$, $N=\n$ and $\epsilon = \eps$\textdegree}
        \label{DomainHoffmann}
    \end{subfigure}
    \hspace{0.05\textwidth}
    \begin{subfigure}[t]{0.28\linewidth}
        \tikzmath{\n = 5;
        \eps = 180/(\n+10);}
        \resizebox{\linewidth}{!}{
        \begin{tikzpicture}
            \definecolor{pastelgruen}{HTML}{acf6c8}
            \fill[pastelgruen] (0,0) circle (3cm);
            \foreach \angle in {0,...,\n}{
                \tikzmath{\start = 360*\angle/\n+\eps;
                \end = 360*(\angle+1)/\n-\eps;}
                \draw[thick, black, domain= \start: \end, fill=white] plot({cos(\start+\end-\x)}, {sin(\start+\end-\x)}) -- plot({2*cos(\x)}, {2*sin(\x)}) -- ({cos(\end)}, {sin(\end)}) ;}
                \draw[thick, black](0,0) circle (3cm);
            \draw[thick, blue] (0,0) --  ({cos(180/\n)},{sin(180/\n)})  node[midway, above, sloped]{$R_1$};
            \draw[thick, blue] (0,0) --  ({2*cos(180*(2*\n-1)/\n)},{2*sin(180*(2*\n-1)/\n)})  node[midway, above, xshift= 12pt, sloped]{$R_2$};
            \draw[thick, blue] (0,0) --  (3,0)  node[midway, above, xshift=25pt, sloped]{$R_3$};
            \draw[thick,blue,decorate,decoration={brace,amplitude=12pt}] ({cos(360*(\n-1)/\n+\eps)},{sin(360*(\n-1)/\n+\eps)}) -- ({cos(360*(\n-1)/\n-\eps)},{sin(360*(\n-1)/\n-\eps)}) node[midway, below,yshift=-12pt, xshift= 6pt]{$\epsilon$};
        \end{tikzpicture}}
        \caption{Domain $D_{N, \epsilon}$ used in \cite{Hoffmann_NodalLine} and \cite{Kennedy2011NodalLineRobin}, with $R_1 = 1$, $R_2 = 2$, $R_3 = 3$, $N=\n$ and $\epsilon = \eps$\textdegree}
        \label{DomainKennedy}
    \end{subfigure}
    \hspace{0.05\textwidth}
    \begin{subfigure}[t]{0.28\linewidth}
        \tikzmath{\n = 40;
        \eps = 180/(\n+10);}
        \resizebox{\linewidth}{!}{
        \begin{tikzpicture}
            \definecolor{pastelgruen}{HTML}{acf6c8}
            \fill[pastelgruen] (0,0) circle (3cm);
            \foreach \angle in {0,...,\n}{
                \tikzmath{\start = 360*\angle/\n+\eps;
                \end = 360*(\angle+1)/\n-\eps;}
                \draw[thick, black, domain= \start: \end, fill=white] plot({cos(\start+\end-\x)}, {sin(\start+\end-\x)}) -- plot({2*cos(\x)}, {2*sin(\x)}) -- ({cos(\end)}, {sin(\end)}) ;}
                \draw[thick, black](0,0) circle (3cm);
        \end{tikzpicture}}
        \caption{Domain $D_{N, \epsilon}$ used in \cite{Hoffmann_NodalLine} and \cite{Kennedy2011NodalLineRobin}, with $R_1 = 1$, $R_2 = 2$, $R_3 = 3$, $N=\n$ and $\epsilon = 3.6$\textdegree}
        \label{DomainKennedyToInfinity}
    \end{subfigure}
    \caption{Domains used to reject Paynes conjecture theoretically}
    \label{DomainRejectTheoretically}
\end{figure}
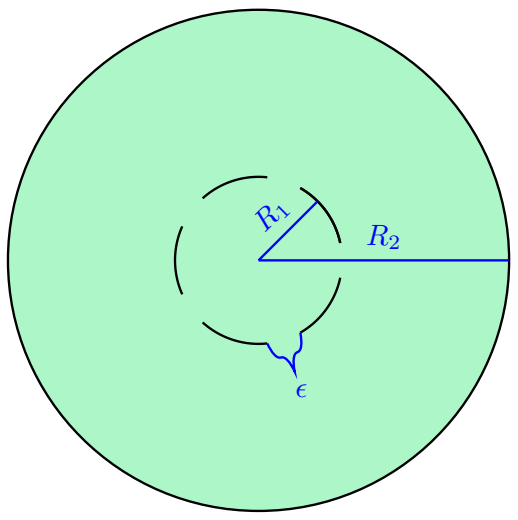
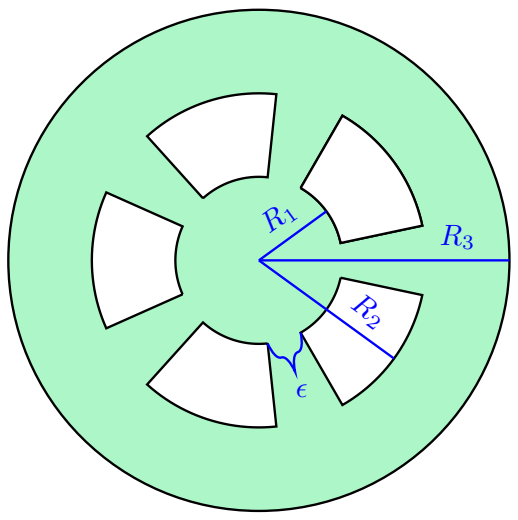
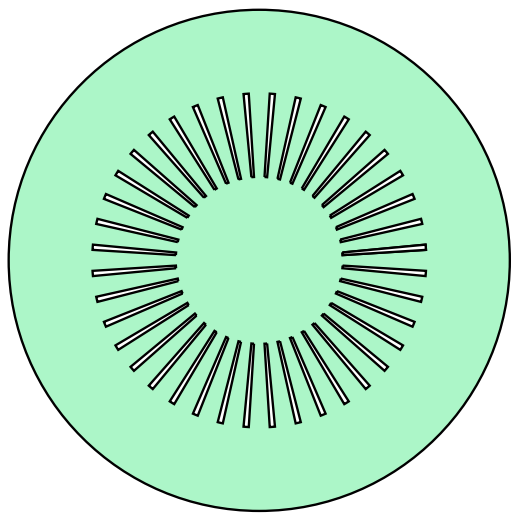

In 1995, Hoffmann-Ostenhof, Hoffmann-Ostenhof and Nadirashvili found a counterexample with a non-simply connected, bounded domain in $\mathbb{R}^2$ for Payne's conjecture. Therefore, they first defined two balls $B_{R_1}$ and $B_{R_2}$ with $B_{R_1} \coloneq \{x \in \mathbb{R}^2 : \|x\| < R_i \}$ and $0 < R_1 < R_2$. Based on these, the annulus $M_{R_1,R_2} \coloneq B_{R_2}\setminus B_{R_1}$ is defined. $R_1$ and $R_2$ have to be defined in such a way that
\begin{equation}
    \lambda_1(B_{R_1}) < \lambda_1(M_{R_1, R_2}) < \lambda_2(B_{R_1}).
\end{equation}
Carving holes into the boundary between the two areas $\partial B_{R_1} \cup \partial M_{R_1, R_2}$ results in the domain  
\begin{equation}
    D_{N,\epsilon} = B_{R_1} \cup M_{R_1, R_2} \cup \bigcup\limits_{J=0}^{N-1} \left\{x  \in \mathbb{R}^2 : r = R_1, \omega \in \left( \frac{2\pi j}{N}-\epsilon, \frac{2\pi j}{N}+\epsilon\right)\right\},
\end{equation}
with $N \in \mathbb{N}$ and $\epsilon <\frac{\pi}{N}$. This is shown in Figure \ref{DomainRejectTheoretically}(\subref{DomainHoffmann}).

With this domain it could be shown theoretically that there exists an $N_0 \in \mathbb{N}$ such that for every $N \geq N_0$ the second eigenfunction $u_2$ has a closed nodal line in $D_{N,\epsilon}$ and because of that never touches the boundary. It could not be found an explicit bound $N_0$. Furthermore, the same result could be shown for the domain displayed in Figure \ref{DomainRejectTheoretically}(\subref{DomainKennedy}) which was as well used by Kennedy in \cite{Kennedy2011NodalLineRobin}.
For larger $N$, it looks, for example, like the one shown in Figure \ref{DomainRejectTheoretically}(\subref{DomainKennedyToInfinity}) \cite{Hoffmann_NodalLine}.

In \cite{Hoffmann_NodalLine} the counterexample is constructed in $\mathbb{R}^2$. The symmetry properties granted in \cite{Hoffmann_NodalLine}, can generally not be used in higher dimensions. However, in \cite{Fournais2001NodalSurfaceClosed} they used a higher dimension generalization of the domain presented in Figure \ref{DomainRejectTheoretically}(\subref{DomainKennedy}). With this, Fournais could show that if there are enough, small cuts, which are equally distributed over the sphere, then the nodal surface can be closed. This is essentially the same statement, except that an upper limit for the minimum quantity of extinguishers can be specified. However, this limit is very large at $10^9$ \cite{Fournais2001NodalSurfaceClosed}.

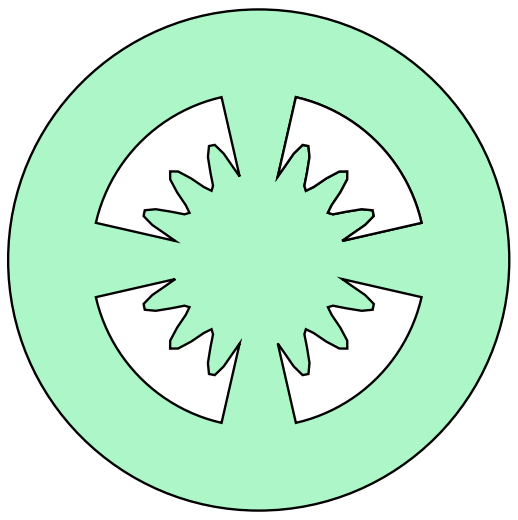
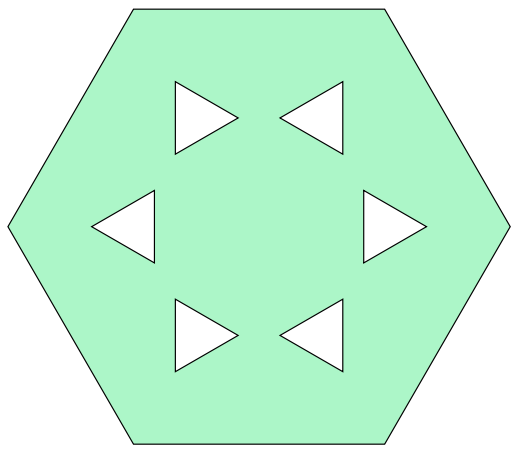
\begin{figure}[h!]
    \centering
    \begin{subfigure}[t]{0.28\textwidth}
        \tikzmath{\n = 4;
        \eps = 180/(\n+10);}
        \resizebox{\linewidth}{!}{
        \begin{tikzpicture}
            \definecolor{pastelgruen}{HTML}{acf6c8}
            \fill[pastelgruen] (0,0) circle (3cm);
            \foreach \angle in {0,...,\n}{
                \tikzmath{\start = 360*\angle/\n+\eps;
                \end = 360*(\angle+1)/\n-\eps;}
                \draw[thick, black, domain= \start: \end, fill=white] plot({(1+(1+cos(\n*\n*(\x)))/(\n))*cos((\x))}, {(1+(1+cos(\n*\n*(\x)))/(\n))*sin((\x)}) -- plot({2*cos(\start+\end-\x)}, {2*sin(\start+\end-\x)}) -- ({(1+(1+cos(\n*\n*(\start)))/(\n))*cos((\start))}, {(1+(1+cos(\n*\n*(\start)))/(\n))*sin((\start)});}
                \draw[thick, black](0,0) circle (3cm);
        \end{tikzpicture}}
        \caption{Domain $D_{N, \epsilon}$ used in \cite{Kennedy2011NodalLineRobin}, with $R_1 = 1$, $R_2 = 2$, $R_3 = 3$, $N=\n$ and $\epsilon = \eps$\textdegree}
        \label{DomainKennedyRobin}
    \end{subfigure}
    \hspace{0.05\textwidth}
    \begin{subfigure}[t]{0.28\textwidth}
        \resizebox{\linewidth}{!}{
        \begin{tikzpicture}
            \definecolor{pastelgruen}{HTML}{acf6c8}
            \node[regular polygon, draw, regular polygon sides = 6,fill = pastelgruen,
    	minimum size = 6cm] (p) at (0,0){};
            \foreach \angle in {0,60,...,300}
            \node[regular polygon, draw, regular polygon sides = 3,fill = white,
    	minimum size = 1cm, rotate= \angle+90/3] (\angle) at (\angle:1.5) {};
        \end{tikzpicture}}
        \caption{Domain $D_{N, \epsilon}$ used in \cite{Dahne_Counterexample}}
        \label{DomainDahne}
    \end{subfigure}
    \caption{Domains used to reject Paynes conjecture}
    \label{DomainRejectTheoreticallyOrNumerically}
\end{figure}

In comparison to our problem, for the Neumann-problem it is well known that the nodal line always touches the boundary. That is the reason why Payne's conjecture occurred. In the paper \cite{Kennedy2011NodalLineRobin}, the Robin problem is considered. It could be shown that, if the Dirichlet part of the boundary conditions is sufficiently small, the eigenfunction behaves like a Neumann eigenfunction in that it always touches the boundary (Proposition 1.3). If the Dirichlet term is bigger, the existence of a nodal domain which is compactly contained in $\Omega$, $\Omega^- \coloneq \{x \in  \Omega: u_2 (x)\leq 0\} \Subset \Omega$, where $U \Subset V$ means $U \subset K \subset V$ for some compact $ K \subset \mathbb{R}^n$, can been shown. This happened nearly in the same way as in \cite{Hoffmann_NodalLine} and \cite{Fournais2001NodalSurfaceClosed}. The first difference is that instead of using the ball $B_{R_1}$ and the annulus $M_{R_1,R_2}$ the annulus $A_{R_2,R_3}$ is used, with $0 < R_1 < R_2 < R_3$. This is also possible in \cite{Hoffmann_NodalLine}, as noted in Remark 1. The other difference is that in the case of the Robin-Problem, instead of the $B_{R_1}$, a sequence of perturbed domains $U_n$ is defined with four different properties. One of them is that $U_n$ has to be symmetric with respect to rotations through angles $\theta= \frac{2k\pi}{n}$, with $k=0,1,\dots, n-1$. But also in this case, holes are carved between $U_n$ and $A_{R_2,R_3}$. An example for one domain can be seen in Figure \ref{DomainRejectTheoreticallyOrNumerically}(\subref{DomainKennedyRobin}).

Paper \cite{Hoffmann_NodalLine} already posed the question of how small the $N_0$ actually has to be for there to be a region with $N_0$ boundary components such that the nodal line does not touch the boundary. Dahne tried to approximate this number in 2021 using the finite element method and found a counterexample with six holes, i.e. $N_0 = 7$ see Figure \ref{DomainRejectTheoreticallyOrNumerically}(\subref{DomainDahne}). However, it was assumed that it is also possible to find counterexamples with fewer holes.

Recently, a preprint \cite{Freitas2025_Paynesnodallineconjecture} occurred by Freitas and Leylekina and tried to answer this question analytically \footnote{The results presented in this paper were obtained independently. The results of Section \ref{sec:6Ellipses} and \ref{sec:5Ellipses} were first presented at the International High Performance Computing Summer School in Lisbon in July 2025 and were generated much earlier for PoF IV (Programmorientierte Förderung) in March 2025.} They claim that they have identified a doubly-connected domain rejecting Payne's conjecture. The main idea of this article is to define a family of domains with a sliding handle, like a thin annulus intersecting a convex domain, such as a rectangle with rounded corners, intersecting the annulus at one side. The authors assert that the nodal line would only touch the inner boundary of the domain if more vertices of the convex domain were on the side of the inner boundary, and vice versa. Accordingly, a sliding procedure has been formulated where the convex area is moved from the left to the right side. With the continuity argument, they argue that at a certain point, the nodal line would nowhere touch the boundary.
However, there is no information about the domain's exact dimensions, which makes it more difficult to prove this result numerically.
Moreover, this article has not yet been reviewed and also not verified numerically.
Additionally, the identified doubly-connected domain rejecting Payne's conjecture was determined through the application of the continuity argument. This suggests that a counterexample exists, though it may be quite unstable. This means that small changes to the domain's exact information could result in Payne's conjecture being fulfilled. Therefore, it may not be possible to obtain these results in practical applications, such as numerical calculations or the production of manufactured parts.
\subsection{Contributions}
\label{sec:Contribution}
The aim of this article is to determine the minimum number of holes that can be verified numerically and used to produce manufactured parts. Additionally, we want to investigate how the parameter areas that reject Payne's conjecture change for a specific domain's shape as the number of holes decreases. Furthermore, we want to figure out how the nodal line looks like, if a counterexample to Payne's can be found and how this shape depends on the domain's parameters.

Motivated by this, we developed an algorithm that can find such counterexamples by itself. Furthermore, we were able to build a high-performance implementation of this method to reduce the runtime significantly. Based on this, we were successful in numerically finding counterexamples with six and five holes. While previous work has focused on reducing this number of holes, the behavior of the nodal line has remained unexplored, a gap this paper addresses. We have evaluated many different nodal lines and have formulated the hypothesis, that the nodal line has always the shape of an hypotrochoid or an affine transformation of an hypotrochoid. Additionally, many other hypotheses occurred, such as: The smaller the number of holes without changing the parametrization, the smaller is the volume of the parameter set for which Payne's conjecture is rejected. Moreover, we could define quartic models which describe the behavior of the outer radius, inner radius and the ratio between both of them of the hypotrochoid depending on the parameters of the domain's shape. With this we could identify parameter areas, where it is impossible to find counterexamples to Payne's conjecture. This helped to argue why it is impossible to identify a counterexample with four holes in our domain shape when the parameters are gridded with a $0.01$ step width. All of the results are calculated using the strongly parallelized code using MPI and ran on the supercomputer JUSUF (Jülich Support for Fenix) \cite{Vieth2021_JUSUF} at the Jülich Supercomputing Centre.
\subsection{Outline}
\label{sec:Outline}
The remainder of this article is organized as follows:
In Section \ref{sec:Methodology} the methods used to find the counterexamples are described. In the first three Subsections (\ref{sec:BIE} -- \ref{sec:CollocationMethod}) the focus is on the numerical methods used to solve the partial differential equation and the eigenvalue problem in Section \ref{sec:Beyn}. Afterwards, the algorithm to find the counterexamples is developed in Subsection \ref{sec:FindParameters} and an efficient implementation is elaborated, see Subsection \ref{sec:Parallel}.
In Section \ref{sec:NumericalResults}, the numerical results generated by the previous explained methods are presented. The Subsection \ref{sec:ReviewDahne} reviews Dahne's findings, while the following two Subsections (\ref{sec:6Ellipses} -- \ref{sec:5Ellipses})
demonstrate our counterexamples with six and five holes. Then, we analyze the parameter sets for which Payne's conjecture is rejected in Subsection \ref{sec:ParametersetsRejectingPayne}. Furthermore, Subsection \ref{sec:BehaviorRatioNodalLine} analyzes how the shape of the nodal line changes with respect to the domain parameters such as the outer and inner radii of the hypotrochoid. In the last Subsection \ref{sec:4_holes}, we investigate whether it is possible to find numerically a counterexample with 4 holes based on our model with a 0.01 grid.
A summary and future work concludes this article, see Section \ref{sec:Conclusion}.

\section{Methodology}
\label{sec:Methodology}
\subsection{Boundary integral equations}
\label{sec:BIE}
As described in Section \ref{sec:RelatedWork}, the Dirichlet eigenvalues of the negative Laplace operator and its corresponding eigenfunctions have to be computed, as described in Equation \eqref{LaplaceEigenvalueProblem}.
The eigenvalue problem can be defined by
\begin{definition}[Dirichlet eigenvalue problem of the negative Laplacian]
    Find the non trivial eigenfunction $u \in C^2\left(\Omega\right) \cap C^1\left(\overline{\Omega}\right)$ for which
    \begin{align}
        \Delta u + \kappa^2 u &= 0 \quad \text{in } \Omega\\
        \text{and } \quad u&=0 \quad \text{on } \Gamma
        \label{EigenvalueProblem}
    \end{align}
    holds, with $\kappa^2$ the Dirichlet eigenvalue \cite[p. 8]{Chavel1984_Eigenvalues}. 
    \label{MainProblem}
\end{definition}
In this case, $\kappa^2$ was used instead of the eigenvalue $\lambda$. This is to show the analogy to the Helmholtz equation.
Whereas in the other papers, mostly the finite element method was used, see Dahne \cite[p. 4]{Dahne_Counterexample}, in this paper, the boundary element method (BEM) is used, which is based on boundary integral equations (BIEs). 
Especially when the considered domains become more complicated, it is difficult to find a suitable triangulation with the finite element method. This problem does not appear when using the BEM because the BIEs reduce the problem from two dimensions to one.
It is known from the reduced wave equation that the wave number is positive and constant $\kappa>0$, see \cite[sec. 2.1]{Colton1998_Inverse}.
For the Helmholtz equation in $\R^2$, the fundamental solution is well known:
\begin{definition}[Fundamental solution]
    The fundamental solution of the Helm\-holtz equation in $\R^2$ is given by
    \begin{equation}
        \Phi_{\kappa}(x,y) = \frac{\mi H_0^{(1)}(\kappa\|x-y\|)}{4}, \qquad x\neq y,
        \label{FundamentalSolution}
    \end{equation}
    where $x,y \in \R^2$ and $H_0^{(1)}$ denotes the first-kind Hankel function of order zero.
\end{definition}
For this reason, this problem can be solved using the BIE. Here, we choose the double layer ansatz to get a BIE of the second kind. 

The double layer approach is defined by
\begin{equation}
     u(x) = \int\limits_{\Gamma} \Psi(y)\frac{\partial\Phi_{\kappa}(x,y)}{\partial \nu(y)}\; \md s(y), \qquad x \in \Omega,
    \label{DoubleLayerAnsatzDomain}
\end{equation}
with the exterior normal $\nu(y)$ \cite[p. 39]{Colton1998_Inverse}.
To solve this, the integrable density function $\Psi$ must first be calculated on the boundary $\Gamma$ using jump relations
\begin{equation}
    -\frac{1}{2}\Psi(x)+ \int\limits_{\Gamma} \Psi(y)\frac{\partial\Phi_{\kappa}(x,y)}{\partial \nu(y)} \; \md s(y) = 0, \qquad x \in \Gamma,
    \label{DoubleLayerAnsatzBoundary}
\end{equation}
refer to \cite[p. 39]{Colton1998_Inverse}.
This is a boundary integral equation of the second kind.
This approach is used in the following, because unlike the single layer approach, the problem is well posed.

Next, as in \cite[sec. 4.2]{Kleefeld_2021HotSpotsConjecture}, we convert the problem into a form that takes the presence of the holes more into account, and thus simplifies the calculation by eliminating terms.

To do this, the integral of the double layer operator in Equation \eqref{DoubleLayerAnsatzDomain} is first split and all functions $u(x)$ for which $x \in \Gamma_i$ applies are defined by $u_i(x)$, just as all exterior normal $\nu(x)$ are defined as $\nu_i(x)$ and all $\Psi(x)$ as $\Psi_i(x)$, where $i=0,\dots q$ and $q$ represents the number of holes within the area $\Omega$. This yields
\begin{align}
    u(x) = \sum\limits_{i=0}^q\, \int\limits_{\Gamma_i} \Psi_i(y)\partial_{\nu_i(y)} \Phi_{\kappa}(x,y) \;\md s(y), \qquad x \in \Omega.
    \label{SplittingIntegralsDomain}
\end{align}

Just as within the domain, the integrals of Equation \eqref{DoubleLayerAnsatzBoundary} are also divided, to obtain the unknown density functions $\Psi_i(x)$ for $i=0,\dots, q$
\begin{align}
    0 =\sum\limits_{i=0}^q -\frac{1}{2}\Psi_i(x)+ \int\limits_{\Gamma_i} \Psi_i(y)\partial_{\nu_i(y)} \Phi_{\kappa}(x,y) \;\md s(y), \qquad x \in \Gamma.
    \label{SplittingIntegralsBoundary}
\end{align}
Afterwards, the area in which $x$ lies is partitioned. For $x_i \in \Gamma_i$, where $i=0,\dots,q$, the terms $\Psi_k(x_i)$, where $k \neq i$ and $k=0,\dots,q$, are omitted, since $x_i \notin \Gamma_k$, the following applies
\begin{align}
    0 =-\frac{1}{2}\Psi_k(x) + \sum\limits_{i=0}^q \int\limits_{\Gamma_i} \Psi_i(y)\partial_{\nu_i(y)} \Phi_{\kappa}(x,y) \;\md s(y), \qquad x \in \Gamma_k.
\end{align}
By defining
\begin{equation}
    \mathrm{D}_{\kappa}^{\Gamma_i \to \Gamma_k} \left[\Psi\right](x)\coloneq \int\limits_{\Gamma_i} \Psi_i(y)\partial_{\nu_i(y)} \Phi_{\kappa}(x,y)\; \md s(y), \qquad x \in \Gamma_k, \qquad i,k = 0,\dots, q,
    \label{DoubleLayerBoundaryFromTo}
\end{equation}
the problem on the boundary can be described by 
\begin{align}
    M_{\kappa} \left[\Psi\right](x) =
    \begin{pmatrix}
        -\frac{1}{2}I+
        \begin{pmatrix}
            \mathrm{D}_{\kappa}^{\Gamma_0 \to \Gamma_0} & \cdots & \mathrm{D}_{\kappa}^{\Gamma_q \to \Gamma_0}\\
            \vdots & \ddots & \vdots\\         
            \mathrm{D}_{\kappa}^{\Gamma_0 \to \Gamma_q} & \cdots & \mathrm{D}_{\kappa}^{\Gamma_q \to \Gamma_q}\\
        \end{pmatrix}
    \end{pmatrix}
    \left[
    \begin{pmatrix}
        \Psi_0\\ \vdots \\ \Psi_q
    \end{pmatrix}\right](x)
    = \begin{pmatrix}
        0\\ \vdots \\ 0
    \end{pmatrix},
    \label{BIESystem}
\end{align}
where $I$ is the identity operator. In the following we use $M\left(\kappa\right) = M_{\kappa}$.
\subsection{Boundary element collocation method}
\label{sec:CollocationMethod}
\begin{figure}[h!]
    \centering
    \begin{tikzpicture}[scale=0.7, transform shape]
        \definecolor{pastelgruen}{HTML}{acf6c8}
        \fill[pastelgruen] (0,0) circle (3cm);
        \tikzmath{\n = 5;
        \nminus= \n-1;}
        \draw[thick, black, fill=white] (0,0) circle (1.5cm);
        \draw[thick, black](0,0) circle (3cm);
        \foreach \angle [count=\xi] in {0,...,\nminus}{
            \tikzmath{\newangle = 360*\angle/\n;
            \nextangle = \newangle+360/\n;
            \meanangle = (\newangle+\nextangle)/2;}
            \draw[line width=1pt] (\newangle:1.4cm) -- (\newangle:1.6cm);
            \node[font=\large] at (-\meanangle:1.cm) {\textsf{$\Delta_{1, \xi}$}};
            \draw[line width=1pt] (\newangle:2.9cm) -- (\newangle:3.1cm);
            \node[font=\large] at (\meanangle:3.5cm) {\textsf{$\Delta_{0, \xi}$}};
            }
    \end{tikzpicture}
    \caption{Subdivision $\Delta_{i,j}$ of the boundary of an area with one hole. Each boundary is separated in $n_f = 5$ pieces, which means $j=1,\dots, 5$ and $i=0,1$, where the outer boundary has always the index $i=0$.}
    \label{Subdivision}
\end{figure}
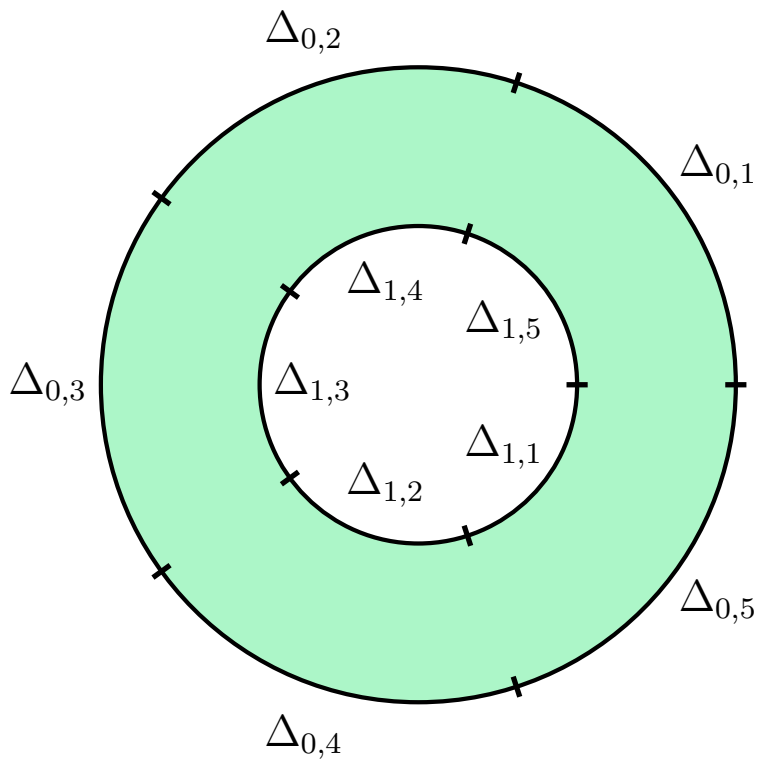

In the next step, this system of BIEs \eqref{BIESystem} is to be solved using the boundary element collocation method as in \cite[sec. 4.3]{Kleefeld_2021HotSpotsConjecture}. To do this, the given domain is first divided into $(q+1) n_f$ subdomains, which are denoted by $\Delta_{i,j}$. Here $q$ stands for the number of holes within the area $\Omega$, $i=0,\dots, q$ and $j=1,\dots, n_f$, where $n_f$ is the amount of pieces in which one boundary is separated. A simple example with just one hole for such a subdivision can be seen in Figure \ref{Subdivision}.
It is important to note that the outer boundary $\Gamma_0$ is always oriented counterclockwise, whereas the inner boundary of the holes $\Gamma_i$, where $i=1,\dots , q$, must be oriented clockwise. Due to the different orientation of the boundaries, the outer normals must be selected with different directions. While the normal points outwards at the outer boundary, it points inwards at the holes. 
This allows each line of the equation \eqref{BIESystem} to be converted to
\begin{align}
    \sum\limits_{i=0}^{q} \sum\limits_{j=1}^{n_f} \int\limits_{\Delta_{i,j}} \Psi_i(y)\partial_{\nu_i(y)}
        \Phi_{\kappa}(x,y)\; \md s(y) -\frac{1}{2}\Psi_k(x) = 0, \qquad x \in \Gamma_k,
    \label{SubdivisionAreas}
\end{align}
where $\Psi_i(y)$ is the unknown density function at the $i$-th boundary $y \in \Gamma_i$.

In the next step, for each subdomain a bijective mapping $m_{i,j}$ can be used, by applying a change of variables. This maps from the standard interval $\sigma =[0,1]$ to $\Delta_{i,j}$ and results in
\begin{equation}
    \begin{aligned}
        &\sum\limits_{i=0}^{q} \sum\limits_{j=1}^{n_f} \int\limits_{\sigma} \Psi_i(m_{i,j}(s))\partial_{\nu_i(m_{i,j}(s))}
        \Phi_{\kappa}(x,m_{i,j}(s)) J_{i,j}(s) \; \md s(s)\\ 
        &-\frac{1}{2}\Psi_k(x) = 0, \qquad x \in \Gamma_k,
        \label{BIEChangeOfVariables}
    \end{aligned}
\end{equation}
with $J_{i,j}(s) = \Vert\partial_s m_{i,j}(s)\Vert$, where $\Vert \cdot \Vert$ is the Euclidean norm.

Since the mapping $m_{i,j}$ cannot be written down explicitly in all cases, a quadratic Lagrange interpolation is carried out for each $m_{i,j}(s)$ denoted by $\widetilde{m}_{i,j}(s)$
\begin{align}
    m_{i,j}(s) \approx \widetilde{m}_{i,j}(s) = \sum\limits_{p=1}^3  v_{i,((p-1) + 2 (j-1))\bmod{(2n_f)}+1} L_p (s),
\end{align}
where $v_{i,h}$ is the $h$-th point on the boundary, containing the vertices and midpoints at the $i$-th boundary.
Afterwards, the unknown density function has to be approximated with a quadratic interpolation, using the generalized Lagrange basis functions $\widetilde{L}_r(s)$, with $r = 1,2,3$. For a detailed explanation, see \cite[p. 8]{Kleefeld_2021HotSpotsConjecture}.
Therefore, we need the definition of the collocation nodes
\begin{align}
    \widetilde{v}_{k,\ell,r} = \widetilde{m}_{k,\ell}(c_r), \qquad \text{with} \qquad r=1,2,3
\end{align}
with 
\begin{align}
    c_1= \alpha, \qquad c_2 = \frac{1}{2}, \qquad \text{and} \qquad c_3 = 1-\alpha, \qquad\text{where} \qquad 0<\alpha< \frac{1}{2}.
\end{align}
The selected $\alpha$ is a constant value. The open interval of this constant ensure that the collocation nodes are always lying within the piece of the boundary. Furthermore, depending on how $\alpha$ is selected, the convergence rate can be improved. 

By inserting the approximations in \eqref{BIEChangeOfVariables} and factoring out the density function from the integral, which is independent of the integration variable, each line of \eqref{BIESystem} can be approximated.

The residual of this approximation is set to zero at the collocation points. Here, $\widetilde{v}_{k,\ell,r}$ describes the $r$-th collocation node at the $\ell$-th boundary, with an integration over the $k$-th boundary.
This is how a linear system of size $3n_f(q+1) \times 3n_f(q+1)$ occurs
\begin{equation}
    \begin{aligned}
        &\sum\limits_{i=0}^{q} \sum\limits_{j=1}^{n_f} \sum\limits_{p=1}^3 \int\limits_{\sigma}  \widetilde{L}_p(s) \partial_{\nu_i(\widetilde{m}_{i,j}(s))}
        \Phi_{\kappa}(\widetilde{v}_{k,\ell,r},\widetilde{m}_{i,j}(s)) \Vert\partial_s \widetilde{m}_{i,j}(s)\Vert \;\md s(s) \Psi_i\left(\widetilde{v}_{i,j,p}\right)\\
        & -\frac{1}{2}\Psi_k(\widetilde{v}_{k,\ell,r}) = 0, \\
    \end{aligned}
\end{equation}
which can be simplified to
\begin{align}
    \sum\limits_{i=0}^{q} \sum\limits_{j=1}^{n_f} \sum\limits_{p=1}^3 a_{k,\ell,r,i,j,p} \Psi_i\left(\widetilde{v}_{i,j,p}\right) -\frac{1}{2}\Psi_k(\widetilde{v}_{k,\ell,r}) = 0,
\end{align}
with the integrals
\begin{align}
    a_{k,\ell,r,i,j,p} = \int\limits_{\sigma}  \widetilde{L}_p(s) \partial_{\nu_i(\widetilde{m}_{i,j}(s))}
    \Phi_{\kappa}(\widetilde{v}_{k,\ell,r},\widetilde{m}_{i,j}(s)) \Vert \widetilde{m}_{i,j}(s)\Vert \;\md s(s).
    \label{BIESystemIntegrals}
\end{align}
Inserting the fundamental solution of the Helmholtz equation and applying the directional derivative leads to the following form of the integrand of the integral \eqref{BIESystemIntegrals}, where the Jacobian cancels out
\begin{align}
    \widetilde{L}_p(s) \frac{\mi k H_1^{(1)}\left(k\sqrt{a^2+b^2}\right)}{4\sqrt{a^2+b^2}}(a\cdot n_1 + b\cdot n_2),
    \label{BIESystemIntegrandIntegrals}
\end{align}
with $ H_1^{(1)}$ the first-kind Hankel function of order one and 
\begin{equation}
    \begin{aligned}
    a &= \left[\widetilde{v}_{k,\ell,r} - \widetilde{m}_{i,j}(s)\right]_1, & b&= \left[\widetilde{v}_{k,\ell,r} - \widetilde{m}_{i,j}(s)\right]_2, \quad &n_1 & = \left[\widetilde{J}_{i,j}(s)\right]_2,\\
    n_2 & = -\left[\widetilde{J}_{i,j}(s)\right]_1, \quad & J_{i,j}(s) &= \partial_s \widetilde{m}_{i,j}(s).
\end{aligned}
\end{equation}
All integrals in Equation \eqref{BIESystem} can be computed numerically with the 21-point adaptive Gauss-Kronrod quadrature \cite{Piessens1983_Quadpack}. 
\subsection{Beyn algorithm}
\label{sec:Beyn}
To solve this non-linear eigenvalue problem in Equation \eqref{BIESystem} the Beyn algorithm can be used.
First, a random matrix $\hat{V} \in \C^{m \times l}$ is generated with $n(\gamma) \leq l \leq m$ such that $W^H\hat{V} \in \C^{n(\gamma), l}$ has the condition rank $n(\gamma)$, which implies that $W$ must have the rank $n(\gamma)$.
Furthermore, it is assumed that the condition rank of $V$ is $n(\gamma)$.

Beyn's algorithm now uses Keldysh's theorem from \cite[Theorem 2.4]{Beyn2012_Beyn} twice, with $f(z) = 1$ and $f(z) = z$, and inserting $\hat{V}$ to reduce the size of the matrices $A_0$ and $A_1$ given by
\begin{align}
    A_0 &= \frac{1}{2\pi\mi}\int\limits_{\gamma} M(z)^{-1}\hat{V} \;\md z = V W^H \hat{V} ,\label{A0}\\
    A_1 &= \frac{1}{2\pi\mi}\int\limits_{\gamma} z M(z)^{-1} \hat{V}\;\md z = V\Lambda W^H \hat{V},
    \label{A1}
\end{align}
where $\Lambda = \diag(\kappa_1,\dots,\kappa_{n(\gamma)})$, $A_0, A_1 \in \C^{m\times l}$ and $\rank(A_0) = n(\gamma)$.
Afterwards, a reduced singular value decomposition for $A_0$, with $V_0 \in \C^{m\times n(\gamma)}$, $\Sigma_0 \in \C^{n(\gamma) \times n(\gamma)}$ and $W_0 \in \C^{l \times n(\gamma)}$ is computed.
Together this gives
\begin{align}
    A_0 = V_0 \Sigma_0 W_0^H.
\end{align}
This equation can be converted into
\begin{equation}
    W^H \hat{V} = S^{-1} V_0^H V_0 \Sigma_0 W_0^H = S^{-1} \Sigma_0 W_0^H.
    \label{WhVhat}
\end{equation}
By using this Equation \eqref{WhVhat}, we can similarly calculate 
\begin{align}
    B \coloneq S \Lambda S^{-1} = V_0^H A_1 W_0 \Sigma_0^{-1}.
\end{align}

To carry out this algorithm, the two integrals from Equation \eqref{A0} and Equation \eqref{A1} must be calculated numerically first, as in \cite[sec. 3.3]{Beyn2012_Beyn}. To do so, first select a contour to be integrated over that has a $2\pi$-periodic smooth parameterization, such as the circle with the parameterization
\begin{align}
    \phi(t) = \mu + R \cos(t) + \mi R \sin(t).
\end{align}
Here $R$ is the radius of the circle and $\mu$ the midpoint.
First, the parameterization is applied to the integrals in Equations \eqref{A0} and \eqref{A1} and afterwards, the integrals are approximated using the trapezoidal rule at the equidistant nodes $t_p = \frac{2\pi p}{N}$, with $ p= 1,\dots, N$. This results in
\begin{align}
    A_{0,N} &= \frac{1}{\mi N}\sum\limits_{p=0}^{N-1} M(\phi(t_p))^{-1} \hat{V}\phi'(t_p)\label{A0Trapez}\\
    A_{1,N} &= \frac{1}{\mi N}\sum\limits_{p=0}^{N-1} \phi(t_p) M(\phi(t_p))^{-1} \hat{V}\phi'(t_p).
    \label{A1Trapez}
\end{align}
Based on these calculations, the truncated singular value decomposition of $A_{0,N}$ has to be done. Subsequently, the matrix $B$ must be calculated, and finally the eigenvalue problem of $B$ can be solved and the eigenvectors $V= V_0 S$, associated with the problem $M(k)= 0$, can be calculated as already described.
\subsection{Algorithm to find domain parameters}
\label{sec:FindParameters}
In this section, the algorithm to find parameters of the domains, for which the nodal line does not touch the boundary, will be described. The algorithm bases on a theoretical result of the paper \cite[Lemma 3.2]{Melas1992_NodalLine}. The second eigenvalue of the domain $\Omega$ is equal to the first eigenvalue of a simply connected subdomain $\Omega_1$, $\lambda_2(\Omega) = \lambda_1 (\Omega_1)$, if the boundary of the subdomain is contained in the nodal line $\mathcal{N}(u_2)$ and the boundary $\Gamma$.

With this theoretical knowledge, an algorithm can be developed to find a suitable parameter set for a given parameterization, which rejects Payne's conjecture. If a domain exists for which the nodal line does not touch the boundary, then it is possible to cut a hole in this domain, with the nodal line as the boundary. Additionally, the first eigenvalue of this new domain has to have the same eigenvalue as the second eigenvalue of the domain without hole.

For a better understanding, in the following, the algorithm which we used to find the numerical results in Section \ref{sec:NumericalResults} will be described with the help of the first found result from Section \ref{sec:6Ellipses}.
At the beginning the parametrizations of the two domains have to be defined.
On the one hand, we chose a circle as the outer boundary. On the other hand, we decided to choose six ellipses in the inner of the big circle. The reason for this is that having boundaries which are twice differentiable, is necessary for using the boundary element method. This results in the parameterization
\begin{equation}
\begin{aligned}
    r_0(t) &= 
    \begin{pmatrix}
      x_m +a_0\cos(t_j)\\
      y_m+b_0\sin(t_j)  
    \end{pmatrix},\\
    r_i(t) &= 
    \begin{pmatrix}
      x_m +a_i\cos(t_j)\cos(\alpha)-b_i\sin(t_j) \sin(\alpha)\\
      y_m + a_i\cos(t_j)\sin(\alpha)+  b_i\sin(t_j) \cos(\alpha)
    \end{pmatrix},
    \qquad i \in 1 \dots q.
\end{aligned}
\end{equation}
For the second domain an additional hole has to be set. The shape of this hole is chosen based on the results in \cite[Fig. 1(b)]{Dahne_Counterexample}. As can be seen in Figure \ref{fig:Dahne_all_results}(\subref{fig:Dahne_13_27}) of this paper, the shape of the nodal line is a geometric curve (roulette). In general, this is a curve, which is created by rolling one curve along another. In our case, it is a hypotrochoid, which will be generated by rolling one circle through another.
The parameterization of the roulette can be developed in the following way.
First, the parametric equation for a hypotrochoid is given by 
\begin{equation}
\begin{aligned}
    c &= (R-r)\cos(t_j)+ar\cos\left(\left(\frac{R}{r}-1\right)t_j\right),\\
    d &= (R-r)\sin(t_j)-ar\sin\left(\left(\frac{R}{r}-1\right)t_j\right),\\
    r_{q+1}(t) &= 
    \begin{pmatrix}
      x_m +c\cos(\alpha)-d \sin(\alpha)\\
      y_m + c\sin(\alpha)+ d \cos(\alpha)
    \end{pmatrix},
\end{aligned}
\label{parameterizationHypotrochoid}
\end{equation}
where $R$ is the radius of the fixed circle, $r$ is the radius of the rolling circle, which rolls inside the bigger circle, $a$ is the by the radius of the revolving circle normalized distance from the center of the revolving circle to the drawing point, $(x_m, y_m)$ is the center of the bigger circle and $\alpha$ is the rotation of the hypotrochoid in counterclockwise direction in radian.
The needed parameters to draw the hypotrochoid in the same way as in Figure \ref{fig:Dahne_all_results}(\subref{fig:Dahne_13_27}) can be calculated with three parameters: $x_o$, the radius of the circumcircle, $x_i$ the radius of the incircle and $A$ the number of loops, which can be read directly from the graph.
Then,
\begin{equation}
    \begin{aligned}
        r &= \frac{x_o + x_i}{2(A-1)},\\
        R &= Ar,\\
        a &= (A-1) \frac{x_o-x_i}{x_o+x_i},
    \end{aligned}
\end{equation}
can be calculated and inserted into Equation \eqref{parameterizationHypotrochoid}.
This results in two different domains.
\begin{figure}[h!]
    \centering
        \begin{subfigure}[b]{0.45\linewidth}
        \begin{tikzpicture}[scale=0.8, transform shape]
            \definecolor{pastelgruen}{HTML}{acf6c8}
            \node[circle, draw,fill = pastelgruen, minimum size = 6cm] (p) at (0,0){};
            \foreach \angle in {0,60,...,300}
                \draw[fill = white] (\angle:1.5) ellipse [x radius = 0.7cm, y radius = 0.15cm, rotate= \angle+90];
            \draw[thick, blue] (0,0) --  (1.5,0)  node[midway, above, sloped]{$r_E$};
            \draw[thick, blue] (-1.5,0) node[above, sloped]{$b$} --  (-1.35,0)  ;
            \draw[thick, blue] (60:1.5) --  ({1.5*cos(60)+0.7*cos(330)},{1.5*sin(60)+0.7*sin(330)})  node[midway, above=2, sloped]{$a$};
        \end{tikzpicture}
        \caption{Start domain to calculate the second eigenvalue}
        \label{DomainCircleEllipses}
        \end{subfigure}
        \hspace{0.05\textwidth}
        \begin{subfigure}[b]{0.4\linewidth}
        \begin{tikzpicture}[scale=0.8, transform shape]
            \tikzmath{\Ahyp=6; \xohyp = 0.7; \xihyp = 0.5;
            \rhyp = (\xohyp + \xihyp)/(2*(\Ahyp-1));
            \Rhyp = \Ahyp*\rhyp;
            \ahyp = (\Ahyp-1)*(\xohyp-\xihyp)/(\xohyp+\xihyp);
            \xmhyp=0; \ymhyp=0;\alphahyp=0;};
            \definecolor{pastelgruen}{HTML}{acf6c8}
            \node[circle, draw,fill = pastelgruen, minimum size = 6cm] (p) at (0,0){};
            \foreach \angle in {0,60,...,300}
                \draw[fill = white] (\angle:1.5) ellipse [x radius = 0.7cm, y radius = 0.15cm, rotate= \angle+90];
            \draw[fill = white, domain=0:360, samples=1000, variable=\t] plot({\xmhyp +((\Rhyp-\rhyp)*cos(\t)+\ahyp*\rhyp*cos((\Rhyp/\rhyp-1)*\t))*cos(\alphahyp)-((\Rhyp-\rhyp)*sin(\t)-\ahyp*\rhyp*sin((\Rhyp/\rhyp-1)*\t))* sin(\alphahyp)
          },{\ymhyp + ((\Rhyp-\rhyp)*cos(\t)+\ahyp*\rhyp*cos((\Rhyp/\rhyp-1)*\t))*sin(\alphahyp)+ ((\Rhyp-\rhyp)*sin(\t)-\ahyp*\rhyp*sin((\Rhyp/\rhyp-1)*\t))* cos(\alphahyp)});
            \draw[thick, blue] (0,0) --  (1.5,0)  node[pos = 0.7,right, above, sloped]{$r_E$};
            \draw[thick, blue] (-1.5,0) node[above, sloped]{$b$} --  (-1.35,0)  ;
            \draw[thick, blue] (60:1.5) --  ({1.5*cos(60)+0.7*cos(330)},{1.5*sin(60)+0.7*sin(330)})  node[midway, above=2, sloped]{$a$};
            \draw[thick, blue] (0,0) --  (0,-0.5)  node[midway,right,sloped,rotate=90, xshift=-3.]{$x_i$};
            \draw[thick, blue] (0,0) --  ({0.7*cos(120)},{0.7*sin(120)})  node[midway,right,sloped,rotate=90, xshift=-2., yshift=-3]{$x_o$};
        \end{tikzpicture}
        \caption{Start domain to calculate the first eigenvalue}
        \label{DomainCircleEllipsesHypotrochoid}
        \end{subfigure}
    \caption{}
    \label{DomainCircleEllipsesHypotrochoid2}
\end{figure} 
It is important to check that the boundaries of the holes do not touch each other and do not touch the outer boundary. For example, if the hypotrochoid touches the inner ellipses, we would adjust the parameter set such that the nodal line of the second eigenfunction resembles a hypotrochoid that touches the inner ellipses. This would be consistent with Payne's conjecture and would not disprove it.
If we increase the parameters of the inner ellipses too much, one domain will split into two, which is not the kind of problem we are looking for. 

Now, the second eigenvalue of the left domain in Figure \ref{DomainCircleEllipsesHypotrochoid2}\subref{DomainCircleEllipses} and the first eigenvalue of the right domain \ref{DomainCircleEllipsesHypotrochoid2}\subref{DomainCircleEllipsesHypotrochoid} have to be calculated. Therefore, the Beyn-Algorithm from Section \ref{sec:Beyn} is used. This algorithm calculates the eigenvalues of the problem defined in Equation \eqref{LaplaceEigenvalueProblem} within a circle, with the middle point laying on the positive real axis. The circle starts close to zero and moves along the positive real axis. This is due to the fact that the eigenvalues of this problem are real positive numbers. Afterwards, the obtained results must be checked to determine whether they are eigenvalues for our defined problem, eigenvalues of another problem or scattering poles, which lie in the negative half axis and are resonance states \cite[ch. 9]{Taylor2011_PDE2}. Filtering the scattering poles, can be done very easily. It is sufficient to delete all results, which have an imaginary part that is smaller than a given $eps = 10^{-5}$. Afterwards, it is necessary to check, whether the eigenvalues correspond to the eigenvalue problem defined in Equation \eqref{EigenvalueProblem}. In some cases, eigenvalues occur which correspond to another eigenvalue problem, where one or more boundaries are ignored. This happens exactly then, when the eigenfunction on those boundaries is approximately equal to zero. This has to be proven first. Moreover, we can use the knowledge of the first two eigenvalues to check, whether we found the correct ones. First of all, it is well known, that first eigenvalues are always simple and second eigenvalues can have maximal a double eigenvalue. With this knowledge, it is easy to delete all eigenvalues which have a bigger multiplicity than one, if no eigenvalue has been found yet. If the first eigenvalue could be detected and the second eigenvalue is searched, then all eigenvalues with a bigger multiplicity than two can be deleted. Next, to figure out whether the given outcome is still a spurious or a proper eigenvalue, it is necessary to calculate the eigenfunction of the current result of the Beyn-algorithm as defined in the Equations \eqref{SplittingIntegralsDomain} and \eqref{DoubleLayerAnsatzDomain}. For this eigenfunction, the number of nodal domains has to be calculated. If we are looking for the first eigenvalue, it is necessary that the number of nodal domains is equal to one; if we are looking for the second eigenvalue it has to be two. If the number of nodal domains is higher, then the outcome is an eigenvalue of another problem. If not all eigenvalues sought have been found, but all results of the Beyn algorithm have been checked, the contour of the algorithm must be shifted along the positive axis.

At the end, the two results have to be compared. As explained in the beginning of this section, the second eigenvalue of the left domain is equal to the first eigenvalue of the right domain, if the shape of the inner hypotrochoid of the seventh hole of the domain in Figure \ref{DomainCircleEllipsesHypotrochoid2}\subref{DomainCircleEllipsesHypotrochoid} has the same shape as the nodal line of the second eigenfunction of the domain in Figure \ref{DomainCircleEllipsesHypotrochoid2}(\subref{DomainCircleEllipses}).
To compare the results, we use the SciPy function \texttt{scipy.optimize.least\_squares}. With this function, the parameters $a$, $b$, $r_E$, $x_o$  and $x_i$, shown in Figure \ref{DomainCircleEllipsesHypotrochoid2}\subref{DomainCircleEllipsesHypotrochoid}, will be adjusted. Due to the reason that changes of the the parameters of the hypotrochoid have an impact on the parameters of the ellipses, it is important to stop the algorithm after a few iterations, set the parameters of the hypotrochoid fix, and afterwards, find the best parameters $a$, $b$ and $r_E$ of the ellipses. 
\subsection{Parallel processing of the method}
\label{sec:Parallel}
In the previous chapters, we explained algorithms to solve the Helmholtz equation, a nonlinear eigenvalue problem and to find a parameter set to reject Payne's conjecture. All of these algorithms were described theoretically with no attention to the implementation of the algorithms. 
Since, inserting more holes into an area enlarges the matrix, obtained by discretizing the double layer operator in Equation \eqref{BIESystem}, and increasing the complexity of the area requires more collocation points to describe it accurately, it is necessary to consider the most efficient implementation method. Our simulations are written in the programming language Python and run on the the supercomputer JUSUF at the Jülich Supercomputing Centre \cite{Vieth2021_JUSUF}. It consists out of $144$ standard compute nodes, where each node has $2\times 64$ cores.  We were able to use several of these nodes to run our calculations, mostly six. This results in $768$ cores, on which independently can be calculated. To use these cores in the most efficient way, we want to introduce a high-performance way to parallelize the listed algorithms. 

The main algorithm is explained in Chapter \ref{sec:FindParameters} and describes how to find the parameter set with which Payne's conjecture can be rejected. As explained in the previous chapter, in each iteration step, two different eigenvalues from two different domains have to be calculated. We decided to calculate each of the domains one by one. The reason for this is that both calculations take a different amount of time. Parallelizing both calculations would end in waiting for the other task, where several cores are doing nothing and are wasting core-hours. It is not possible to estimate a factor that indicates how much faster one task is in relation to another. Because of this, we are using a serial approach here. Based on the residuum of the first and the second eigenvalue of the two domains, the least square algorithm changes the parameter set and the calculation of the two eigenvalues starts again. The causality between the individual steps necessitates again a serial approach. 

In that main algorithm, we have to calculate many times the first and the second eigenvalue of two different domains. This is done in a subroutine, based on Beyn's algorithm, see Section \ref{sec:Beyn}. To find the first and the second eigenvalue, we have to place circles with their midpoint along the positive real axis. Within these circles all eigenvalues lying in this area are calculated. By ordering all eigenvalues, one can find the first and the second eigenvalue of a specific domain. It is possible to parallelize the computations for the different circles. We have decided against it. The reason for this is, that the eigenvalues of interest often lie within the first few circles. If we parallelize it, we have to choose at the beginning a specific amount of circles $n$ and split all the available cores by the amount of circles. If the searched eigenvalue lies within the $n$-th circle, there is no problem, and this substep would be good parallelized. The problem is that one can not know exactly where the searched eigenvalue occurs. So it is possible that the searched eigenvalue lies within the circle $n-1$ or even before. Then, we would have a reduced amount of cores to do the calculations for each circle in comparison, if we would know the exact number. This means a runtime extension, which is not necessary. By keeping this part serial, we avoid this problem. Furthermore, the next steps can be parallelize very efficiently and can be reduced to very small subcomputations. This means that the difference between the serial computation of the different circles and the parallelization of the circles by knowing the exact number of needed circles $n$ would be negligible. We chose our circles with a radius $R=0.5$ and used $N=24$ points at the boundary of each circle. After computing the eigenvalues of the first circle, we shifted the circle by $0.9$ along the positive real axis. The Beyn-algorithm calculates eigenvalues correctly, which lie in the inner part of the circle. By shifting the circle a little bit less than twice the radius, we avoid that eigenvalues lie on the boundary of the circle and cannot be recognized because of that.

As mentioned in the previous paragraph, the eigenvalues within one circle have to be calculated. Therefore, we take $N$ values from the boundary of the circle, which are the wave numbers $\kappa_i$, with $i = 0, \dots, N-1$. For each wave number a matrix has to be calculated, see Figure \ref{BeynSubmatrices}. These matrices can be calculated in parallel. In each element of the matrices, integrals have to be calculated numerically. We are using \texttt{scipy.integrate.quad\_vec} which uses the adaptive Gauss-Kronrod 21-point rule. The computation of each matrix is quiet expensive, especially if the domain has multiple holes $q$. For each boundary, the number of faces $n_{f_i}$, with $i=0,\dots,q$, have to be defined. As illustrated in Figure \ref{ParallelizeBeyn1}, a total of $9N q^2 n_f^2$ integrations must be performed for all matrices if the domain consists two holes and the boundaries are discretized with $n_{f_i} = n_f$ faces. 
\tikzmath{\n = 4;}
\begin{figure}[h!]
    \centering
    \begin{tikzpicture}
        \definecolor{pastelgruen}{HTML}{acf6c8}
        \definecolor{pastelrosa}{HTML}{800080}
        \pgfkeys{/pgf/number format/relative*={1}}
        \foreach \i in {1,...,\n}{
            \tikzmath{\midpoint = \i*0.9;}
            \fill[pastelgruen] (\midpoint,0) circle (0.5cm);}
        \foreach \i in {1,...,\n}{
            \tikzmath{\midpoint = 0.9*\i;}
            \draw[thick, black](\midpoint,0) circle (0.5cm);
            \draw[black] (\midpoint,-0.1) -- (\midpoint,0.1) node [below=0.15cm, scale=0.5] {\pgfmathprintnumber{\midpoint}};}
        \draw [->] (0,0) -- (5,0) node [below] {$\mathrm{Re}(\kappa)$};
        \draw [->] (0,-0.75) -- (0,0.75) node [above] {$\mathrm{Im}(\kappa)$};
        \foreach \angle in {0,60,...,300}{
            \tikzmath{\xcross = 0.9+0.5*cos(\angle);
            \ycross = 0.5*sin(\angle);}
            \draw (\xcross,\ycross) node[pastelrosa] {x};
        }
        \tikzmath{\xfirstarrow = 0.9+0.5*cos(240);
        \xsecondarrow = 0.9+0.5*cos(300);
        \yfirstarrow = 0.5*sin(240);
        \ysecondarrow = 0.5*sin(300);}
        \draw [->, thick, pastelrosa] (\xfirstarrow,\yfirstarrow-0.1) --node [left] {$\kappa_i$} (0.3,-1.2) ;
        \draw [->, thick, pastelrosa] (\xsecondarrow+0.1,\ysecondarrow) --node [below left] {$\kappa_j$} (4.6,-1.2) ;
        \node at (0.4,-2)[scale =0.58] {
        $\left(
        \begin{array}{ccc}
          \int\limits_{\Gamma_0} f_0(\kappa_i) \mathrm{d} s(y)  & \cdots & \int\limits_{\Gamma_q}f_0(\kappa_i) \mathrm{d} s(y) \\
          \vdots & \ddots & \vdots \\
          \int\limits_{\Gamma_0} f_q(\kappa_i) \mathrm{d} s(y)& \cdots & \int\limits_{\Gamma_0}f_q(\kappa_i) \mathrm{d} s(y)
        \end{array}
        \right)$
      };
      \node at (4.7,-2)[scale =0.58] {
        $\left(
        \begin{array}{ccc}
          \int\limits_{\Gamma_0} f_0(\kappa_j) \mathrm{d} s(y)  & \cdots & \int\limits_{\Gamma_q}f_0(\kappa_j) \mathrm{d} s(y) \\
          \vdots & \ddots & \vdots \\
          \int\limits_{\Gamma_0} f_q(\kappa_j) \mathrm{d} s(y)& \cdots & \int\limits_{\Gamma_0}f_q(\kappa_j) \mathrm{d} s(y)
        \end{array}
        \right)$
      };
    \end{tikzpicture}  
    \caption{Calculation of the submatrices $M\left(\kappa_i\right)$ for each $i \in 0,\dots N-1$ to find all eigenvalues within one circle}
    \label{BeynSubmatrices}
\end{figure}
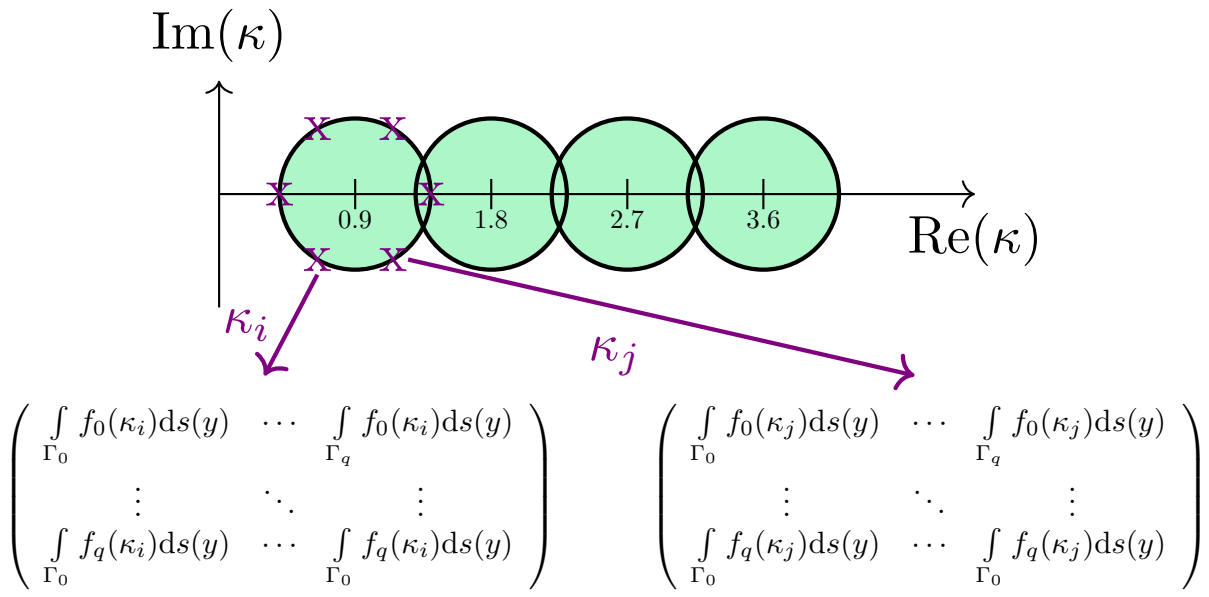

\begin{figure}[h!]
    \centering
    \resizebox{\linewidth}{!}{
    \begin{tikzpicture}
        \def\squarelength{0.08}
        \def\distancesquares{0.12}
        \def\nx{7}
        \def\ny{7}
        \fill[blue!40] (0,0) rectangle ++(2,3);
        \fill[blue!40] (2,1.625) rectangle ++(1,1.375);
        \fill[cyan!40] (2,0) rectangle ++(1,1.625);
        \fill[cyan!40] (3,0) rectangle ++(1,3);
        \fill[cyan!40] (4,0.25) rectangle ++(1,2.75);
        \fill[green!40] (4,0) rectangle ++(1,0.25);
        \fill[green!40] (5,0) rectangle ++(1,3);
        \fill[yellow!40] (9,0) rectangle ++(1,3);
        \fill[yellow!40] (10,1.625) rectangle ++(1,1.375);
        \fill[orange!40] (10,0) rectangle ++(1,1.625);
        \fill[orange!40] (11,0) rectangle ++(1,3);
        \fill[orange!40] (12,0.25) rectangle ++(2,2.75);
        \fill[red!40] (12,0) rectangle ++(1,0.25);
        \fill[red!40] (13,0) rectangle ++(2,3);
        \foreach \N in {0,1,3,4} {
            \foreach \m in {0,...,2} {
                \foreach \n in {0,...,2} {
                    \foreach \i in {0,...,\nx} {
                        \foreach \j in {0,...,\ny} {
                            \draw[gray, very thin] (\N*3+\m+\i*\distancesquares+0.042, \n+\j*\distancesquares+0.042) rectangle ++(\squarelength, \squarelength);
                        }
                    }
                }
            }
        }
        \draw[step = 1cm, black, thick] (0,0) grid (6,3);
        \draw[step = 3cm, black, very thick] (0,0) grid (6,3);
        \node (midpoints) at (7.5,1.5) {\scalebox{1.4}{$\bullet\ \bullet\ \bullet$}};
        \draw[step = 1cm, black, thick] (9,0) grid (15,3);
        \draw[step = 3cm, black, very thick] (9,0) grid (15,3);
        \draw[thick,black,decorate,decoration={brace,mirror,amplitude=12pt}] (0,0) -- (3,0) node[midway, below,yshift=-12pt]{$M\left(\kappa_0\right)$};
        \draw[thick,black,decorate,decoration={brace,mirror,amplitude=12pt}] (3,0) -- (6,0) node[midway, below,yshift=-12pt]{$M\left(\kappa_1\right)$};
        \draw[thick,black,decorate,decoration={brace,mirror,amplitude=12pt}] (9,0) -- (12,0) node[midway, below,yshift=-12pt]{$M\left(\kappa_{N-2}\right)$};
        \draw[thick,black,decorate,decoration={brace,mirror,amplitude=12pt}] (12,0) -- (15,0) node[midway, below,yshift=-12pt]{$M\left(\kappa_{N-1}\right)$};
        \draw[thick,black,decorate,decoration={brace,amplitude=12pt}] (9,3) -- (10,3) node[midway, above,yshift=12pt]{$3 n_{f_0}$};
        \draw[thick,black,decorate,decoration={brace,amplitude=12pt}] (9,2) -- (9,3) node[midway, left,xshift=-12pt]{$3 n_{f_0}$};
        
        \draw[-{Latex[length=3mm, width=3mm]}, line width=2pt]
        (1.5,3.6) -- (1.5,3.1);
        \draw[-{Latex[length=3mm, width=3mm]}, line width=2pt]
        (3.5,3.6) -- (3.5,3.1);
        \draw[-{Latex[length=3mm, width=3mm]}, line width=2pt]
        (5.5,3.6) -- (5.5,3.1);
        \node (above1) at (3.5,4.5) {Integration along the boundary:};
        \node (above2) at (1.5,4.) {$\Gamma_1$};
        \node (above3) at (3.5,4.) {$\Gamma_0$};
        \node (above4) at (5.5,4.) {$\Gamma_2$};
        \draw[-{Latex[length=3mm, width=3mm]}, line width=2pt]
        (15.6,2.5) -- (15.1,2.5);
        \draw[-{Latex[length=3mm, width=3mm]}, line width=2pt]
        (15.6,1.5) -- (15.1,1.5);
        \draw[-{Latex[length=3mm, width=3mm]}, line width=2pt]
        (15.6,0.5) -- (15.1,0.5);
        \node (right1) at (14,3.5) {$x$ is on the boundary:};
        \node (right2) at (16,2.5) {$\Gamma_2$};
        \node (right3) at (16,1.5) {$\Gamma_1$};
        \node (right4) at (16,0.5) {$\Gamma_0$};
    \end{tikzpicture}}
    \caption{Distribution of the different rows within the different submatrices to different cores}
    \label{ParallelizeBeyn1}
\end{figure}
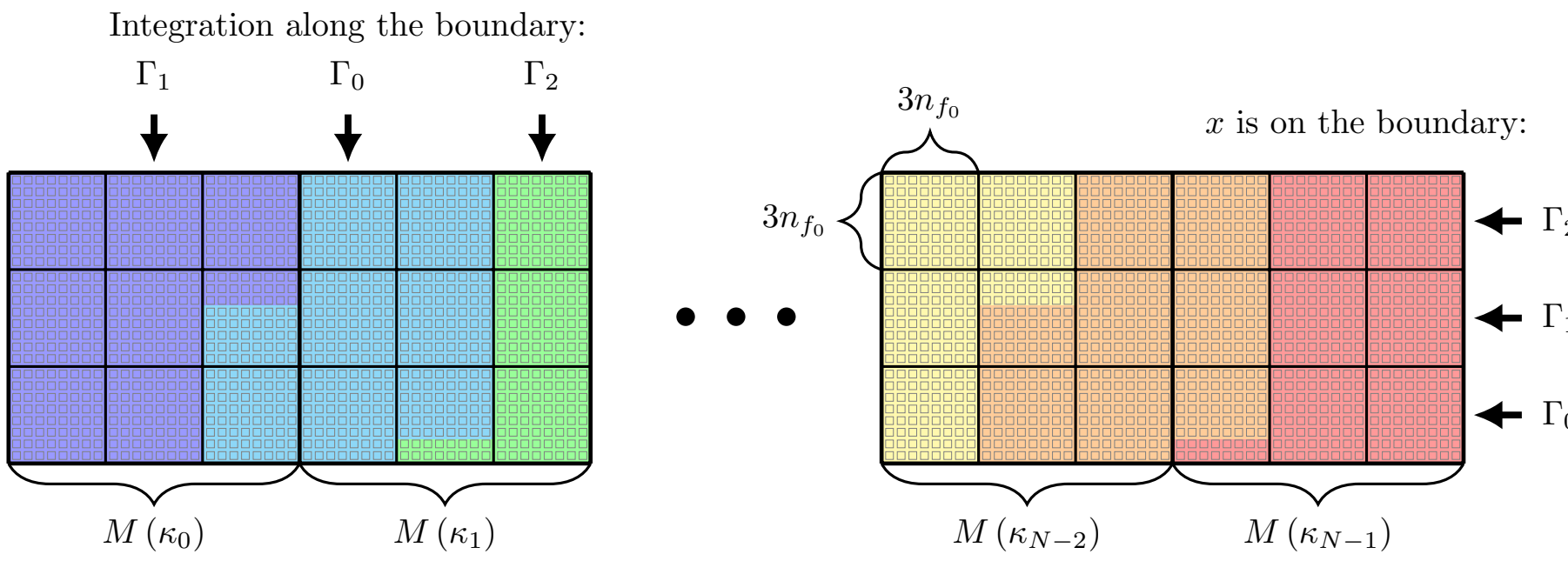

In the following, two ways to parallelize the integration within the matrices $M\left(\kappa_i\right)$ will be explained. Both have different advantages and disadvantages and are better for different circumstances, which will be explained later.
For the first method, the result of the paper \cite[sec. 1.2.3]{Abele2020_Eigenvalues} is used. This says that the required time to evaluate each entry of $M(\kappa)$ for a specific $\kappa$ is equal for all rows and the compute time per row is higher in the main diagonal, where a singularity of the Helmholtz-function occurs.
Based on this, the first algorithm to calculate the matrices is the following.
Previously, the number of rows that each process has to calculate must be determined. This should be distributed as evenly as possible. In Figure \ref{ParallelizeBeyn1}, the distribution to different processes can be seen. Afterwards, the information, which rows should be calculated, is scattered to each of the processes. 
Then, on each process the partial matrices are calculated and are gathered back to the root process. Here, it is scattered again, but this time not to all processes, just to the sub-communicator, which includes $N$ processes. In these processes the right hand-side of the equation
\begin{equation}
     M(\phi(t_r)) x =  \hat{V}\phi'(t_r)\, , \qquad \forall r \in 0,\dots N-1
\end{equation}
is calculated and afterwards, the Equation \eqref{A1Trapez} is solved. 
MPI for Python (MPI4Py) is a Message Passing Interface (MPI) and is used for the whole communication of this algorithm.

\begin{figure}[h!]
    \centering
    
    \begin{subfigure}[t]{0.48\textwidth}
        \includegraphics[width=1\linewidth]{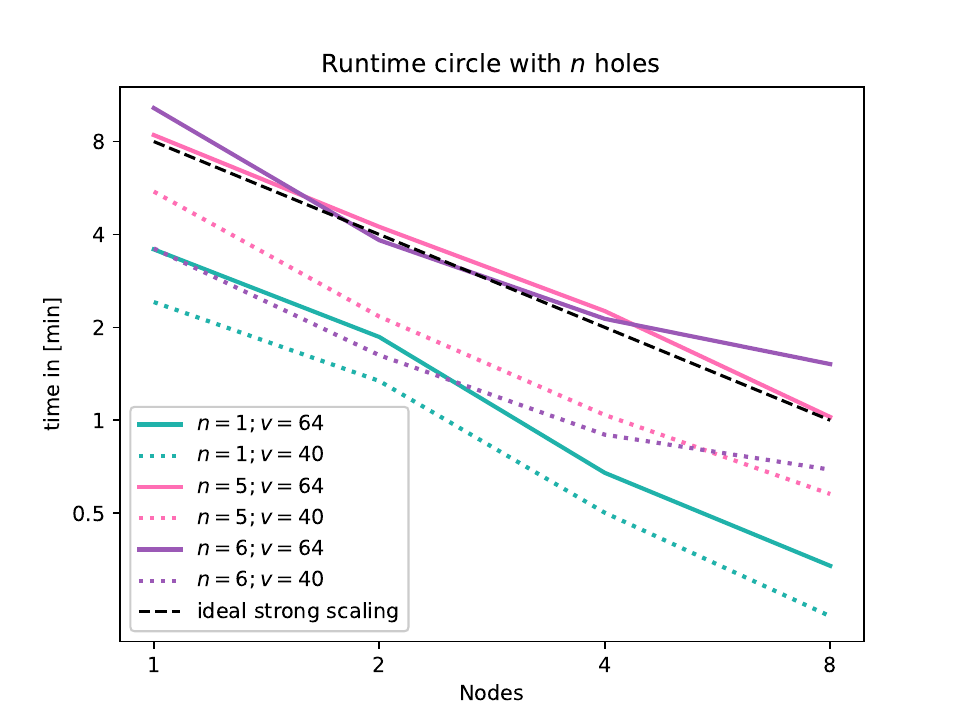}
        \caption{This result is obtained by scattering the different lines which have to be calculated by each process at the beginning}
        \label{fig:Runtime_1}
    \end{subfigure}
    \begin{subfigure}[t]{0.48\textwidth}
        \includegraphics[width=1\linewidth]{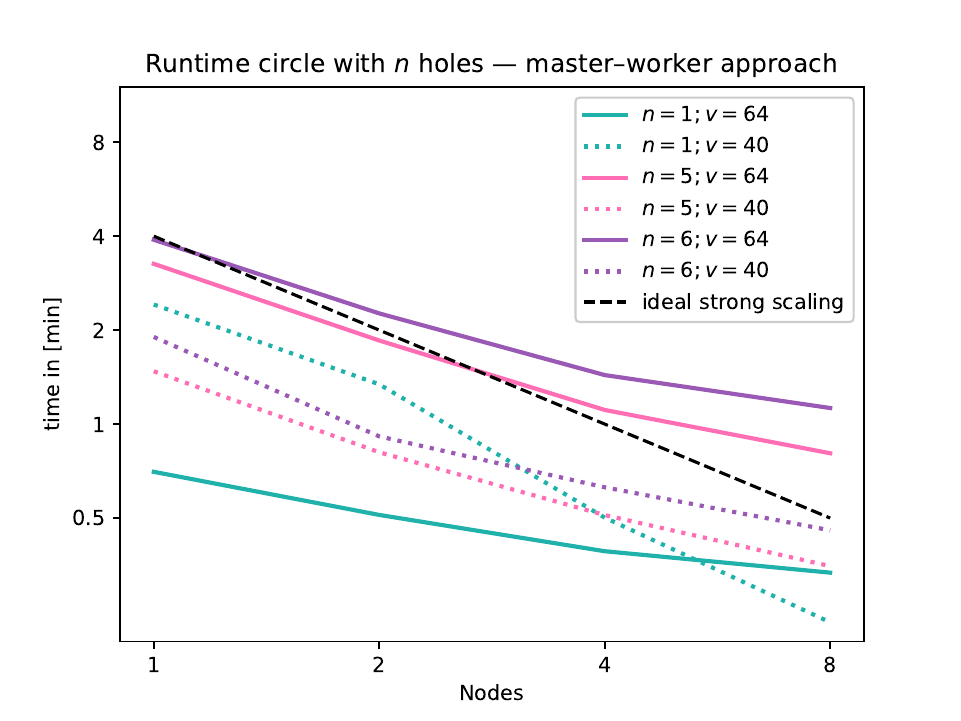}
        \caption{This result is obtained by using the master worker approach.}
        \label{fig:Runtime_master_workers}
    \end{subfigure}
    \caption{Strong scaling for circle with one, five and six holes with 40 and 64 faces}
    \label{fig:Runtime_all}
\end{figure}

By using this approach, a big reduction of the time could be achieved. As can be seen in Figure \ref{fig:Runtime_all}(\subref{fig:Runtime_1}), in the case of one and five holes, the amount of calculation time is halved when the number of nodes increases. A comparison with the black dashed line shows that nearly ideal strong scaling was achieved. In contrast, the runtime for the area with six holes is reduced by slightly less than the scaling factor of $0.57$ when the number of holes is doubled. The algorithm, however, does not take into account that the calculation time per row depends on the given wave number $\kappa$. 

Hence, we need to optimize our previously described algorithm by adding a master-worker approach. This means that we split all given processes in a specific number of master processes, here four, and workers. Afterwards, we assign a specific number of lines in the submatrices to one task, here $15$. Then, all tasks are assigned to a master process. This happens in a circular manner such that each master process task contains every wave number in order to guarantee an approximately equal distribution of working time. As in the previously described algorithm, we scatter a specific number of lines to each worker-process at the beginning. If one of the worker-processes finished its calculations, it sends the result back to the corresponding master-process. This master-process sends the next task to the worker-process. This repeats until all tasks have been completed. With this method, it is not necessary to know at the outset how the individual tasks need to be distributed in order to achieve an even distribution of work. The result can be seen in Figure \ref{fig:Runtime_all}(\subref{fig:Runtime_master_workers}). In comparison to the previous method it can be seen that for smaller number of nodes the calculation time is much smaller if we use the master-worker approach. For more processes, the communication is higher and the slope is not as steep as in the previous algorithm. It can therefore be said that the master-worker approach should be chosen for a small number of nodes and the second approach for a large number of processes.

Besides of the parallelizing part, which was done with \texttt{MPI4Py}, we used \texttt{NumPy} and \texttt{SciPy} in Python to implement the method in a high-performance way. They are written C code and are well-optimized \cite{Harris2020_NumPy, Bressert2012_Scipy, Virtanen2020_SciPy}.
\section{Numerical results}
\label{sec:NumericalResults}
Dahne \cite{Dahne_Counterexample} could show numerically with the finite element method that it is possible to find a counterexample to Payne's nodal line conjecture with six holes, and asks the question, whether it is possible with less holes, moreover what the least amount of holes is.
In this section, we want to present numerically calculated counterexamples with a smaller number of holes and analyze the behavior of the nodal line and the parameter sets rejecting Payne's conjecture for different numbers of holes. 
These results are based on the algorithms explained in the previous chapters.
In Section \ref{sec:ReviewDahne}, the results of Dahne's paper will be reproduced using the BEM method from Section \ref{sec:CollocationMethod}, and a new result with six holes will be illustrated in Section \ref{sec:6Ellipses}. Then, new results with 5 holes will be presented in Section \ref{sec:5Ellipses}. Furthermore, we will focus on the area in which the domain's parameters can lie without the nodal line touching the domain's boundary, refer to Section \ref{sec:ParametersetsRejectingPayne}. Moreover, we will analyze in Section \ref{sec:BehaviorRatioNodalLine}, how the shape of the nodal line rejecting Payne's conjecture changes, varying the parameter set of the holes within the domain. This will give us a bigger knowledge about the behavior and the shape of the nodal line. This knowledge can be used to find easier other counterexamples to Payne's conjecture with fewer holes, or other parameterizations, as was attempted in Section \ref{sec:4_holes}.
The corresponding Python code is available in the GitHub repository \href{https://github.com/mattfeldt/nodal-lines}{\texttt{https://github.com/mattfeldt/nodal\_lines}}.
\subsection{Review of Dahne's finding}
\label{sec:ReviewDahne}
Dahne \cite{Dahne_Counterexample} could find a domain with six holes for which the nodal domain nowhere touches the boundary with the finite element method, and asked the question, what the least amount of holes could be, to reject Payne's conjecture. Based on this, we first try to reproduce these results with the boundary element collocation method. Therefore, we previously used the same domain like in the paper, also seen in Figure \ref{DomainRejectTheoreticallyOrNumerically}(\subref{DomainDahne}). Here, the parameters for the hexagon are given with side length $1$. The six triangles in the middle have the height $\frac{6}{27}$ and centers are placed at a distance of $\frac{11}{27}$ of the origin. That are the numbers written in the paper \cite{Dahne_Counterexample} directly under Figure 1(a). In contrast, Section 2 ``Finding a candidate'' states that the distance should be $\frac{14}{27}$. We proofed both parameter sets and came to the following outcomes.
\begin{figure}[h!]
    \centering
    \begin{subfigure}{0.32\textwidth}
        \includegraphics[width=1\linewidth]{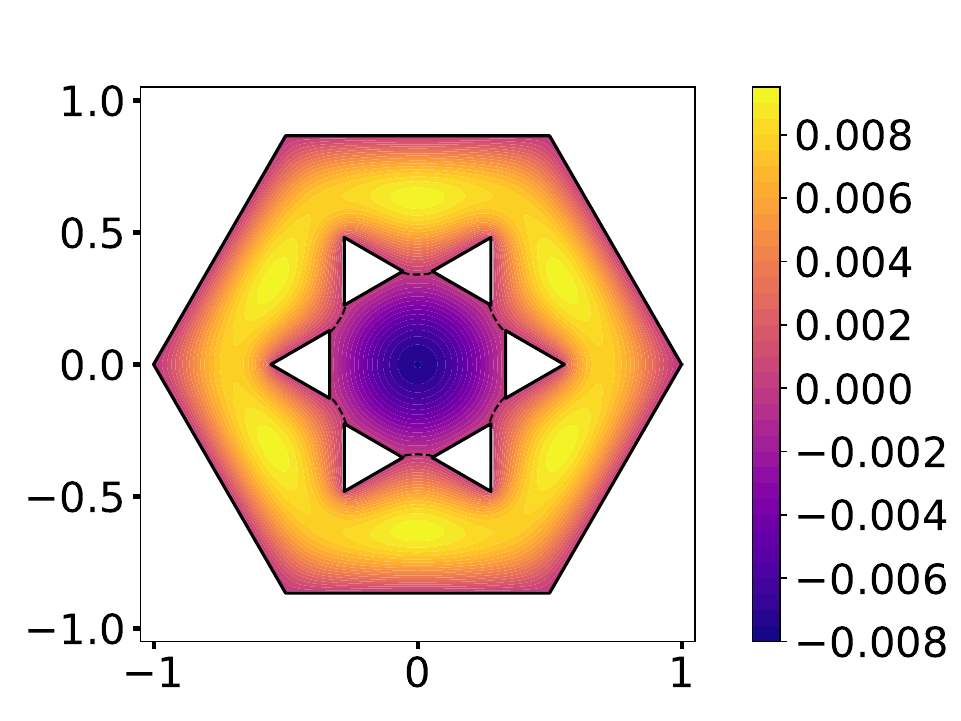}
        \caption{Centers of the triangles are placed at a distance of $\frac{11}{27}$ from the origin with $\lambda_2 \approx 49.8069$.}
        \label{fig:Dahne_11_27}
    \end{subfigure}
    \begin{subfigure}{0.32\textwidth}
        \includegraphics[width=1\linewidth]{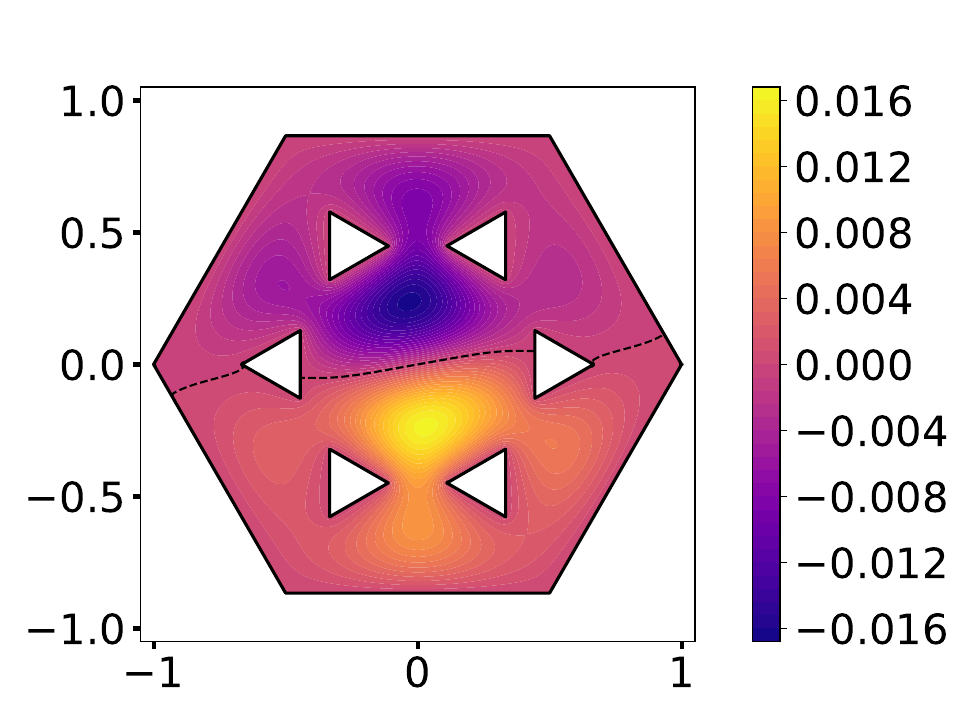}
        \caption{Centers of the triangles are placed at a distance of $\frac{14}{27}$ from the origin with $\lambda_{2,3} \approx 61.7958$.}
        \label{fig:Dahne_14_27}
    \end{subfigure}
    \begin{subfigure}{0.32\textwidth}
        \includegraphics[width=1\linewidth]{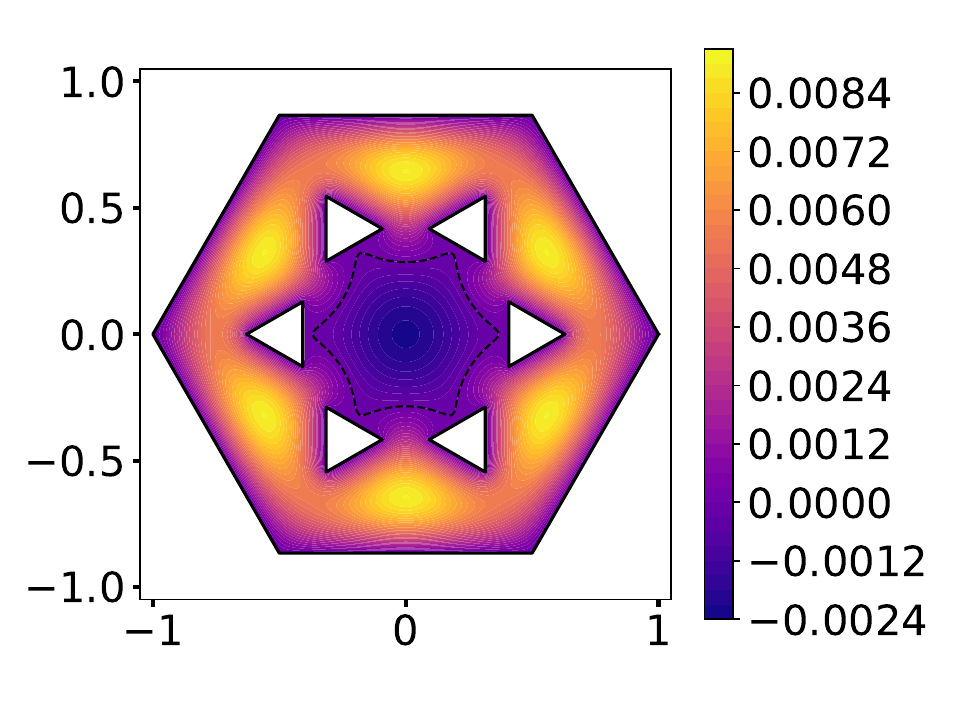}
        \caption{Centers of the triangles are placed at a distance of $\frac{13}{27}$ from the origin with $\lambda_2 \approx 63.2692$.}
        \label{fig:Dahne_13_27}
    \end{subfigure}
    \caption{Real part of the second eigenfunction $\operatorname{Re}(u_2)$ in the inner part of an hexagon with six triangular regions as holes with the dashed nodal line.}
    \label{fig:Dahne_all_results}
\end{figure}
The first case with a distance of $\frac{11}{27}$ is presented in Figure \ref{fig:Dahne_all_results}(\subref{fig:Dahne_11_27}). The nodal line does not touch the boundary of the hexagon. Instead, it connects all the triangles' boundaries in a circular fashion, splitting the domain into an inner and an outer part. 
The first three eigenvalues are $\lambda_1 \approx 47.4390$, $\lambda_2 \approx 49.8069$ and $\lambda_{3,4} \approx 51.6727$, where the first two are simple eigenvalues and the third eigenvalue is a double.

The second case is shown in Figure \ref{fig:Dahne_all_results}(\subref{fig:Dahne_14_27}). Here, the distance from the origin to the center of each triangle is $\frac{14}{27}$. As can be seen, the nodal line splits the domain into two halves. It starts at the boundary of the hexagon and goes nearly straight throw the domain, touching two triangles along the way. The first four eigenvalues are $\lambda_1 \approx 25.762$, $\lambda_{2,3} \approx 61.7958$, and $\lambda_4 \approx 70.1135$ where the first is a simple eigenvalue, the second eigenvalue is a double and the fourth is again a simple. 

Both cases fulfill Payne's conjecture and are not counterexamples.
Therefore, we introduce a third case with a distance of $\frac{13}{27}$. Choosing this distance enables us to find a solution with six holes that does not touch the boundary. When taking a deeper look into the Figure \ref{fig:Dahne_all_results}(\subref{fig:Dahne_13_27}), we see that the nodal line looks like a hypotrochoid. This is why, in Section \ref{sec:FindParameters}, we chose the hypotrochoid as an additional hole in the algorithm to find future domains for calculating the first eigenvalue. The first three eigenvalues are $\lambda_1 \approx 31.0635$, $\lambda_2 \approx 63.2692$ and $\lambda_{3,4} \approx 63.7979$, where the first two are simple eigenvalues and the third eigenvalue is a double.
\subsection{Six ellipses}
\label{sec:6Ellipses}
Based on the results of Section \ref{sec:FindParameters}, we attempt to find additional results with six holes that have smooth boundaries. 
Therefore, we define two domains.
The first domain is defined as a disc with six elliptical holes inside, what can be seen in Figure \ref{DomainCircleEllipsesHypotrochoid2}(\subref{DomainCircleEllipses}). For this domain the second eigenvalue will be calculated. The second domain has an additional hole bounded by a hypotrochoid in the interior and can be seen in Figure \ref{DomainCircleEllipsesHypotrochoid2}(\subref{DomainCircleEllipsesHypotrochoid}). It is used to calculate the first eigenvalue. As explained in Section \ref{sec:FindParameters} we use an algorithm that reduces the distance of the two eigenvalues. If they are close enough it is probably possible to find a solution for this kind of domain.
We chose the initial parameter set $a_0=2$, $b_0=2$, $a_i=0.4$, $b_i = 0.1$ for $i=1,\dots,6$, $r_E= 1.1$, $x_o = 0.6$, and $x_i=0.5$.
In Table \ref{tab:Tabular_circle_six_ellipses} the variation of the parameters over the iterations is visualized.
\begin{table}[h!]
    \centering
    \resizebox{\linewidth}{!}{
    {\definecolor{azure(colorwheel)}{rgb}{0.0, 0.5, 1.0}
    \rowcolors{3}{azure(colorwheel)!30}{white}
    \scriptsize
    \begin{tabular}{ |p{2.3cm}|p{2.3cm}|p{2.3cm}|p{2.3cm}|p{2.3cm}|p{2.3cm}|  }
    \hline
    \multicolumn{6}{|c|}{Parameters List} \\
    \hline
    $a$ & $b$ & $r_E$ & second eigen-\newline value(s) without \newline hypotrochoid & first eigenvalue \newline with hypotrochoid & error\\
    \hline
    0.4 & 0.1 & 1.1 & 12.75443368 \newline -9.98626342-08$\mathrm{i}$; \newline 12.75443368\newline -9.51668044e-08$\mathrm{i}$ & 14.15666489\newline -1.31673094e-06$\mathrm{i}$  & 1.40223121\\
    0.4 & 0.10000001490116 & 1.1 & 12.75443404 \newline-9.51540354e-08$\mathrm{i}$;  \newline 12.75443404 \newline -9.90266654e-08$\mathrm{i}$ & 14.15666526 \newline-1.31673113e-06$\mathrm{i}$ & 1.40223122\\
    0.40000001490116 & 0.1 & 1.1 & 12.75443383\newline -9.47759266e-08$\mathrm{i}$; \newline 12.75443383 \newline -9.33267139e-08$\mathrm{i}$ & 14.15666519 \newline -1.31673094e-06$\mathrm{i}$ & 1.40223136\\
    0.4 & 0.1 & 1.1000000163913 & 12.75443325 \newline-9.26709261e-08$\mathrm{i}$; \newline 12.75443325\newline -9.64124502e-08$\mathrm{i}$ & 14.15666534 \newline-1.31673112e-06$\mathrm{i}$ & 1.40223209\\
    0.39189243830382 & 0.15592081791325 & 1.07828534328197 & 14.63365226 \newline-1.64817744e-06$\mathrm{i}$ & 14.62972897\newline -1.64556862e-07$\mathrm{i}$ & 0.00392329\\
    0.39189243830382 & 0.15592083281441 & 1.07828534328197 & 14.63365264\newline -1.63853124e-06$\mathrm{i}$ & 14.62972928\newline -1.64556855e-06$\mathrm{i}$ & 0.00392336\\
    0.39189245320498 & 0.15592081791325 & 1.07828534328197 & 14.63365257\newline -1.64435837e-06$\mathrm{i}$ & 14.6297292\newline -1.64556861e-06$\mathrm{i}$ & 0.00392337\\
    0.39189243830382 & 0.15592081791325 & 1.07828535934967 & 14.63365272\newline -1.63873887e-06$\mathrm{i}$ & 14.62972935\newline -1.64556847e-06$\mathrm{i}$ & 0.00392337\\
    0.39126299243500 & 0.15584574005610 & 1.07497307444190 & 14.52885098\newline -1.60168056e-06$\mathrm{i}$ & 14.52507493 \newline-1.60186807e-06$\mathrm{i}$ & 0.00377605\\
    \hline
    \end{tabular}}}
    \caption{Parameter set, results of the first and second eigenvalue and given error between the results changing during the algorithm defined in Section \ref{sec:FindParameters}}
    \label{tab:Tabular_circle_six_ellipses}
\end{table}
Previously, the first and second eigenvalue for the two problems are calculated for the first parameter set. Afterwards, each individual parameter is tweaked slightly, and the effect of this change on the error between the two eigenvalues is observed. Based on this, all parameters are updated and the described process begins again.
While the second eigenvalue is always a double in the first process, it is a simple in the subsequent ones. This already indicates that this is a possible counterexample.
The corresponding eigenfunctions can be observed in Figure \ref{DomainCircleEllipsesHypotrochoid2}. 

\begin{figure}[h!]
    \centering
    \begin{subfigure}{0.32\textwidth}
        \includegraphics[width=1\linewidth]{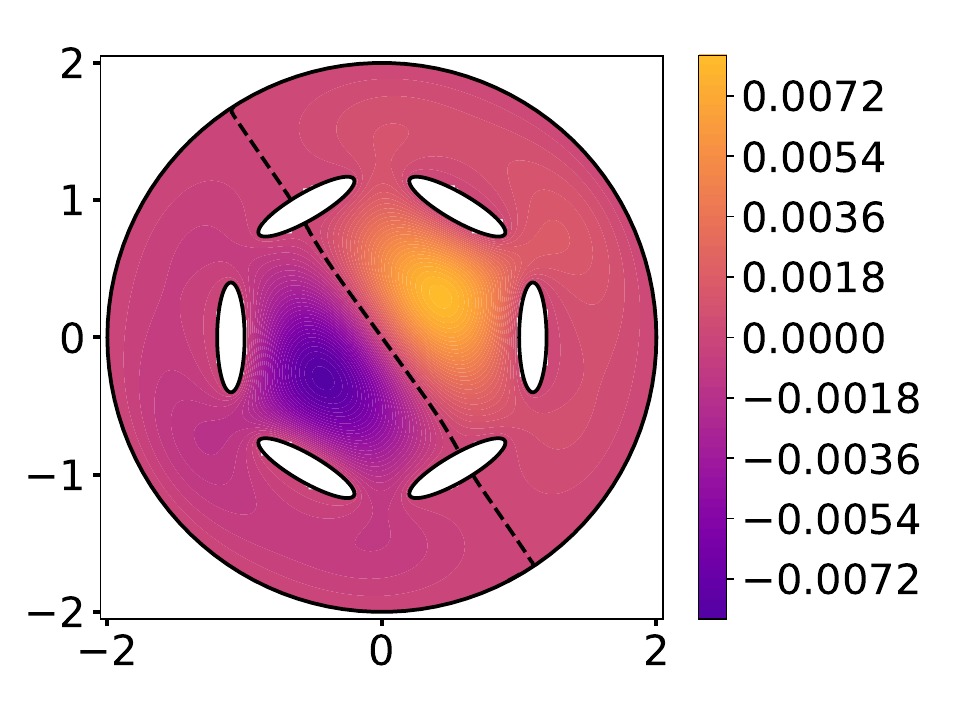}
        \caption{Linear color bar \newline $\begin{aligned}
            a_i &=0.4\\ b_i &= 0.1 \\ r_E &=1.1 \\ \lambda_{2,3} &\approx 12.7544
        \end{aligned}$}
        \label{fig:Circle_six_ellipses_first}
    \end{subfigure}
    \begin{subfigure}{0.32\textwidth}
        \includegraphics[width=1\linewidth]{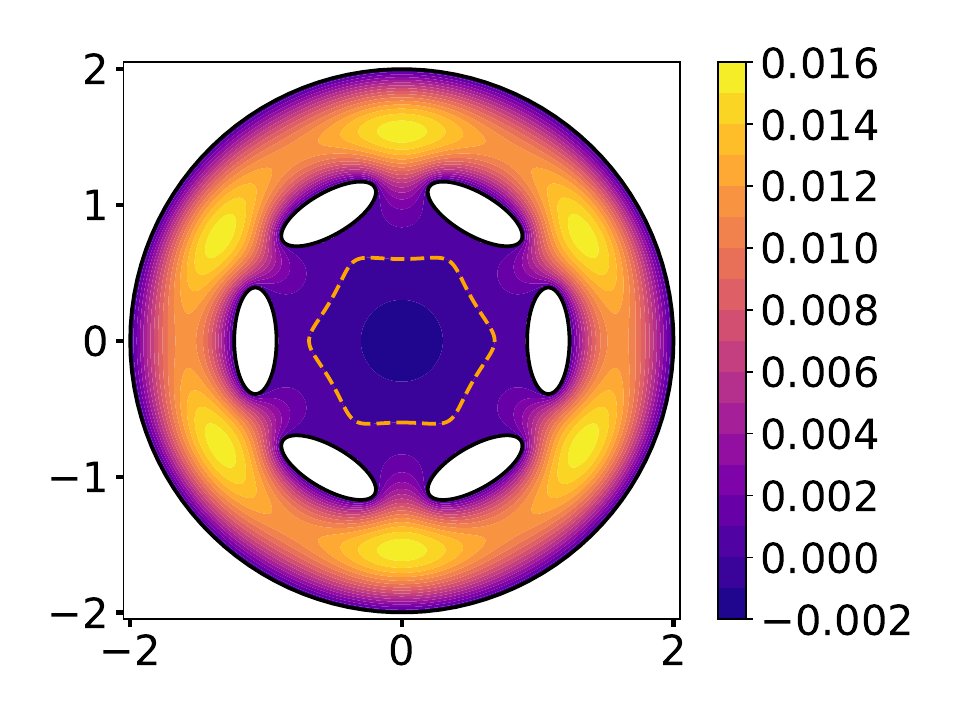}
        \caption{Linear color bar \newline $\begin{aligned}
        a_i&=0.39189243830382\\b_i &= 0.15592081791325\\  r_E&=1.07828534328197 \\ \lambda_2 &\approx 14.6335
        \end{aligned}$}
        \label{fig:Circle_six_ellipses_second}
    \end{subfigure}
    \begin{subfigure}{0.32\textwidth}
        \includegraphics[width=1\linewidth]{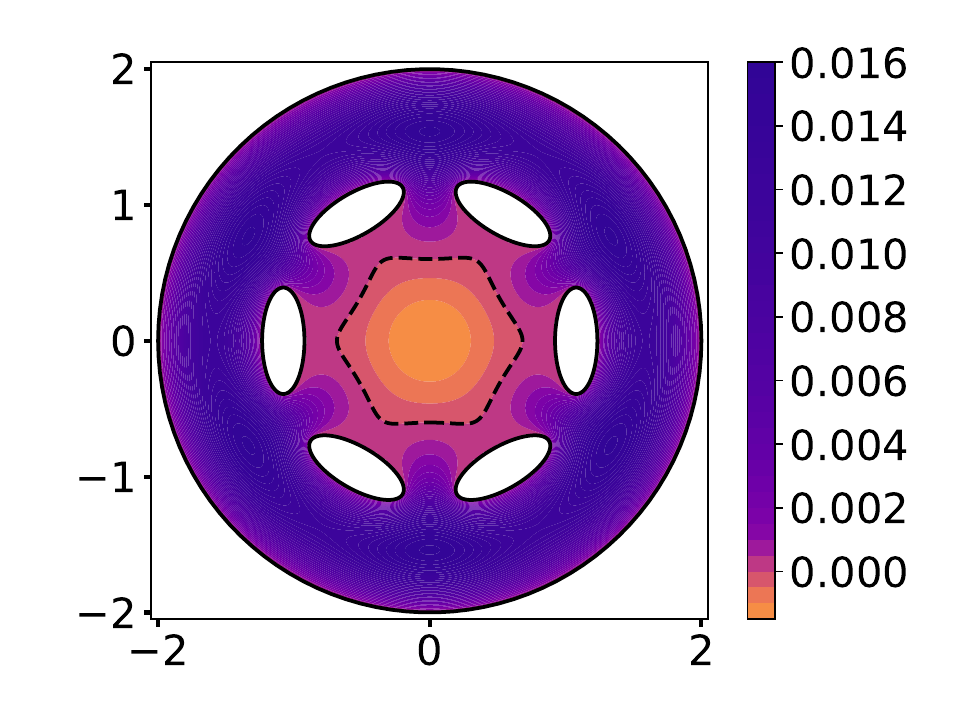}
        \caption{Modified color bar \newline $\begin{aligned}
            a_i&=0.39189243830382 \\b_i &= 0.15592081791325\\ r_E &=1.07828534328197 \\ \lambda_2 &\approx 14.6335
            \end{aligned}$}
        \label{fig:Circle_six_ellipses_second_symlogplot}
    \end{subfigure}

    \begin{subfigure}{0.32\textwidth}
        \includegraphics[width=1\linewidth]{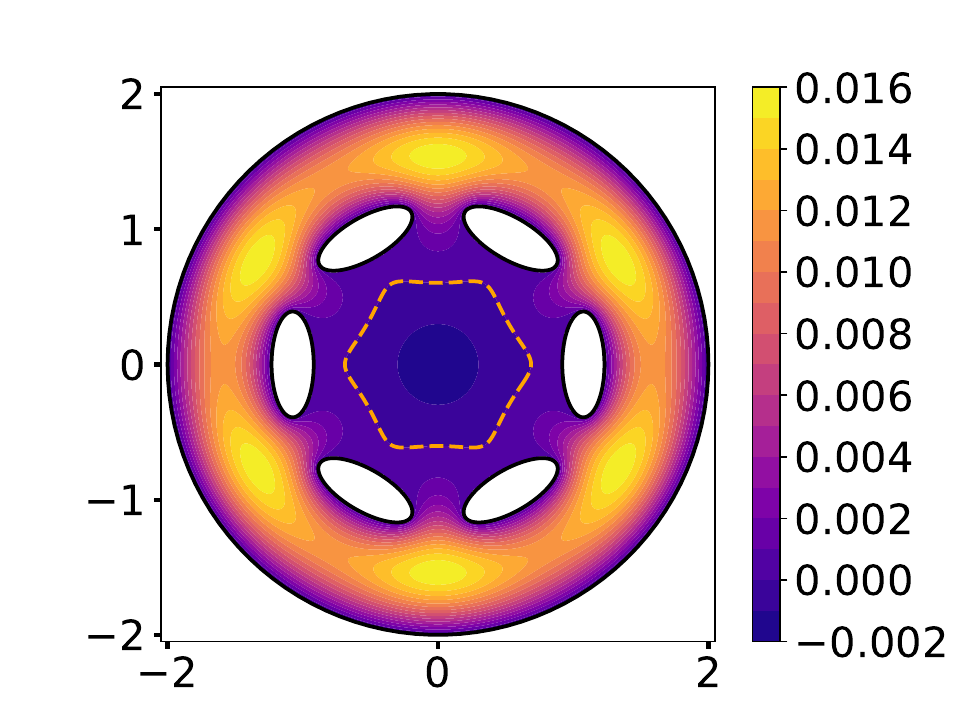}
        \caption{Linear color bar \newline $\begin{aligned}
            a_i&=0.39126299243500 \\b_i &= 0.15584574005610\\ r_E &=1.07497307444190 \\ \lambda_2 &\approx 14.5287
            \end{aligned}$}
        \label{fig:Circle_six_ellipses_third}
    \end{subfigure}
    \begin{subfigure}{0.32\textwidth}
        \includegraphics[width=1\linewidth]{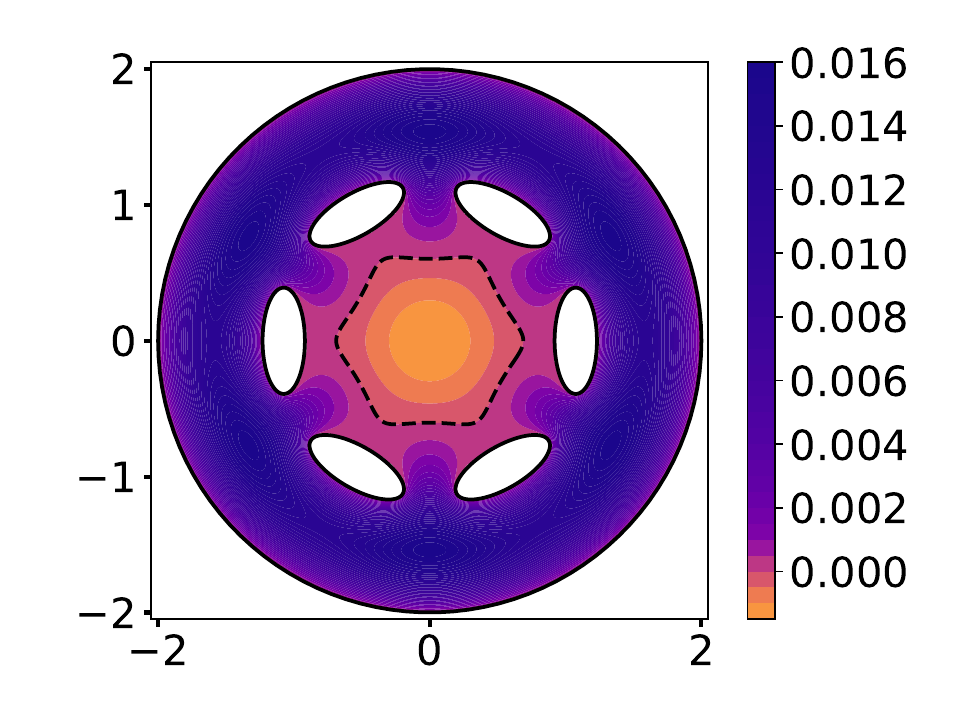}
        \caption{Modified color bar \newline $\begin{aligned}
            a_i&=0.39126299243500 \\b_i &= 0.15584574005610\\ r_E &=1.07497307444190 \\ \lambda_2 &\approx 14.5287
            \end{aligned}$}
        \label{fig:Circle_six_ellipses_third_symlogplot}
    \end{subfigure}
    \caption{Real part of the second eigenfunction $\operatorname{Re}(u_2)$ in the inner part of a disc with six elliptical holes with different parameter sets for the elliptical holes with $i=1,\dots,6$ and a linear or a modified color bar, which means a mix of a linear and a logarithmic color bar. Here, the nodal line is shown as a dashed line.}
    \label{fig:Circle_six_ellipses_all_results}
\end{figure}
The second eigenfunction for the initial parameter set can be seen in Figure \ref{fig:Circle_six_ellipses_all_results}(\subref{fig:Circle_six_ellipses_first}). Here, the first four corresponding eigenvalues are $\lambda_1 \approx 5.0988$, $\lambda_{2,3} \approx 12.7544$, and $\lambda_{4} \approx 14.1683$ and the nodal line of the real part of second eigenfunction touches the boundary. In the last four Figures from \ref{DomainCircleEllipsesHypotrochoid2}(\subref{fig:Circle_six_ellipses_second}) to \ref{DomainCircleEllipsesHypotrochoid2}(\subref{fig:Circle_six_ellipses_third_symlogplot}) the nodal line of the second eigenfunction does not touch the boundary and again looks like a hypotrochoid, this time with the parameter set $a_i=0.39189243830382$, $b_i = 0.15592081791325$, and $r_E =1.07828534328197$ just as $a_i=0.39126299243500$, $b_i = 0.15584574005610$, and $r_E =1.07497307444190$. For a better understanding how the second eigenfunction behaves in the inner part, where the nodal line can be seen, we added the two Figures \ref{DomainCircleEllipsesHypotrochoid2}(\subref{fig:Circle_six_ellipses_second_symlogplot}) and \ref{DomainCircleEllipsesHypotrochoid2}(\subref{fig:Circle_six_ellipses_third_symlogplot}). These two plots show the results for the same parameter sets as in \ref{DomainCircleEllipsesHypotrochoid2}(\subref{fig:Circle_six_ellipses_second}) and \ref{DomainCircleEllipsesHypotrochoid2}(\subref{fig:Circle_six_ellipses_third}), but the color bars have been modified such that values close to 0 are distributed linearly and values farther from 0 are distributed logarithmically. For Figures \ref{DomainCircleEllipsesHypotrochoid2}(\subref{fig:Circle_six_ellipses_second}) and \ref{DomainCircleEllipsesHypotrochoid2}(\subref{fig:Circle_six_ellipses_second_symlogplot}), the corresponding first three eigenvalues are $\lambda_1 \approx 5.9207$, $\lambda_{2} \approx 14.6335$, and $\lambda_{3,4} \approx 14.6515$. The first four eigenvalues for the Figures \ref{DomainCircleEllipsesHypotrochoid2}(\subref{fig:Circle_six_ellipses_third}) and \ref{DomainCircleEllipsesHypotrochoid2}(\subref{fig:Circle_six_ellipses_third_symlogplot}) are $\lambda_1 \approx 5.9632$, $\lambda_2 \approx 14.5287$ and $\lambda_{3,4} \approx 14.6607$, where the first two are simple eigenvalues and the third eigenvalue is a double. The results allow us to show that the algorithm we developed is capable of finding new areas to refute Payne's hypothesis. Furthermore, this result demonstrates that it is with more than one parameter set possible to reject Payne's conjecture. The question appears, whether it is possible to find closed areas within the parameters where the nodal line never touches the boundary. Moreover, it would be interesting to find out, how these areas behave when the number of holes is reduced. 
\subsection{Five ellipses}
\label{sec:5Ellipses}
\begin{figure}[h!]
    \centering
    \begin{subfigure}{0.32\textwidth}
        \includegraphics[width=1\linewidth]{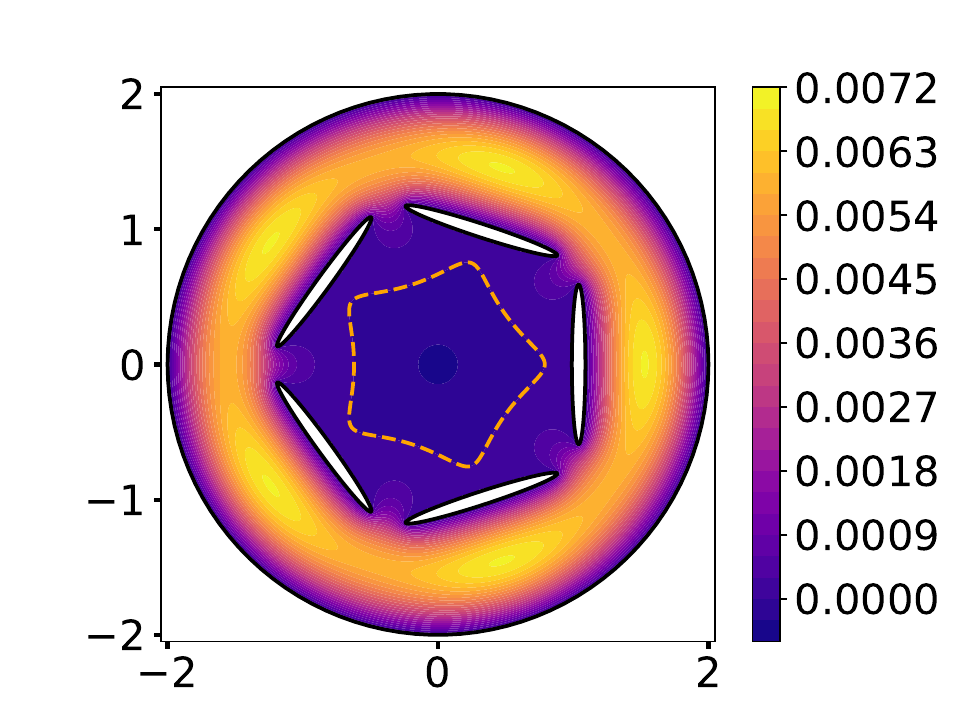}
        \caption{$\begin{aligned}
            a_i&=\SI{0.59}{cm} \\ b_i &= \SI{0.05}{cm}\\r_E &=\SI{1.04}{cm}\\\lambda_2 &\approx 13.0931
            \end{aligned}$}
        \label{fig:Circle_five_ellipses_first}
    \end{subfigure}
    \begin{subfigure}{0.32\textwidth}
        \includegraphics[width=1\linewidth]{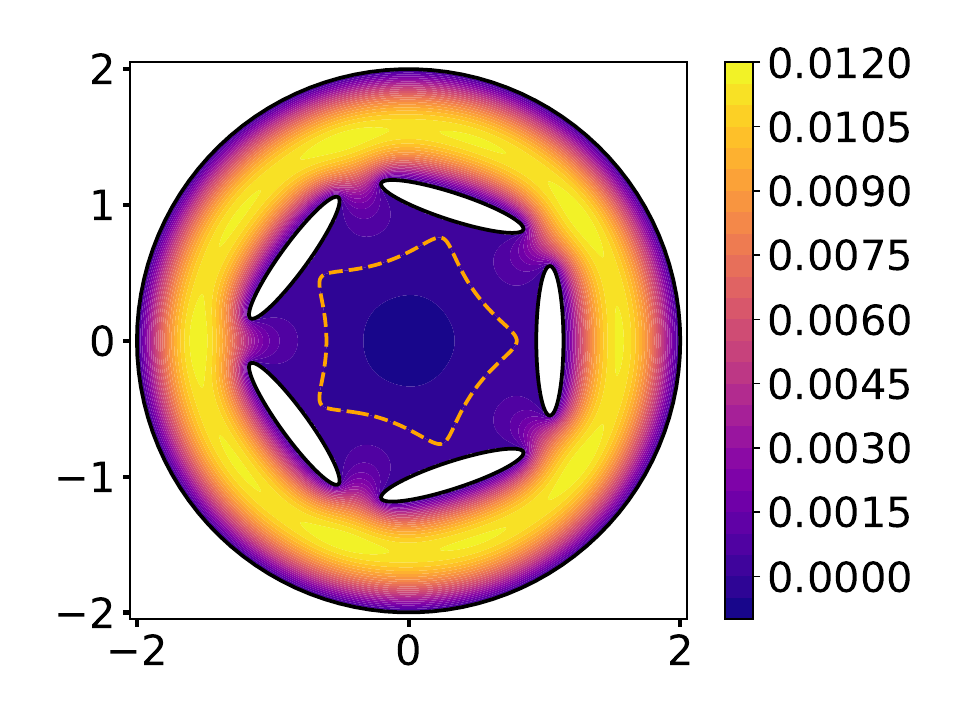}
        \caption{$\begin{aligned}
            a_i &=\SI{0.55}{cm}\\ b_i &= \SI{0.1}{cm} \\ r_E &=\SI{1.04}{cm}\\\lambda_2 &\approx 13.7139
        \end{aligned}$}
        \label{fig:Circle_five_ellipses_second}
    \end{subfigure}
    \begin{subfigure}{0.32\textwidth}
        \includegraphics[width=1\linewidth]{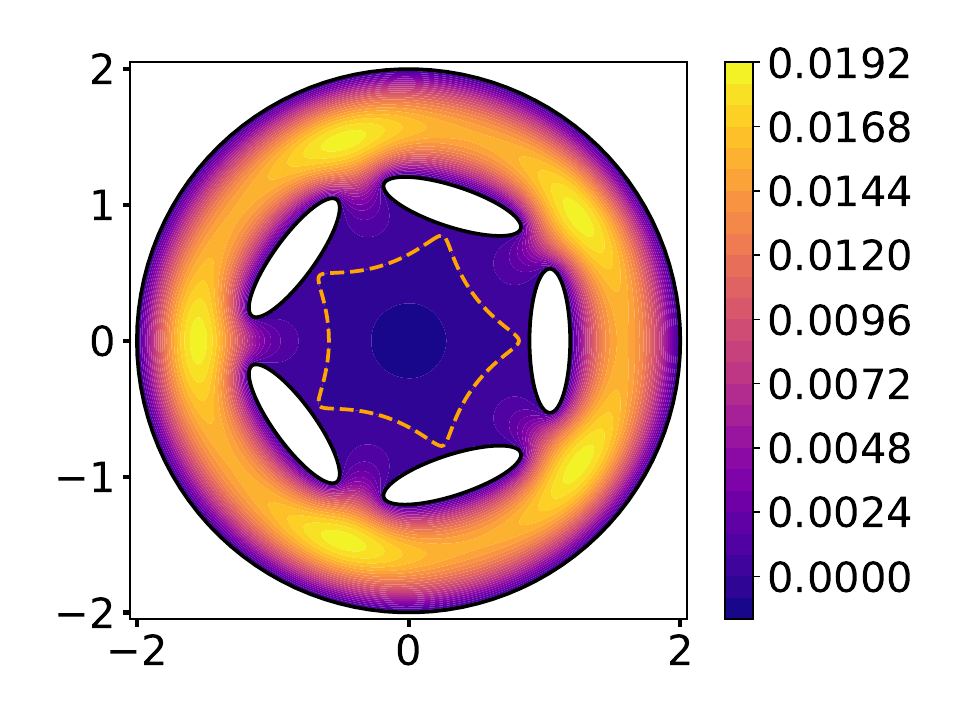}
        \caption{$\begin{aligned}
        a_i&=\SI{0.53}{cm}\\ b_i &= \SI{0.15}{cm}\\ r_E&=\SI{1.04}{cm}\\\lambda_2 &\approx 14.5012
        \end{aligned}$}
        \label{fig:Circle_five_ellipses_third}
    \end{subfigure}
    \caption{Results of a circle with five ellipses for different parameter sets for the ellipses with $i=1,\dots,5$. Here, the nodal line is shown as a dashed line.}
    \label{fig:Circle_five_ellipses_all_results}
\end{figure}

In this subsection, we show the results achieved by trying to reduce the number of holes. Therefore, we used the same method as in the previous Section \ref{sec:6Ellipses}. Depending on the initial values we got diverse results. Some of the results can be seen in Figure \ref{fig:Circle_five_ellipses_all_results}. 
Here, we have three results with the same distance of the ellipses to the origin of the circle $r_E=\SI{1.04}{cm}$. The only thing that changes is the length $a$ and the width $b$ of the ellipses. In all cases, the nodal line appears to be a hypotrochoid and does not touch the boundary. The shape of the nodal line changes by varying the length and width of the ellipses. The smaller the length $a$ and the bigger the width $b$, the more pronounced the lobes of the hypotrochoid become, which are the outward bulges of the curve. Moreover, the bigger the length and the smaller the width of the ellipses, the closer are the second and the third eigenvalues next to each other. The eigenvalues corresponding to Figure \ref{fig:Circle_five_ellipses_all_results}(\subref{fig:Circle_five_ellipses_first}) are $\lambda_1 \approx 5.2265$, $\lambda_2 \approx 13.0931$ and $\lambda_{3,4} \approx 13.1311$. For Figure  \ref{fig:Circle_five_ellipses_all_results}(\subref{fig:Circle_five_ellipses_second}) the corresponding eigenvalues are $\lambda_1 \approx 5.7034$, $\lambda_2 \approx 13.7139$ and $\lambda_{3,4} \approx 14.0249$. The corresponding eigenvalues for Figure \ref{fig:Circle_five_ellipses_all_results}(\subref{fig:Circle_five_ellipses_third}) are $\lambda_1 \approx 6.2780$, $\lambda_2 \approx 14.5012$ and $\lambda_{3,4} \approx 14.855$. If we increase the length $a$ and decrease the width $b$ of Figure \ref{fig:Circle_five_ellipses_all_results}(\subref{fig:Circle_five_ellipses_first})'s domain even more, the second and third eigenvalue would switch, and the nodal line of the second eigenfunction would touch the boundary.

Hence, we have more than one counterexample, similar to the domain with six holes, see Section \ref{sec:6Ellipses}. We expect to have a closed area of parameter sets, for which the nodal line does not touch the boundary. In the next Section \ref{sec:ParametersetsRejectingPayne}, we analyze the areas of the parameter sets for five and six holes, to get knowledge about the behavior of the nodal line. For example, we want to investigate whether the set of parameter sets rejecting Payne's conjecture is really a closed area, or how the set of parameter sets changes for different number of holes and for different distances $r_E$.
\begin{figure}[h!]
    \centering
        \begin{subfigure}[b]{0.3\textwidth}
        \begin{tikzpicture}[scale=2, transform shape]
            \tikzmath{\Ahyp=6; \xohyp = 0.3; \xihyp = 0.7;
            \rhyp = (\xohyp + \xihyp)/(2*(\Ahyp-1));
            \Rhyp = \Ahyp*\rhyp;
            \ahyp = (\Ahyp-1)*(\xohyp-\xihyp)/(\xohyp+\xihyp);
            \xmhyp=0; \ymhyp=0;\alphahyp=0;};
            \definecolor{pastelgruen}{HTML}{acf6c8}
            \draw[domain=0:360, samples=1000, variable=\t, fill = pastelgruen] plot({\xmhyp +((\Rhyp-\rhyp)*cos(\t)+\ahyp*\rhyp*cos((\Rhyp/\rhyp-1)*\t))*cos(\alphahyp)-((\Rhyp-\rhyp)*sin(\t)-\ahyp*\rhyp*sin((\Rhyp/\rhyp-1)*\t))* sin(\alphahyp)
          },{\ymhyp + ((\Rhyp-\rhyp)*cos(\t)+\ahyp*\rhyp*cos((\Rhyp/\rhyp-1)*\t))*sin(\alphahyp)+ ((\Rhyp-\rhyp)*sin(\t)-\ahyp*\rhyp*sin((\Rhyp/\rhyp-1)*\t))* cos(\alphahyp)});
            \draw[thick, blue] (0,0) --  (0,-\xihyp)  node[midway,right,sloped,rotate=90, xshift=0, yshift=-1.5, scale = 0.6, transform shape]{$\scriptstyle x_i$};
            \draw[thick, blue] (0,0) --  ({\xohyp*cos(120)},{\xohyp*sin(120)})  node[midway,right,sloped,rotate=90, xshift=-1.5, yshift=-1, scale = 0.6, transform shape]{$\scriptstyle x_o$};
        \end{tikzpicture}
        \caption{$x_o = 0.3$, $x_i = 0.7$}
        \label{DomainHypotrochoid1}
        \end{subfigure}
        \begin{subfigure}[b]{0.3\textwidth}
        \begin{tikzpicture}[scale=2, transform shape]
            \tikzmath{\Ahyp=6; \xohyp = 0.7; \xihyp = 0.5;
            \rhyp = (\xohyp + \xihyp)/(2*(\Ahyp-1));
            \Rhyp = \Ahyp*\rhyp;
            \ahyp = (\Ahyp-1)*(\xohyp-\xihyp)/(\xohyp+\xihyp);
            \xmhyp=0; \ymhyp=0;\alphahyp=0;};
            \definecolor{pastelgruen}{HTML}{acf6c8}
            \draw[domain=0:360, samples=1000, variable=\t, fill = pastelgruen] plot({\xmhyp +((\Rhyp-\rhyp)*cos(\t)+\ahyp*\rhyp*cos((\Rhyp/\rhyp-1)*\t))*cos(\alphahyp)-((\Rhyp-\rhyp)*sin(\t)-\ahyp*\rhyp*sin((\Rhyp/\rhyp-1)*\t))* sin(\alphahyp)
          },{\ymhyp + ((\Rhyp-\rhyp)*cos(\t)+\ahyp*\rhyp*cos((\Rhyp/\rhyp-1)*\t))*sin(\alphahyp)+ ((\Rhyp-\rhyp)*sin(\t)-\ahyp*\rhyp*sin((\Rhyp/\rhyp-1)*\t))* cos(\alphahyp)});
            \draw[thick, blue] (0,0) --  (0,-0.5)  node[midway,right,sloped,rotate=90, xshift=0, yshift=-1.5, scale = 0.6, transform shape]{$\scriptstyle x_i$};
            \draw[thick, blue] (0,0) --  ({0.7*cos(120)},{0.7*sin(120)})  node[midway,right,sloped,rotate=90, xshift=-1.5, yshift=-1, scale = 0.6, transform shape]{$\scriptstyle x_o$};
        \end{tikzpicture}
        \caption{$x_o = 0.7$, $x_i = 0.5$}
        \label{DomainHypotrochoid2}
        \end{subfigure}
        \begin{subfigure}[b]{0.3\textwidth}
        \begin{tikzpicture}[scale=2, transform shape]
            \tikzmath{\Ahyp=6; \xohyp = 0.7; \xihyp = 0.7;
            \rhyp = (\xohyp + \xihyp)/(2*(\Ahyp-1));
            \Rhyp = \Ahyp*\rhyp;
            \ahyp = (\Ahyp-1)*(\xohyp-\xihyp)/(\xohyp+\xihyp);
            \xmhyp=0; \ymhyp=0;\alphahyp=0;};
            \definecolor{pastelgruen}{HTML}{acf6c8}
            \draw[domain=0:360, samples=1000, variable=\t, fill = pastelgruen] plot({\xmhyp +((\Rhyp-\rhyp)*cos(\t)+\ahyp*\rhyp*cos((\Rhyp/\rhyp-1)*\t))*cos(\alphahyp)-((\Rhyp-\rhyp)*sin(\t)-\ahyp*\rhyp*sin((\Rhyp/\rhyp-1)*\t))* sin(\alphahyp)
          },{\ymhyp + ((\Rhyp-\rhyp)*cos(\t)+\ahyp*\rhyp*cos((\Rhyp/\rhyp-1)*\t))*sin(\alphahyp)+ ((\Rhyp-\rhyp)*sin(\t)-\ahyp*\rhyp*sin((\Rhyp/\rhyp-1)*\t))* cos(\alphahyp)});
            \draw[thick, blue] (0,0) --  (0,-\xihyp)  node[midway,right,sloped,rotate=90, xshift=0, yshift=-1.5, scale = 0.6, transform shape]{$\scriptstyle x_i$};
            \draw[thick, blue] (0,0) --  ({0.7*cos(120)},{0.7*sin(120)})  node[midway,right,sloped,rotate=90, xshift=-1.5, yshift=-1, scale = 0.6, transform shape]{$\scriptstyle x_o$};
        \end{tikzpicture}
        \caption{$x_o = 0.7$, $x_i = 0.7$}
        \label{DomainHypotrochoid3}
        \end{subfigure}
    \caption{Shape of a hypotrochoid for different parameters $x_0$ and $x_i$}
    \label{HypotrochoidDifferentParameters}
\end{figure}
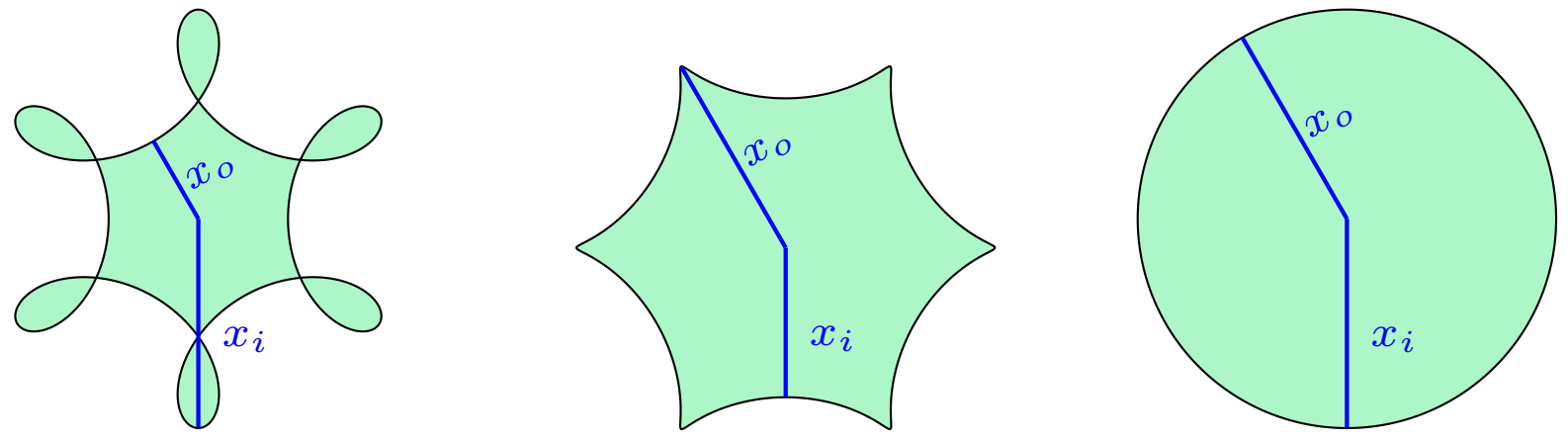 

Furthermore, we can see that the hypotrochoid changes the shape, by varying the parameters of the ellipses. The ratio of the inner and outer radius is limited by $1-\frac{2}{q} \leq \frac{x_i}{x_o} \leq 1$. The lower limit can be determined by examining the geometry of a roll curve with the radius of the fixed circle $R$, the radius of the rolling circle $r$, and the distance of the drawing point to the origin of the rolling curve $d$. The normed distance of the drawing point to the origin of the rolling curve has to be $a = \frac{d}{r} \leq 1$ to ensure that it is not a prolate hypotrochoid with self-intersecting loops. The distance $d$ can be calculated by $d = \frac{x_o - x_i}{2}$, the radius of the inner circle is given by $r = \frac{R}{q}$ and the outer radius can be calculated by $R = \frac{x_o + x_i}{2} + r$. Together we get the limit $\frac{(x_o - x_i) (q-1)}{x_o + x_i} \leq 1$ and with this the lower limit $1-\frac{2}{q} \leq \frac{x_i}{x_o}$. If the ratio would be smaller, the path of the hypotrochoid would get self-intersecting loops, see Figure \ref{HypotrochoidDifferentParameters}(\subref{DomainHypotrochoid1}). In Figure \ref{HypotrochoidDifferentParameters}(\subref{DomainHypotrochoid2}) a typical hypotrochoid with a ratio within the boundaries is displayed. Figure \ref{HypotrochoidDifferentParameters}(\subref{DomainHypotrochoid3}) shows the hypotrochoid with the ratio $\frac{x_i}{x_o} = 1$, which is a circle. If the ratio would be higher it would look again like Figure \ref{HypotrochoidDifferentParameters}(\subref{DomainHypotrochoid2}). The only difference would be the swapping of the names of the radii, $x_i$ and $x_o$.
How theses changes of the shape of the hypotrochoid depend on the parameter set of the ellipses and how this knowledge can help to find other solutions with for example less holes, will be analyzed in Section \ref{sec:BehaviorRatioNodalLine}.
\subsection{Analyzing the parameter sets rejecting Payne's conjecture}
\label{sec:ParametersetsRejectingPayne}
In this section, the parameter sets for which Payne's conjecture is rejected will be described and analyzed numerically.
For six and five elliptical holes within a disc, the parameters $a$, $b$, $r_E$ which define the length, width, and distance of the ellipses to the origin will be varied, and the domains for which Payne's conjecture is rejected will be plotted. From now on, we call these domains $\mathcal{R}_{r_E}$. To achieve this, we previously made a grid search for different length, width, and distances of the ellipses with a step size of $0.1$. Afterwards, we used the Moore-neighbor tracing algorithm \cite[sec. 4]{Pradhan2010_Moore} to get a finer boundary with step size $0.02$.

Some parameter sets are not considered, as they result in overlaps between individual boundaries, which is not permitted by definition. We refer to Theorem \ref{TheoremOverlap}.
\begin{theorem}
    Let there be $q>2$ ellipses of the same shape, each rotated around the origin with the distance $r_E$, the length $a$, and the width $b$, as shown in Figure \ref{DomainIntersectionEllipses}. If the parameters satisfy the condition $b < \sqrt{r_E^2 - \frac{a^2}{\tan^2(\frac{\alpha}{2})}}$, with $\frac{\alpha}{2} = \frac{\pi}{q}$, then there is no overlap between the rotated ellipses.\label{TheoremOverlap}
    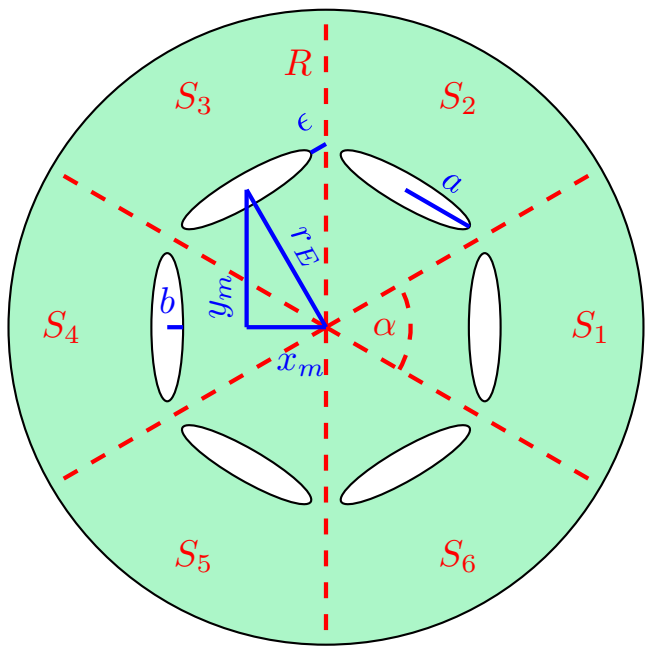
\begin{figure}[h!]
    \centering
        \begin{tikzpicture}[scale=0.7, transform shape]
            \definecolor{pastelgruen}{HTML}{acf6c8}
            \node[circle, draw,fill = pastelgruen, minimum size = 6cm] (p) at (0,0){};
            \foreach \angle in {0,60,...,300}
                \draw[fill = white] (\angle:1.5) ellipse [x radius = 0.7cm, y radius = 0.15cm, rotate= \angle+90];
            \foreach \angle in {30,90,...,330}
                \draw[thick, red, dashed] (0,0) -- ({3*cos(\angle)},{3*sin(\angle)});
            \draw[thick, blue] (0,0) --  ({-1.5*sin(30)},{1.5*cos(30)})  node[midway, above=-1, sloped]{$r_E$};
            \draw[thick, blue] (0,0) --  ({-1.5*sin(30)},0)  node[midway, below=3.5, xshift=4, sloped]{$x_m$};
            \draw[thick, blue] ({-1.5*sin(30)},0) --  ({-1.5*sin(30)},{1.5*cos(30)})  node[midway, above, xshift = -10, sloped]{$y_m$};
            \draw[thick, blue] (-1.5,0) node[above, sloped]{$b$} --  (-1.35,0)  ;
            \draw[thick, blue] (60:1.5) --  ({1.5*cos(60)+0.7*cos(330)},{1.5*sin(60)+0.7*sin(330)})  node[midway, above=2, sloped]{$a$};
            \draw[thick, blue] ({-1.5*sin(30)+0.7*cos(30)},{1.5*cos(30)+0.7*sin(30)}) --  ({-1.5*sin(30)+(tan(30)*1.5)*cos(30)},{1.5*cos(30)+(tan(30)*1.5)*sin(30)})  node[midway, above=2, sloped]{$\epsilon$};
            \draw[thick, red, dashed] (-30:0.8) arc(-30:30:0.8) node[midway, left]{$\alpha$};
            \node[left, red, thick] at (0,2.5){$R$};
            \node[red, thick] at (0:2.5){$S_1$};
            \node[red, thick] at (60:2.5){$S_2$};
            \node[red, thick] at (120:2.5){$S_3$};
            \node[red, thick] at (180:2.5){$S_4$};
            \node[red, thick] at (240:2.5){$S_5$};
            \node[red, thick] at (300:2.5){$S_6$};
        \end{tikzpicture}
    \caption{Domain is a circle with radius $R$ centered at the origin with $q= 6$ ellipses, which must lie within the subsections $S_i$, with $i= 1,\dots,6$}.
    \label{DomainIntersectionEllipses}
\end{figure}
\end{theorem}
\begin{proof}
Each ellipse must be a subspace of the circular sector $S_i \coloneq \{ (r, \phi) \;| \;0 \leq r \leq R,\; \phi_i \leq \phi \leq \phi_i + \alpha\}$ with $i = 1, \dots, q$ and $\phi_i = \frac{(2i-3) \pi}{q}$ to avoid overlaps between the rotated ellipses.
To investigate for which width $b$ of an ellipse with a given length $a$ and distance $r_E$ the ellipse intersects the boundary of the circular sector, we take a deeper look into the third circular sector $S_3$ in Figure \ref{DomainIntersectionEllipses}. Here, the right boundary is on the ordinate. In cases with a different number of circles, where this is not the case, we can just rotate the circle with the ellipses, without changing something of the result. Now, the parameterization of the ellipse centered at $\left(x_m, y_m\right)$ is given by 
\begin{align}
    \begin{pmatrix}
        x_E\\y_E
    \end{pmatrix}
    =
    \begin{pmatrix}
        -x_m + a \cos\left(t\right)\cos\left(\frac{\alpha}{2}\right) - b \sin\left(t\right)\sin\left(\frac{\alpha}{2}\right)\\
        y_m + a \cos\left(t\right)\sin\left(\frac{\alpha}{2}\right) + b \sin\left(t\right)\cos\left(\frac{\alpha}{2}\right)
    \end{pmatrix},
    \qquad t \in [0,2\pi)
\end{align}
and it is clear, that the first coordinate of the ellipse $x_E$ has to be smaller than zero, otherwise it intersects the boundary of the circular sector
\begin{align}
    -r_E \sin\left(\frac{\alpha}{2}\right) + a \cos\left(t\right)\cos\left(\frac{\alpha}{2}\right) - b \sin\left(t\right)\sin\left(\frac{\alpha}{2}\right) \stackrel{!}{<} 0 \,,
\end{align} with $x_m = r_E \sin\left(\frac{\alpha}{2}\right)$.
In the next step, we have to expand for $b$
\begin{align}
      b &< \frac{a}{\tan\left(\frac{\alpha}{2}\right)\tan\left(t\right)} - \frac{r_E}{\sin\left(t\right)}.
      \label{UpperBoundOfB}
\end{align}
Hence, we have to calculate the minimum
\begin{equation}
    \min_{t \in \left(\frac{3\pi}{2},\; 2\pi\right)}\left( \frac{a}{\tan\left(\frac{\alpha}{2}\right)\tan\left(t\right)} - \frac{r_E}{\sin\left(t\right)}\right).
\end{equation}
This we can do by calculating the derivative with respect to $t$ and set the result equal to zero
\begin{equation}
    \begin{aligned}
        &\frac{-a\sin^2\left(t\right)\tan\left(\frac{\alpha}{2}\right) + r_E \cos\left(t\right)\tan^2\left(\frac{\alpha}{2}\right) - a \cos^2\left(t\right)\tan\left(\frac{\alpha}{2}\right)}{\sin^2\left(t\right)\tan^2\left(\frac{\alpha}{2}\right)} = 0\\
        &\Leftrightarrow \qquad t_{1,2}= \pm \arccos\left(\frac{a}{r_E \tan\left(\frac{\alpha}{2}\right)}\right),
    \end{aligned}
\end{equation}
where $t_{1}= \arccos\left(\frac{a}{r_E \tan\left(\frac{\alpha}{2}\right)}\right)$ is the maximum and $t_{2}= - \arccos\left(\frac{a}{r_E \tan\left(\frac{\alpha}{2}\right)}\right)$ is the minimum. This can be easily verified by checking the second derivative with respect to $t$ and inserting afterwards first $t_1$and then $t_2$.
By inserting this result in \eqref{UpperBoundOfB} we can calculate that
\begin{equation}
    \begin{aligned}
        b &< \sqrt{r_E^2  - \frac{a^2}{\tan^2\left(\frac{\alpha}{2}\right)}}.
    \end{aligned}
\end{equation}
\end{proof}
The parameter sets that are not considered are indicated by hatched areas in the resulting plots for five holes in Figure \ref{fig:Parametersets_5} and for six holes in Figure \ref{fig:Parametersets_6}.
\begin{figure}[h!]
    \centering
    \begin{subfigure}{0.32\textwidth}
        \includegraphics[width=1\linewidth]{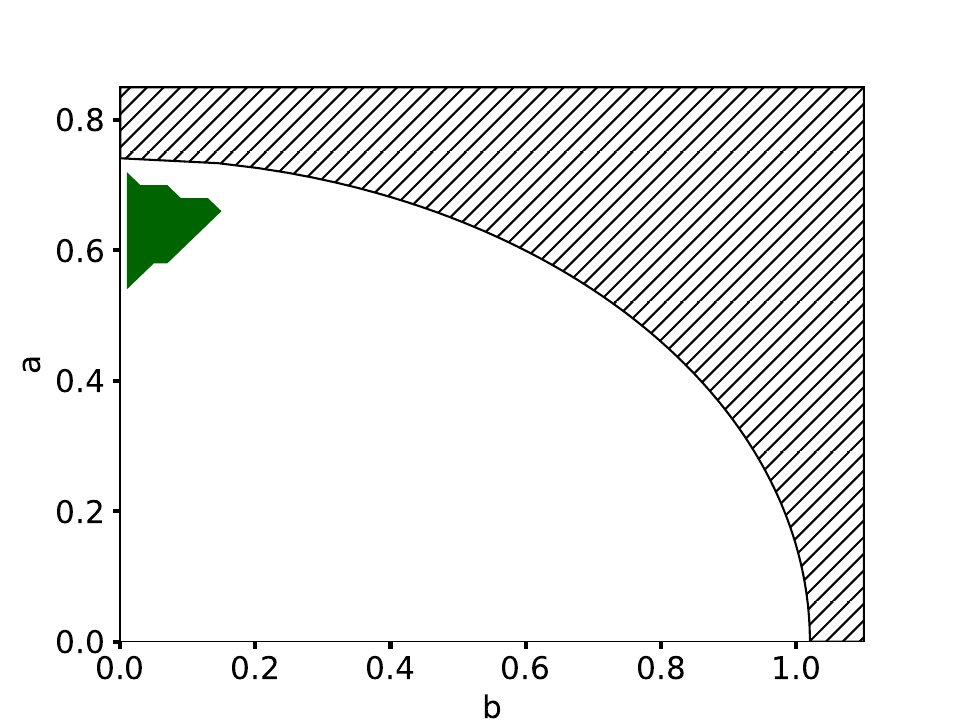}
        \caption{$r_E = 1.02$}
        \label{fig:Parametersets_5_dist_02}
    \end{subfigure}
    \begin{subfigure}{0.32\textwidth}
        \includegraphics[width=1\linewidth]{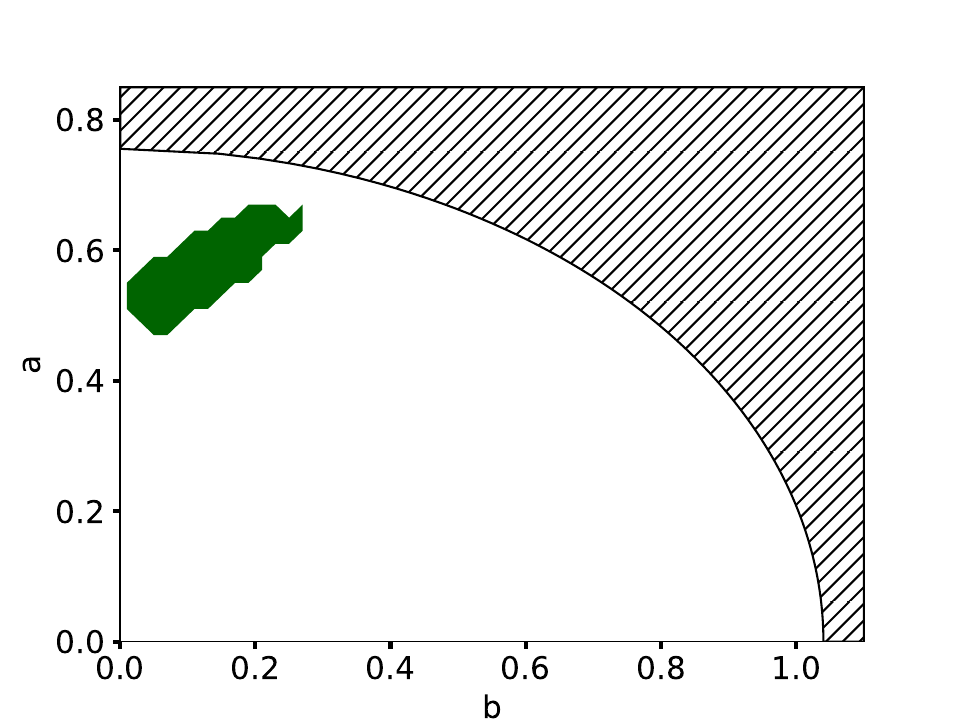}
        \caption{$r_E = 1.04$}
        \label{fig:Parametersets_5_dist_04}
    \end{subfigure}
    \begin{subfigure}{0.32\textwidth}
        \includegraphics[width=1\linewidth]{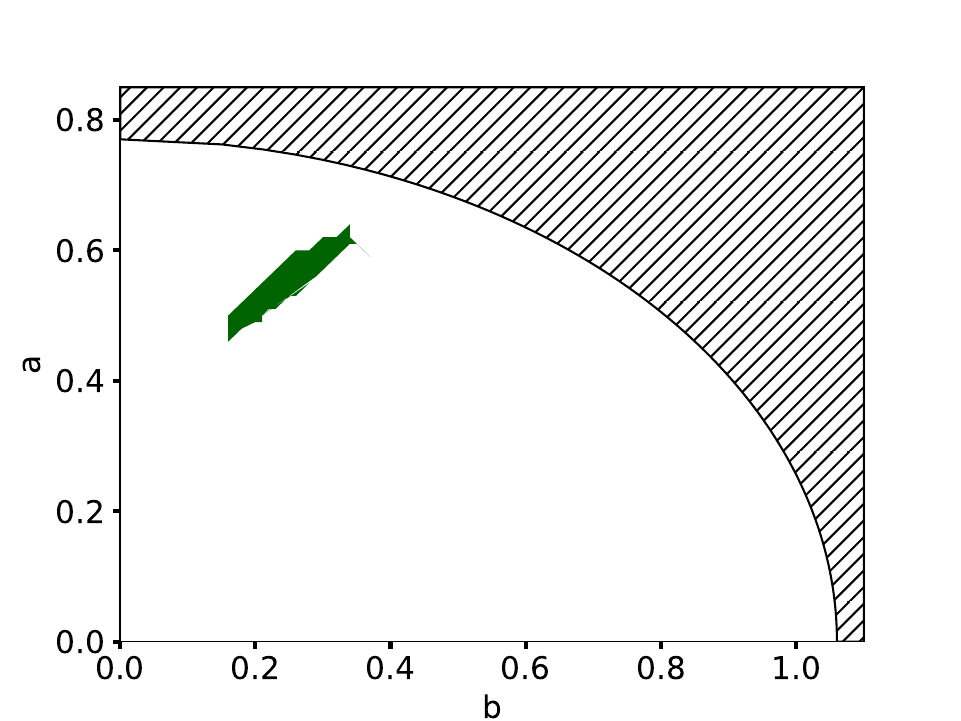}
        \caption{$r_E = 1.06$}
        \label{fig:Parametersets_5_dist_06}
    \end{subfigure}
    \begin{subfigure}{0.32\textwidth}
        \includegraphics[width=1\linewidth]{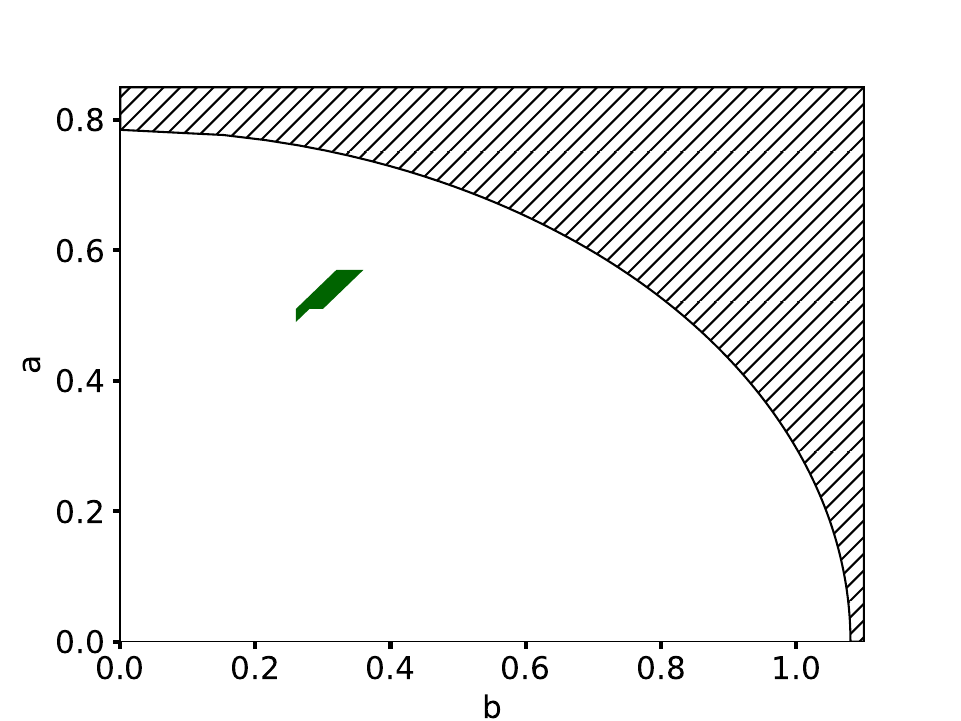}
        \caption{$r_E = 1.08$}
        \label{fig:Parametersets_5_dist_08}
    \end{subfigure}
    \caption{Parameter sets (green) rejecting Payne's conjecture for $5$ holes, with different distances of the ellipses to the origin $r_E$. Parameter sets (black hatched) are not considered due to Theorem \ref{TheoremOverlap}.}
    \label{fig:Parametersets_5}
\end{figure}

\begin{figure}[h!]
    \centering
    \begin{subfigure}{0.32\textwidth}
        \includegraphics[width=1\linewidth]{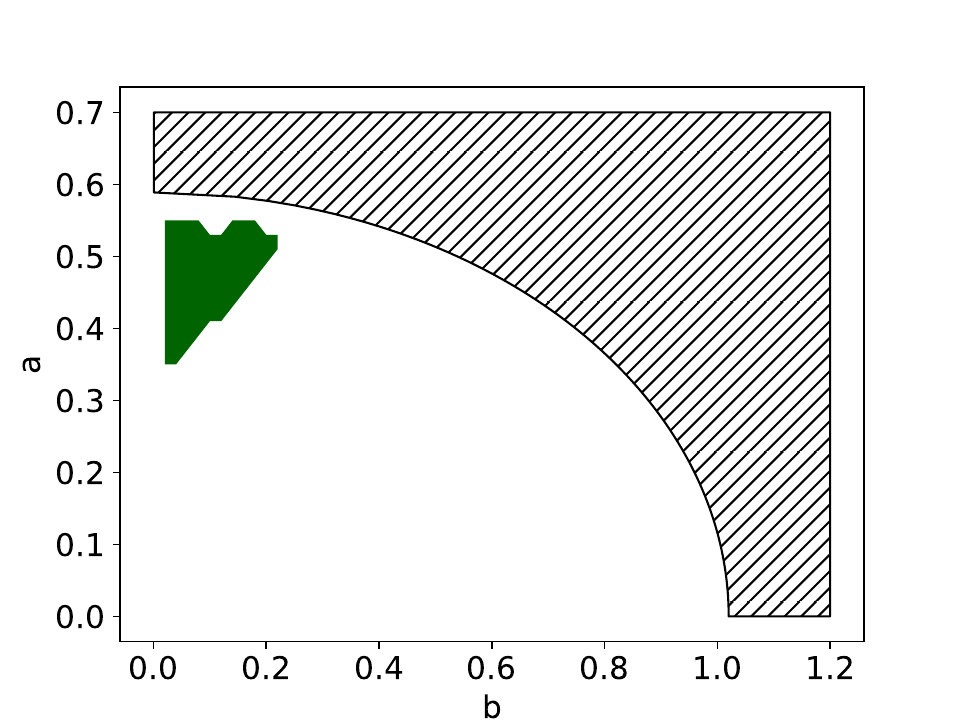}
        \caption{$r_E = 1.02$}
        \label{fig:Parametersets_6_dist_02}
    \end{subfigure}
    \begin{subfigure}{0.32\textwidth}
        \includegraphics[width=1\linewidth]{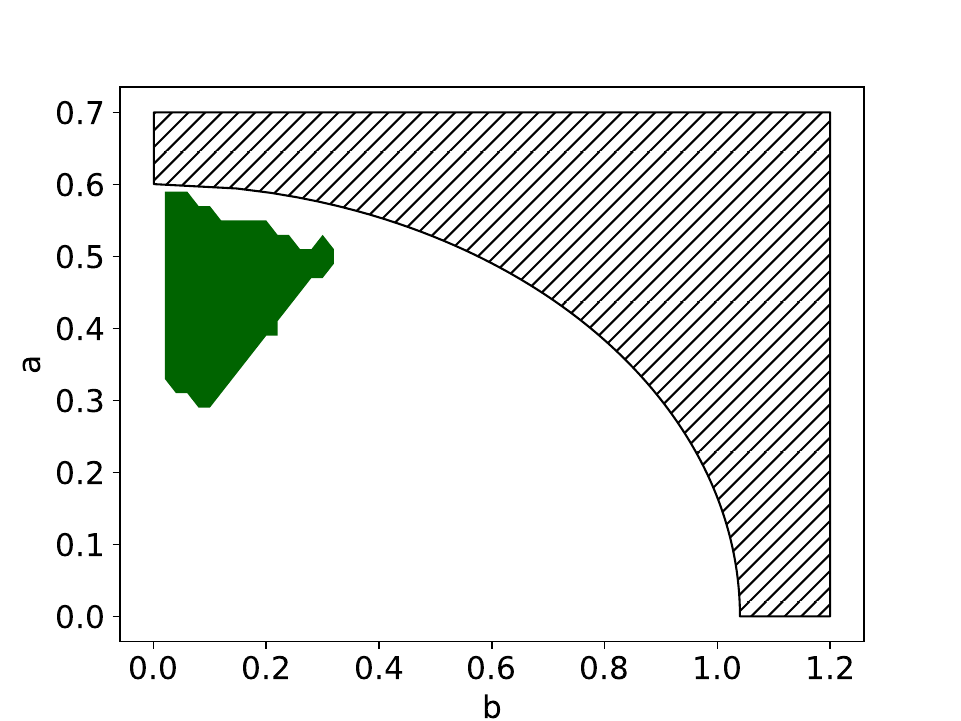}
        \caption{$r_E = 1.04$}
        \label{fig:Parametersets_6_dist_04}
    \end{subfigure}
    \begin{subfigure}{0.32\textwidth}
        \includegraphics[width=1\linewidth]{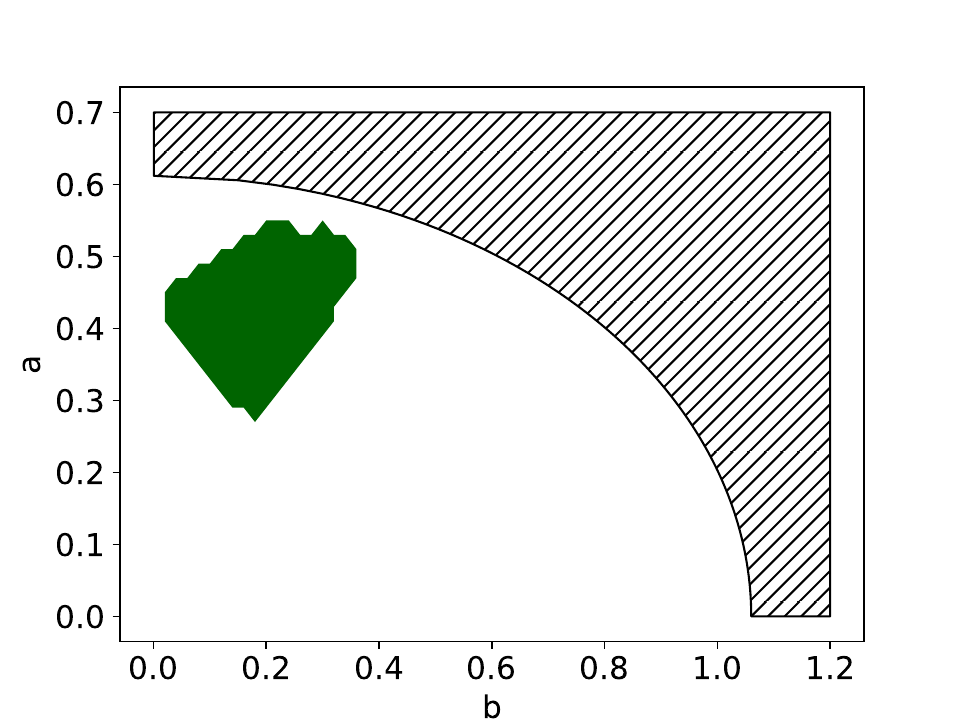}
        \caption{$r_E = 1.06$}
        \label{fig:Parametersets_6_dist_06}
    \end{subfigure}
    \begin{subfigure}{0.32\textwidth}
        \includegraphics[width=1\linewidth]{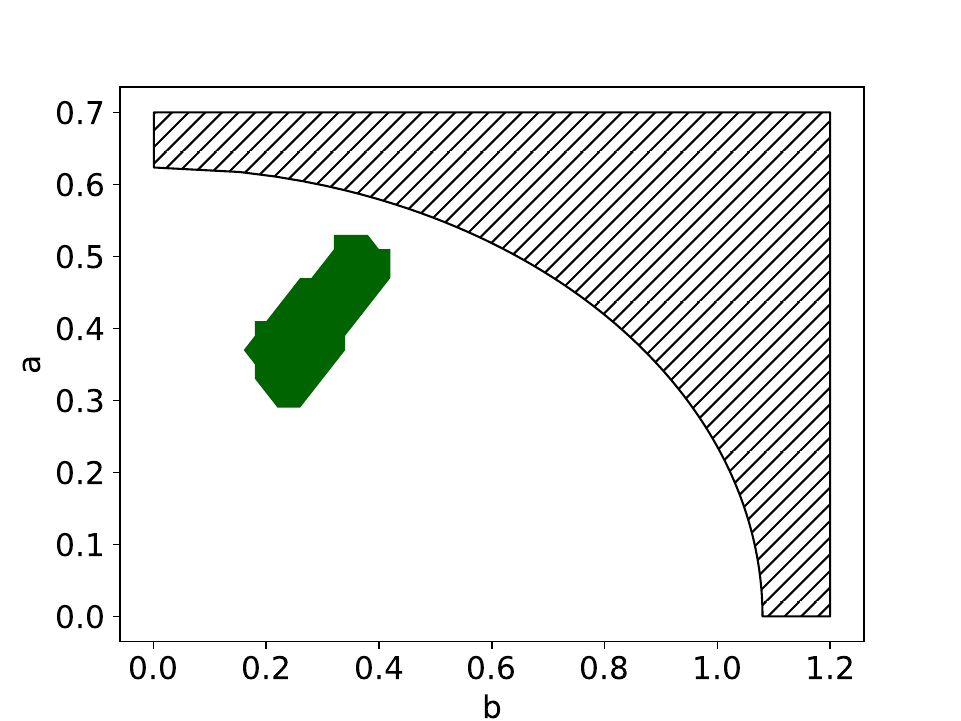}
        \caption{$r_E = 1.08$}
        \label{fig:Parametersets_6_dist_08}
    \end{subfigure}
    \begin{subfigure}{0.32\textwidth}
        \includegraphics[width=1\linewidth]{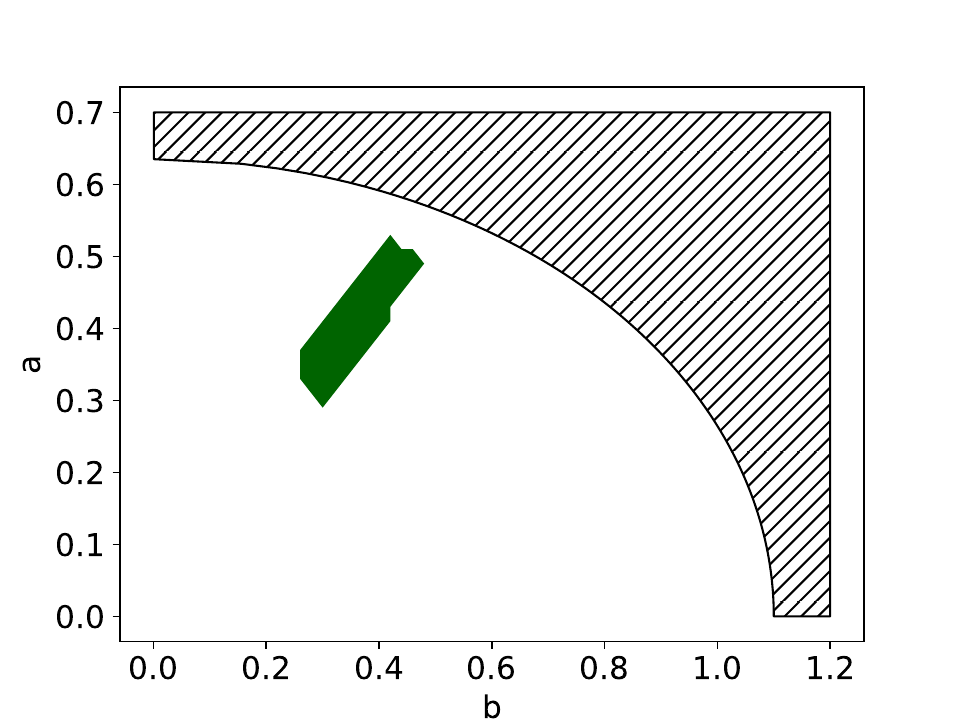}
        \caption{$r_E = 1.10$}
        \label{fig:Parametersets_6_dist_10}
    \end{subfigure}
    \begin{subfigure}{0.32\textwidth}
        \includegraphics[width=1\linewidth]{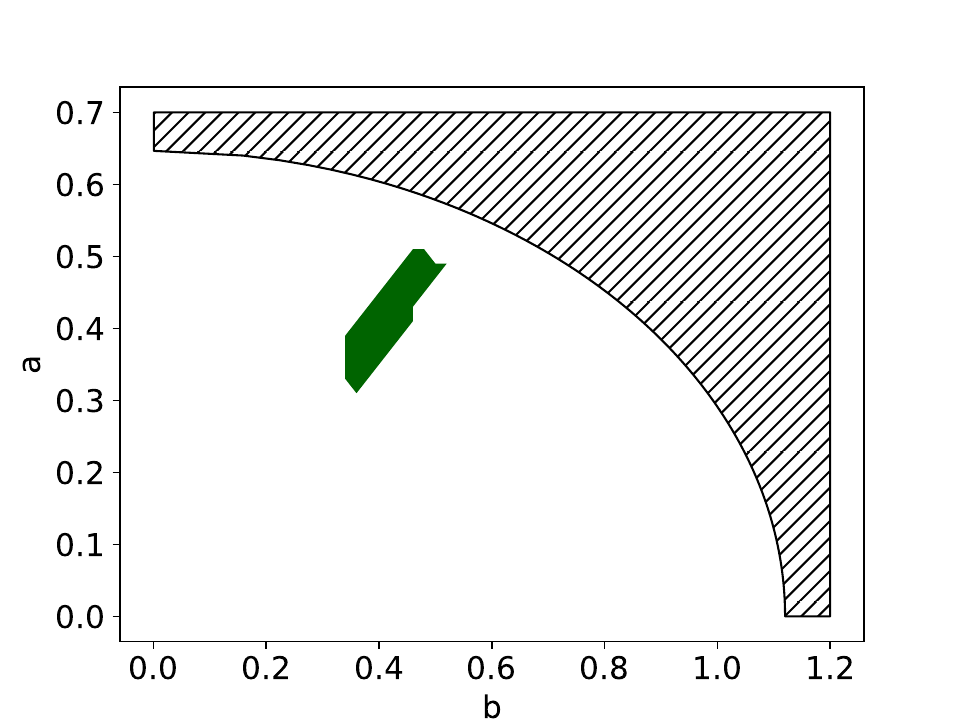}
        \caption{$r_E = 1.12$}
        \label{fig:Parametersets_6_dist_12}
    \end{subfigure}
    \begin{subfigure}{0.32\textwidth}
        \includegraphics[width=1\linewidth]{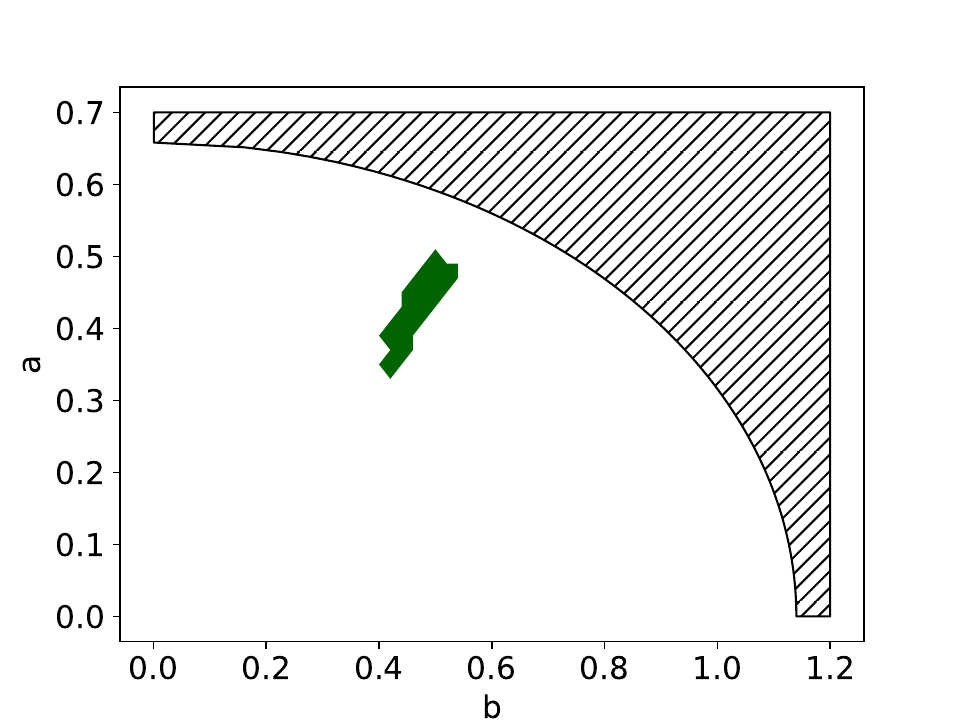}
        \caption{$r_E = 1.14$}
        \label{fig:Parametersets_6_dist_14}
    \end{subfigure}
    \caption{Parameter sets (green) rejecting Payne's conjecture for $6$ holes, with different distances of the ellipses to the origin $r_E$. Parameter sets (black hatched) are not considered due to Theorem \ref{TheoremOverlap}.}
    \label{fig:Parametersets_6}
\end{figure}

Based on these observations, we conjecture the following hypotheses. Hypotheses (3)-(6) consider domains with the same overall geometry. Only the number of holes and their individual parameters are varied. All unspecified parameters are fixed.
\begin{enumerate}
    \item If a nodal line rejecting Payne's conjecture in two dimensions exists, then it always has the shape of a hypotrochoid or an affine transformation of a hypotrochoid.
    \item For a given parameterization centered with a fixed distance around the origin, there exist just a one simply-connected parameter set, rejecting Payne's conjecture.
    \item The fewer holes, the smaller is the range of distance parameters for which Payne's conjecture is rejected.
    \item The fewer holes, the smaller is the volume of the parameter set $(a,b,r_E)$ for which Payne's conjecture is rejected.
    \item If there exist parameters for a specific distance, such that Payne's conjecture is rejected, then the bigger the distance of the holes to the origin, the bigger is the width of the holes in average.
    \item If there exist parameters for a specific distance, such that Payne's conjecture is rejected, then the bigger the distance of the holes to the origin, the smaller is the parameter set which rejects Payne's conjecture, if the region does not intersect the boundary of the possible area.
\end{enumerate}
\subsection{Behavior of the shape of the nodal line}
\label{sec:BehaviorRatioNodalLine}
The questions, we want to answer in this section are: How does the shape of the nodal line change by varying the parameters of the ellipses $r_E$, $a$, $b$, and $q$ and what is the maximum of the minimum pairwise distances between the nodal line and the boundary of the holes for five and six ellipses. Furthermore, we want to find out whether it is possible to get a knowledge about the parameter domains, before calculating them. This might help finding different domains with a different number of holes rejecting Payne's conjecture, with a minimum amount of computational time.
\begin{figure}[h!]
    \centering
    \begin{subfigure}[t]{0.24\textwidth}
        \includegraphics[width=1\linewidth]{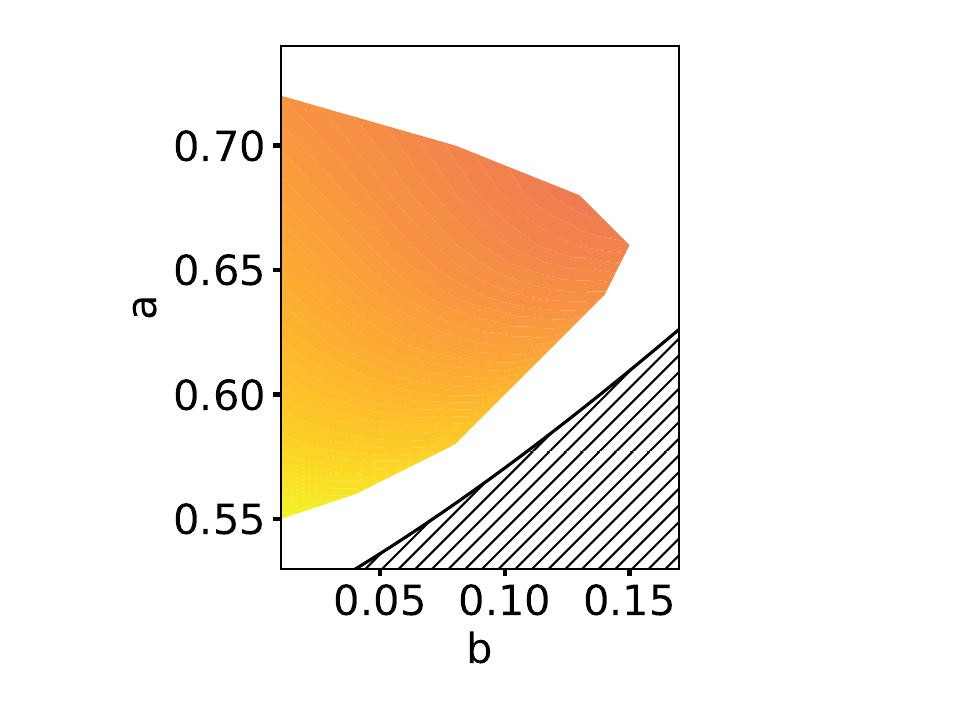}
        \caption{$q = 5$; $r_E = 1.02$}
        \label{fig:Xa_5_dist_02}
    \end{subfigure}
    \begin{subfigure}[t]{0.24\textwidth}
        \includegraphics[width=1\linewidth]{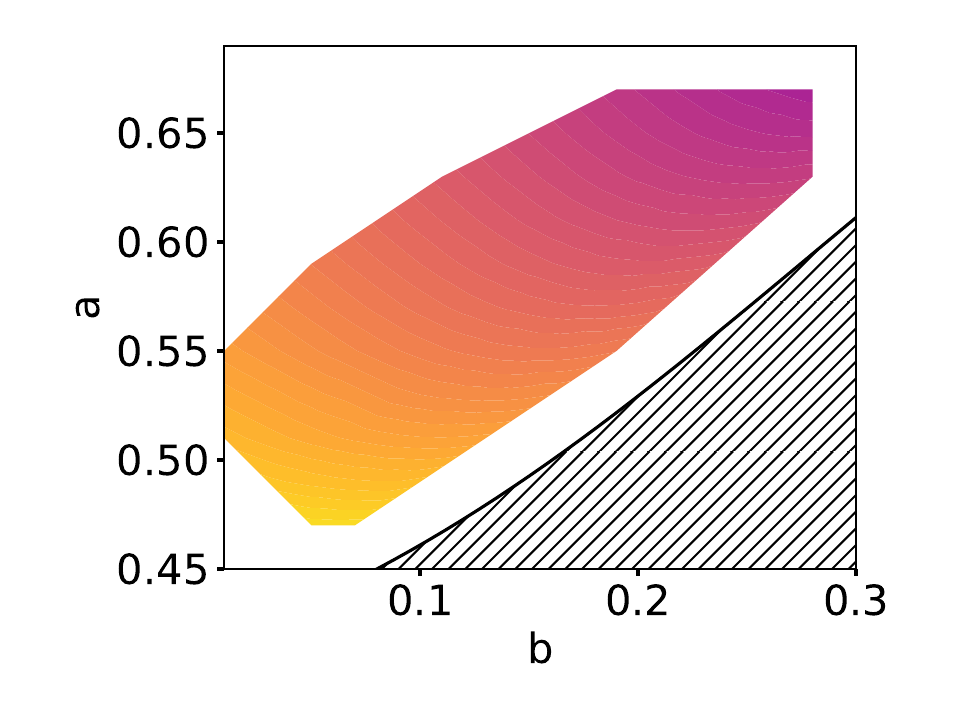}
        \caption{$q=5$; $r_E = 1.04$}
        \label{fig:Xa_5_dist_04}
    \end{subfigure}
    \begin{subfigure}[t]{0.24\textwidth}
        \includegraphics[width=1\linewidth]{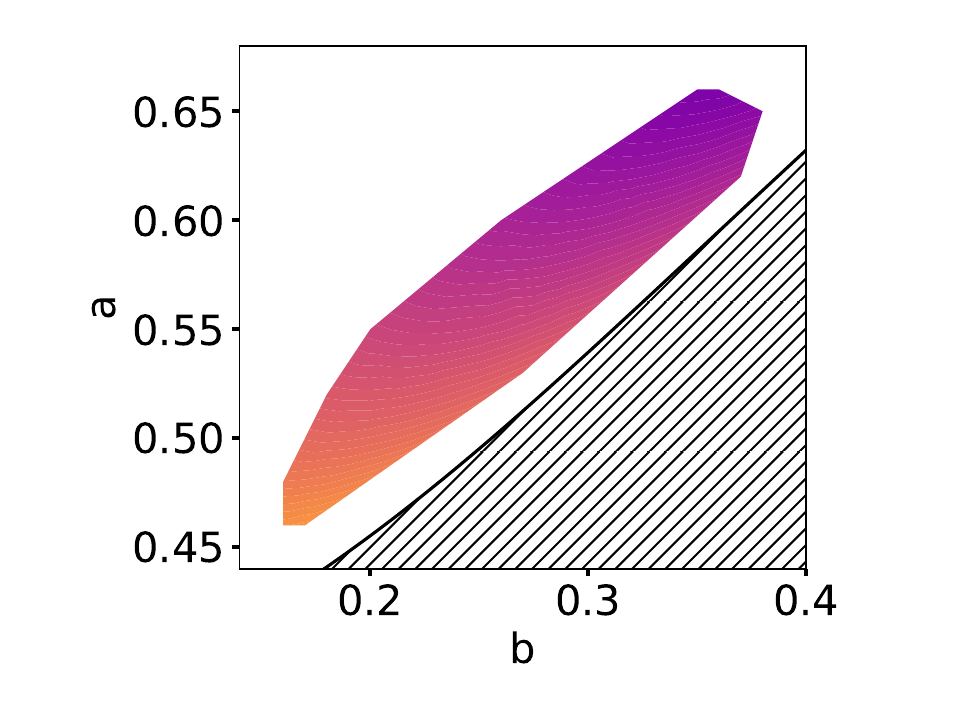}
        \caption{$q=5$; $r_E = 1.06$}
        \label{fig:Xa_5_dist_06}
    \end{subfigure}
    \begin{subfigure}[t]{0.24\textwidth}
        \includegraphics[width=1\linewidth]{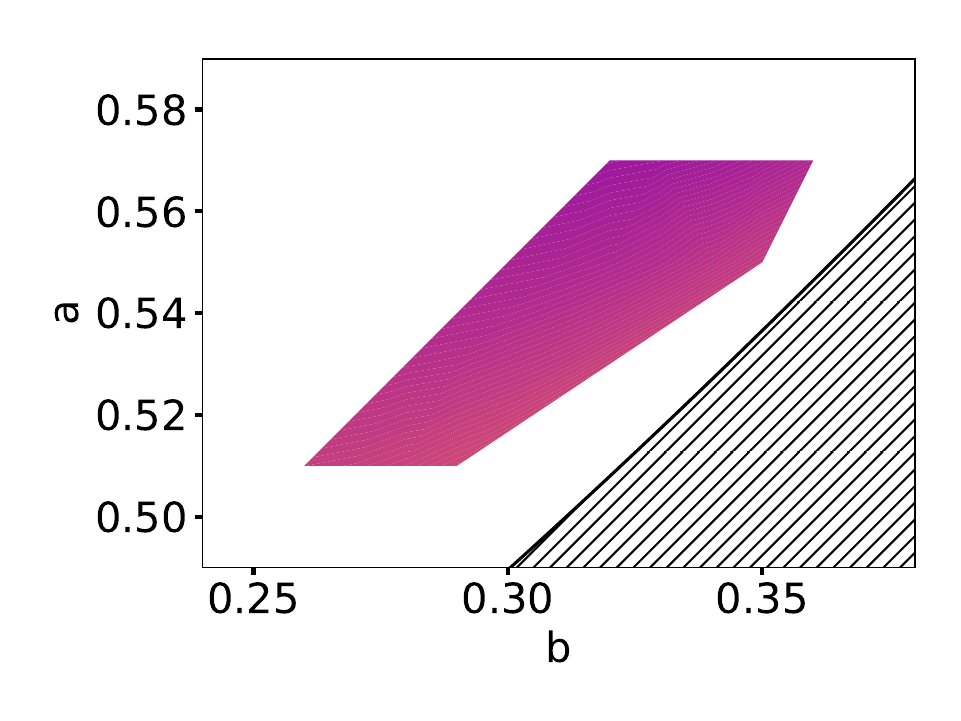}
        \caption{$q=5$; $r_E = 1.08$}
        \label{fig:Xa_5_dist_08}
    \end{subfigure}
    \begin{subfigure}[t]{0.24\textwidth}
        \includegraphics[width=1\linewidth]{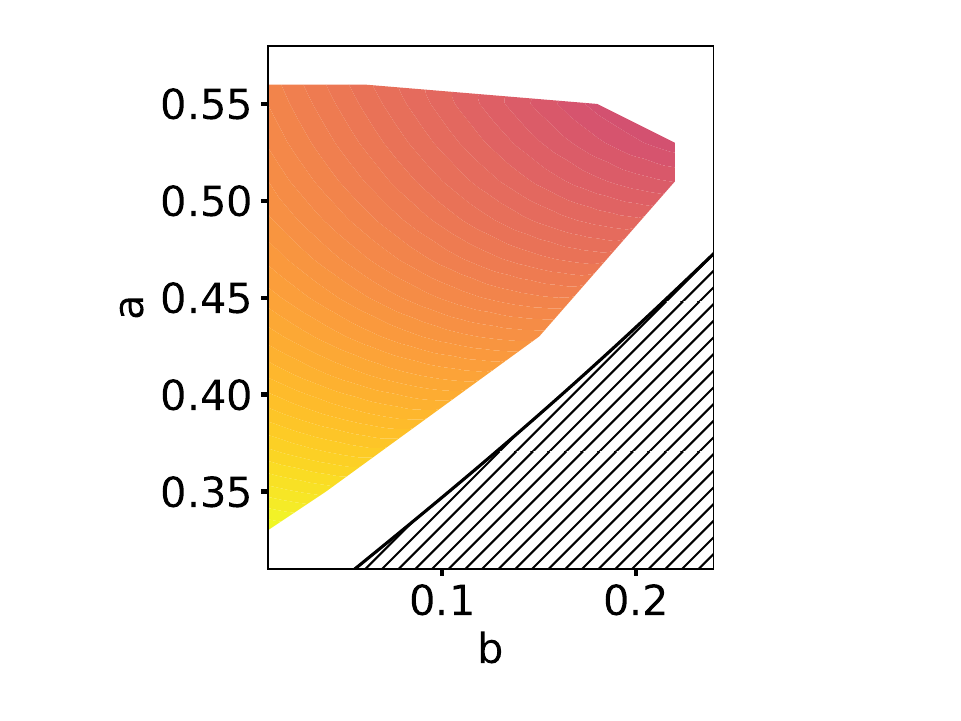}
        \caption{$q=6$; $r_E = 1.02$}
        \label{fig:Xa_6_dist_02}
    \end{subfigure}
    \begin{subfigure}[t]{0.24\textwidth}
        \includegraphics[width=1\linewidth]{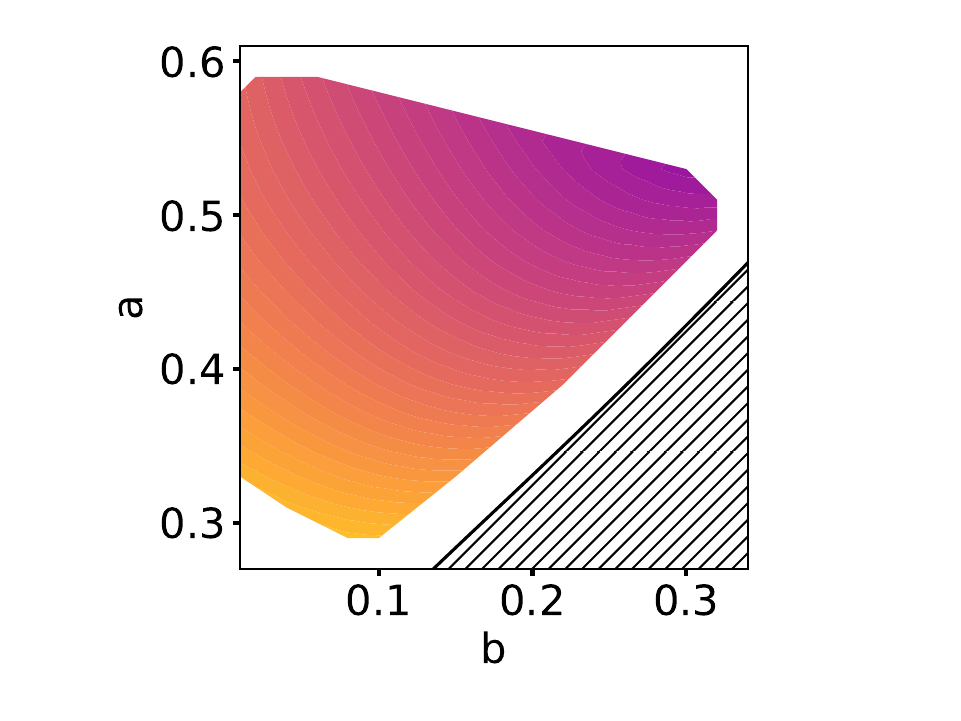}
        \caption{$q=6$; $r_E = 1.04$}
        \label{fig:Xa_6_dist_04}
    \end{subfigure}
    \begin{subfigure}[t]{0.24\textwidth}
        \includegraphics[width=1\linewidth]{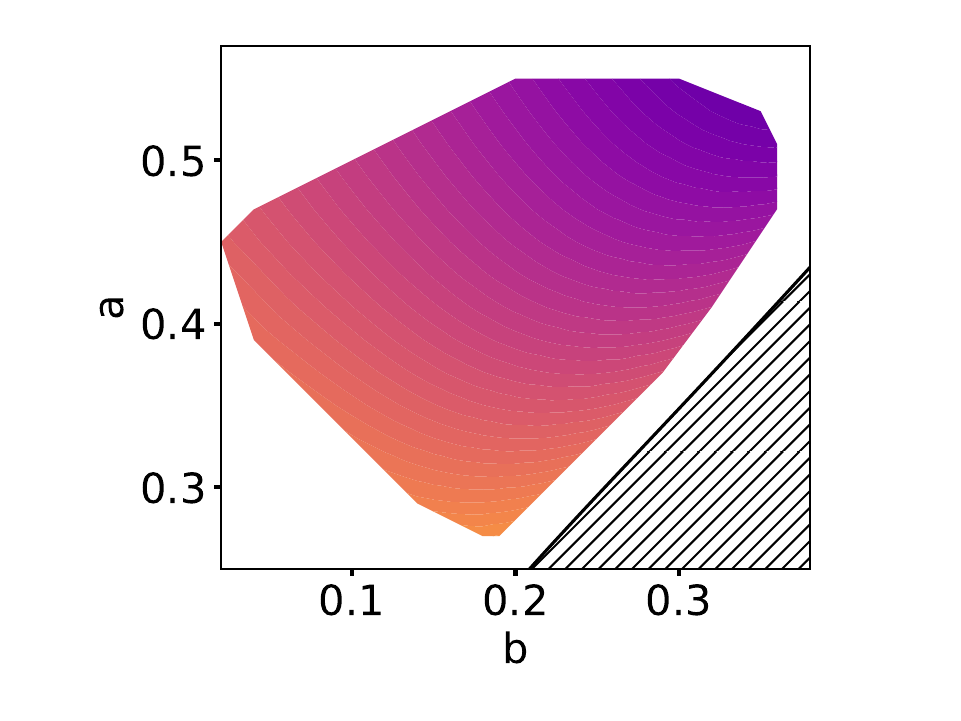}
        \caption{$q=6$; $r_E = 1.06$}
        \label{fig:Xa_6_dist_06}
    \end{subfigure}
    \begin{subfigure}[t]{0.24\textwidth}
        \includegraphics[width=1\linewidth]{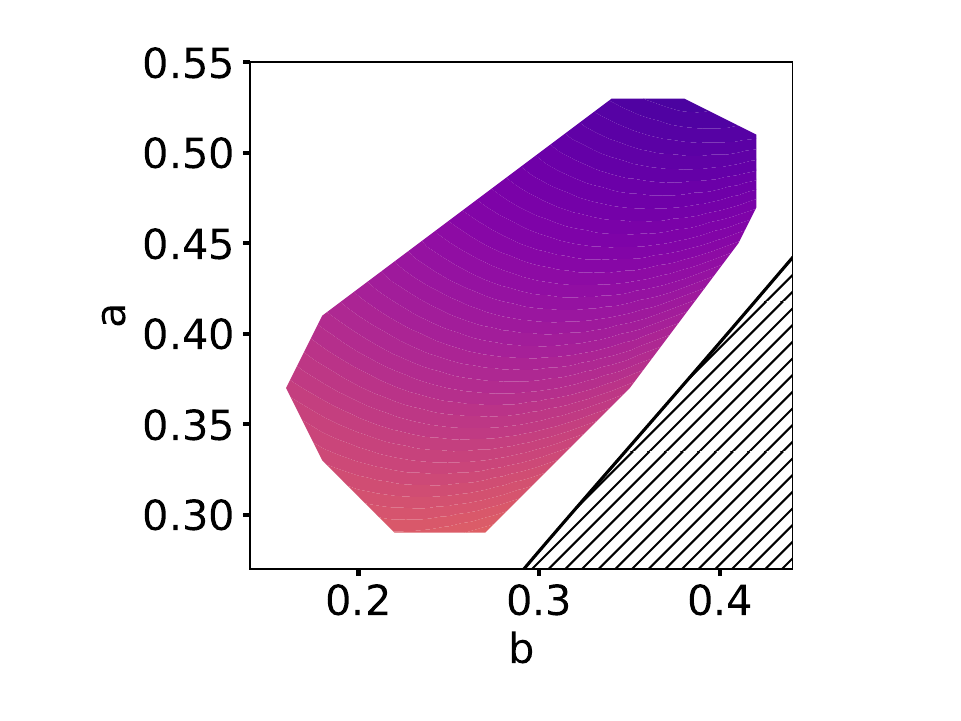}
        \caption{$q=6$; $r_E = 1.08$}
        \label{fig:Xa_6_dist_08}
    \end{subfigure}
    \begin{subfigure}[t]{0.24\textwidth}
        \includegraphics[width=1\linewidth]{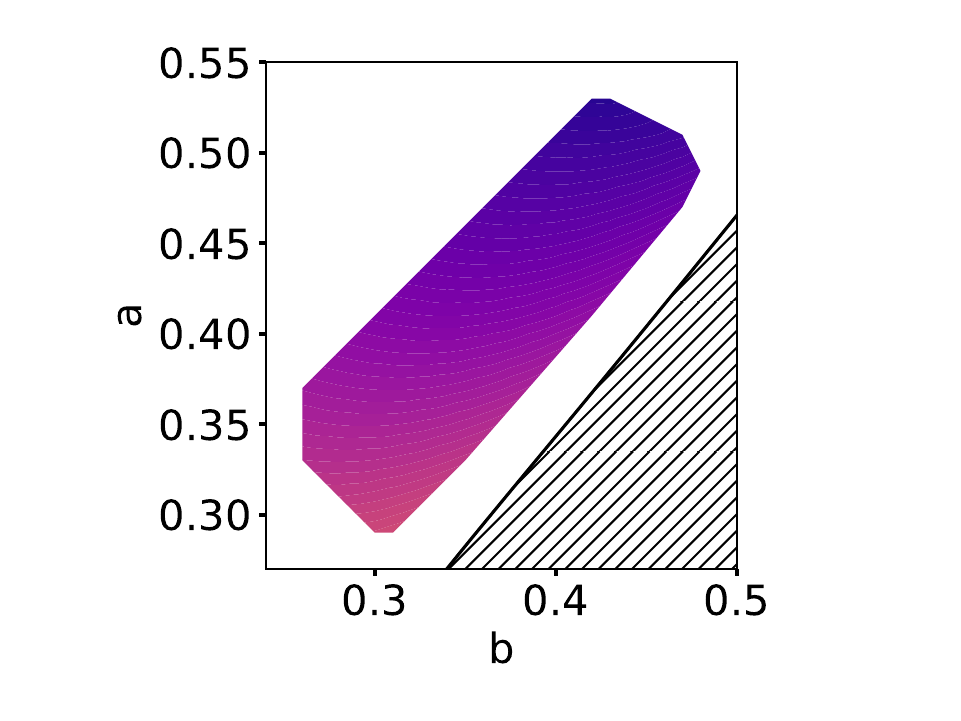}
        \caption{$q=6$; $r_E = 1.1$}
        \label{fig:Xa_6_dist_10}
    \end{subfigure}
    \begin{subfigure}[t]{0.24\textwidth}
        \includegraphics[width=1\linewidth]{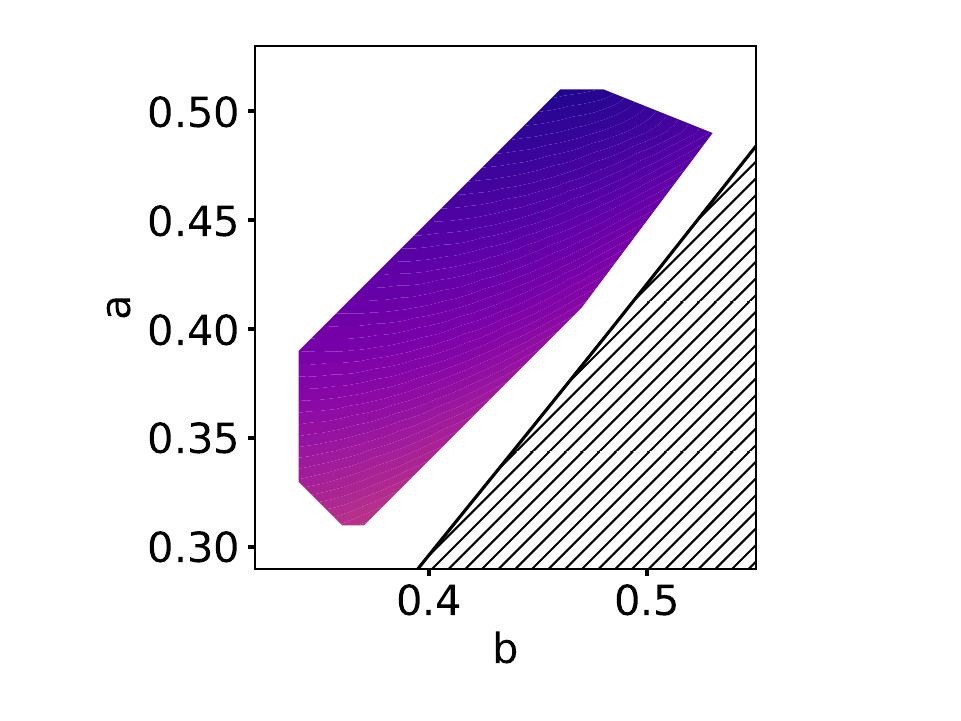}
        \caption{$q=6$; $r_E = 1.12$}
        \label{fig:Xa_6_dist_12}
    \end{subfigure}
    \begin{subfigure}[t]{0.24\textwidth}
        \includegraphics[width=1\linewidth]{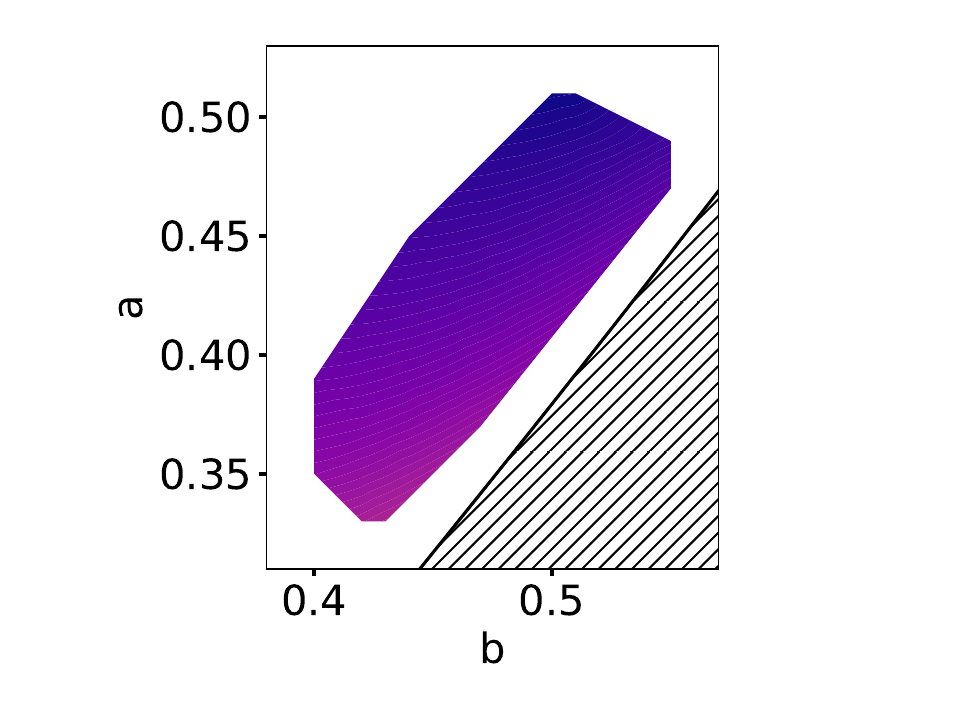}
        \caption{$q=6$; $r_E = 1.14$}
        \label{fig:Xa_6_dist_14}
    \end{subfigure}
    \begin{subfigure}[t]{0.24\textwidth}
        \centering
        \includegraphics[height=\heightof{\includegraphics[width=1\linewidth]{Xa_6_ellipses_distance_1_14.pdf}}]{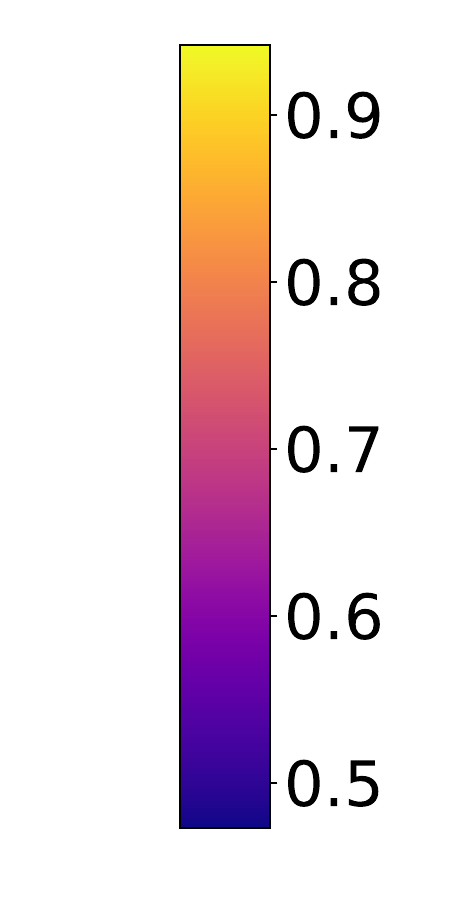}
        \caption{Colorbar}
        \label{fig:Xa_colorbar}
    \end{subfigure}
    \caption{Outer radius of the hypotrochoid based on model \eqref{ModelFitHypotrochoid} for a different amount of holes $q$ and distances of the holes to the origin $r_E$. The colorbar for these plots is visualized in (\subref{fig:Xa_colorbar})}
    \label{fig:Xa_model}
\end{figure}

To answer these questions, we use the calculated nodal lines generated in Section \ref{sec:ParametersetsRejectingPayne} utilized to find the parameter domains rejecting Payne's conjecture with six holes and five holes, see Figures \ref{fig:Parametersets_6} and \ref{fig:Parametersets_5}, and the knowledge from the previous Section \ref{sec:ParametersetsRejectingPayne}, that the nodal line is always an affine transformation of an hypotrochoid. First, we fit a hypotrochoid against the nodal line by using the least squares method with different starting points. With these results, we can find the outer radius $x_o$, the inner radius $x_i$, and the ratio $\frac{x_i}{x_o}$ of the hypotrochoid for each nodal line for a specific parameter set $r_E$, $a$, $b$ and $q$.
Next, we fit for each single parameter $x_i$, $x_o$, and $\frac{x_i}{x_o}$ a model against these results, which is quartic for the parameters $a$, $b$, and $r_E$, but linear for the number of holes $q$, given by 
\begin{equation}
\begin{aligned}
    M_1(r_E, a, b) =& \left(r_E^4, \, a^4,\, b^4,\, r_E^3 a,\, r_E^3 b, \,a^3 r_E, \,a^3 b, \,b^3 r_E,\, b^3 a, \,r_E^2 a^2,\,r_E^2 a b, \right.\\
    &r_E^2 b^2, \, a^2 r_E b, \,a^2 b^2,\, b^2 r_E a, \,r_E^3, \,a^3, \,b^3, \,r_E^2 a, \,r_E^2 b,\, a^2 r_E,\, a^2 b,  \\
    &\left.  b^2 r_E,\, b^2 a, \,r_E a b,\, r_E^2, \,a^2, \,b^2,\, r_E a,\, r_E b,\, a b, \,r_E, \,a,\, b, \,1\right),
    \label{SubModelFitHypotrochoid}
\end{aligned}
\end{equation}
\begin{equation}
    M(r_E, a, b, q) = \left(a_0, \dots , a_{34} \right) M_1(r_E, a, b)^{\top} q + \left(a_{35}, \dots , a_{69} \right) M_1(r_E, a, b)^{\top}.
    \label{ModelFitHypotrochoid}
\end{equation}
The corresponding parameters for the three models can be found in the GitHub repository \href{https://github.com/mattfeldt/nodal-lines}{\texttt{https://github.com/mattfeldt/nodal\_lines}} in the tabular in file \href{https://github.com/mattfeldt/nodal-lines/src/parameters_models_xa_xi_ratio.csv}{\texttt{src/parameters\_models\_xa\_xi\_ratio.csv}}.
\begin{figure}[h!]
    \centering
    \begin{subfigure}[t]{0.24\textwidth}
        \includegraphics[width=1\linewidth]{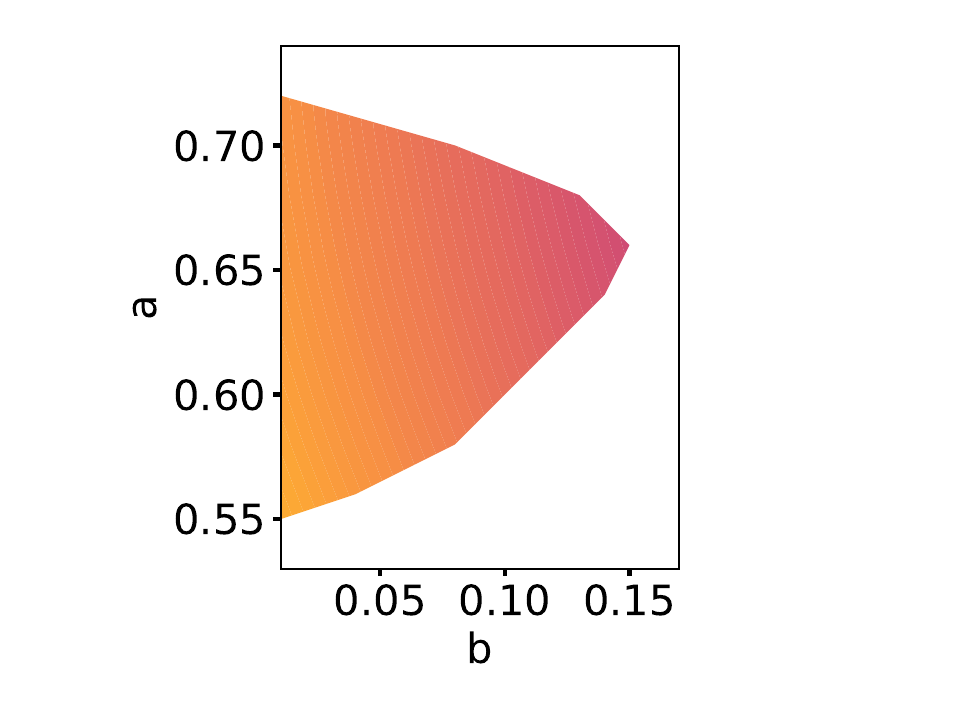}
        \caption{$q=5$; $r_E = 1.02$}
        \label{fig:Xi_5_dist_02}
    \end{subfigure}
    \begin{subfigure}[t]{0.24\textwidth}
        \includegraphics[width=1\linewidth]{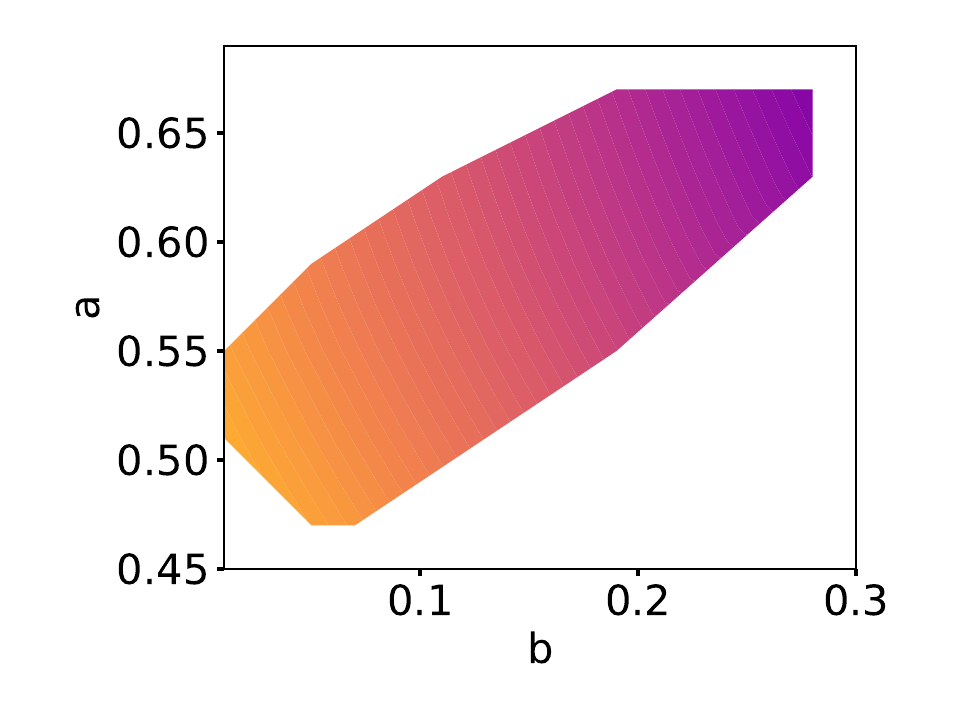}
        \caption{$q=5$; $r_E = 1.04$}
        \label{fig:Xi_5_dist_04}
    \end{subfigure}
    \begin{subfigure}[t]{0.24\textwidth}
        \includegraphics[width=1\linewidth]{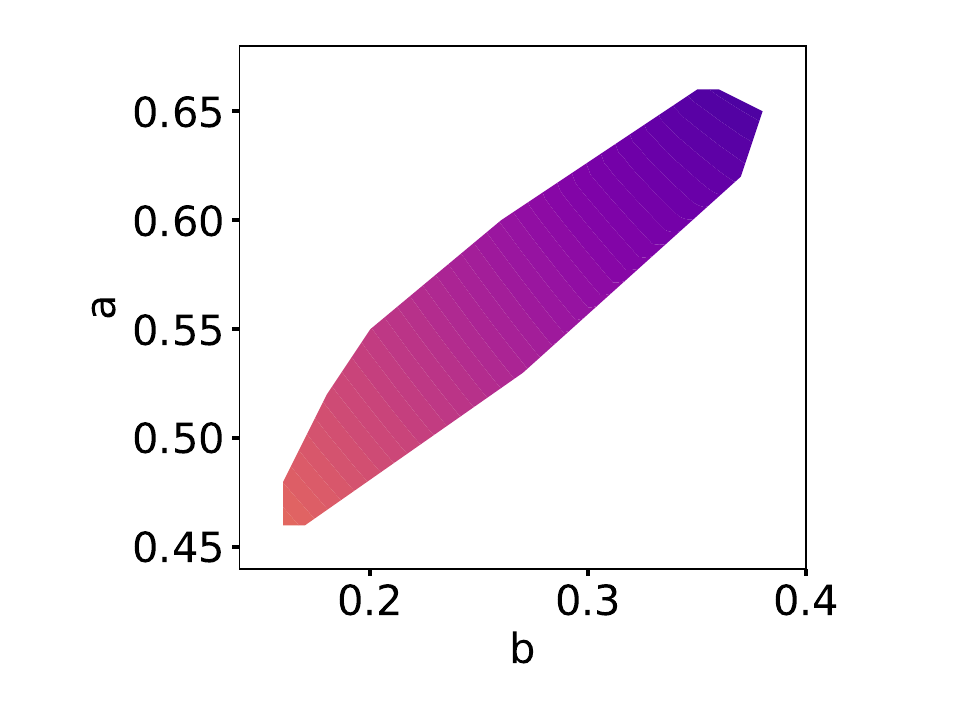}
        \caption{$q=5$; $r_E = 1.06$}
        \label{fig:Xi_5_dist_06}
    \end{subfigure}
    \begin{subfigure}[t]{0.24\textwidth}
        \includegraphics[width=1\linewidth]{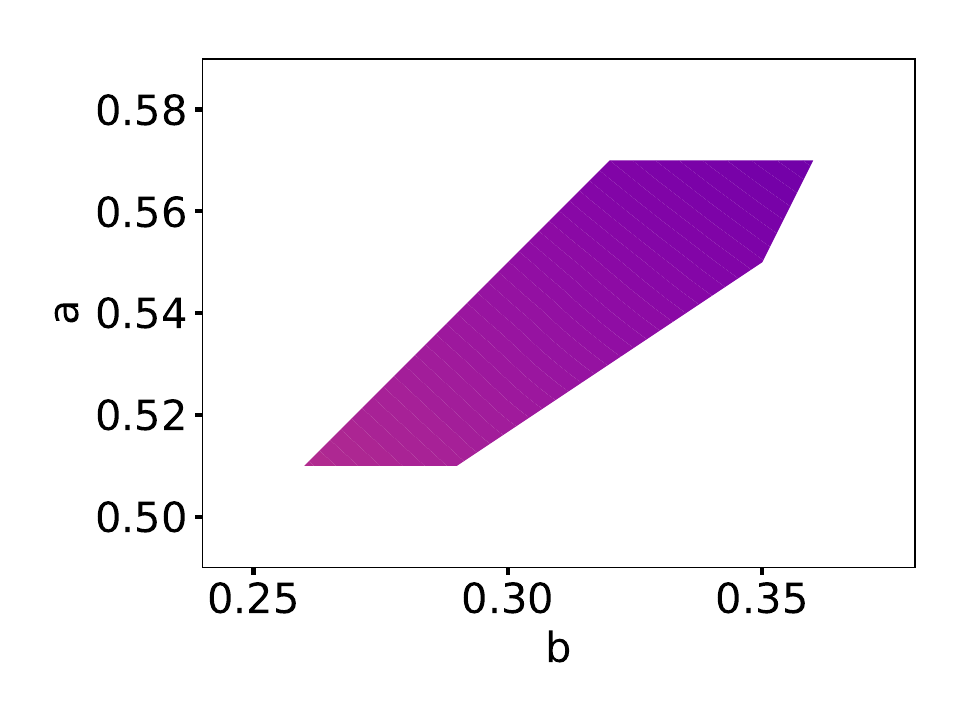}
        \caption{$q=5$; $r_E = 1.08$}
        \label{fig:Xi_5_dist_08}
    \end{subfigure}
    \begin{subfigure}[t]{0.24\textwidth}
        \includegraphics[width=1\linewidth]{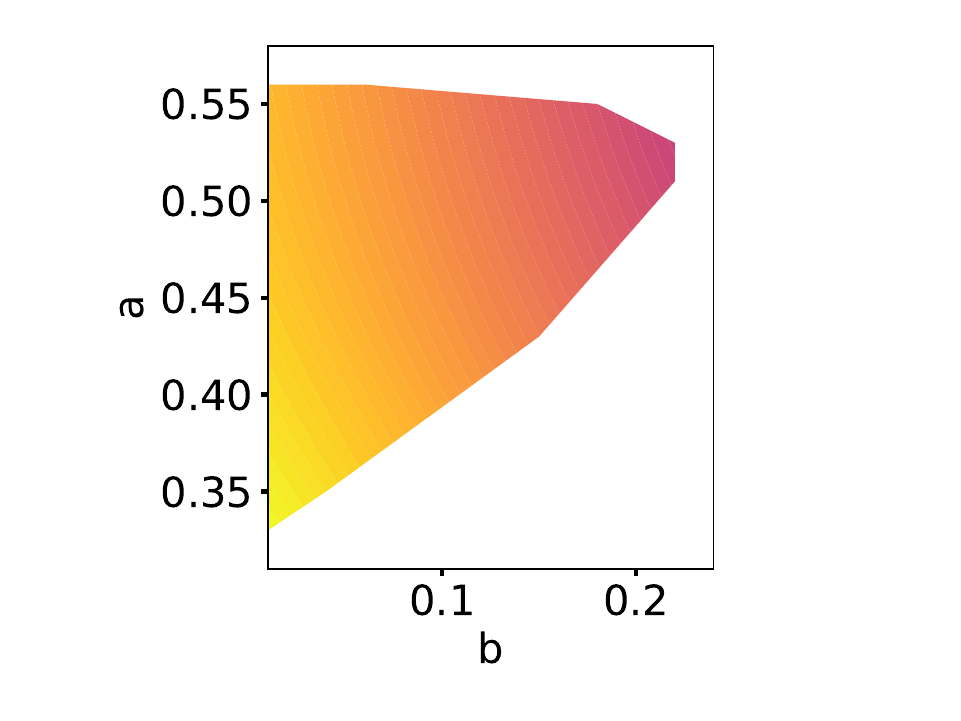}
        \caption{$q=6$; $r_E = 1.02$}
        \label{fig:Xi_6_dist_02}
    \end{subfigure}
    \begin{subfigure}[t]{0.24\textwidth}
        \includegraphics[width=1\linewidth]{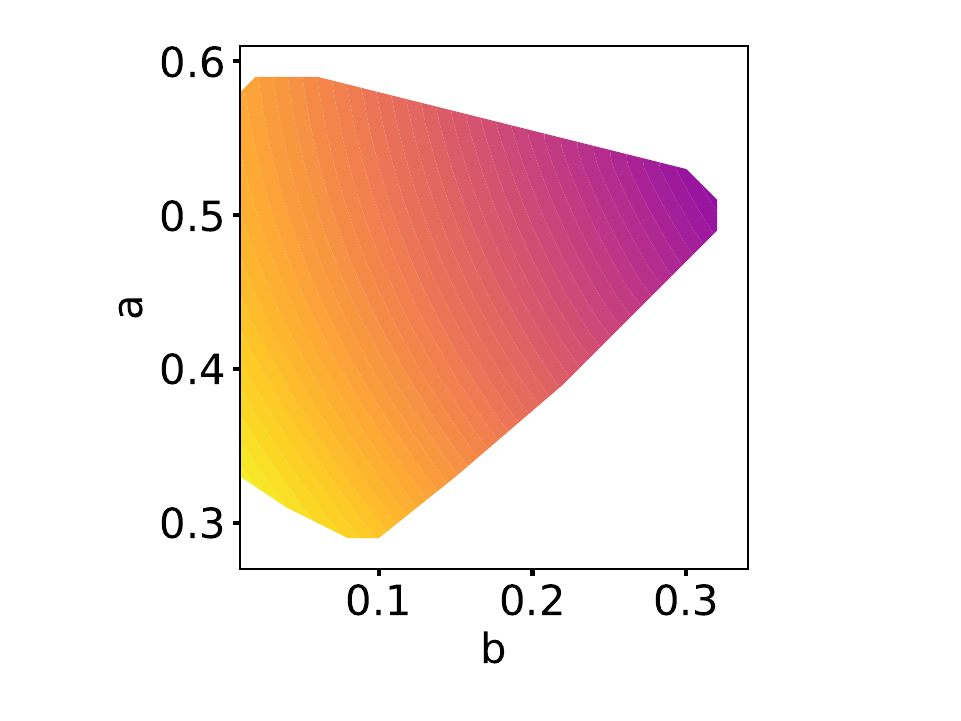}
        \caption{$q=6$; $r_E = 1.04$}
        \label{fig:Xi_6_dist_04}
    \end{subfigure}
    \begin{subfigure}[t]{0.24\textwidth}
        \includegraphics[width=1\linewidth]{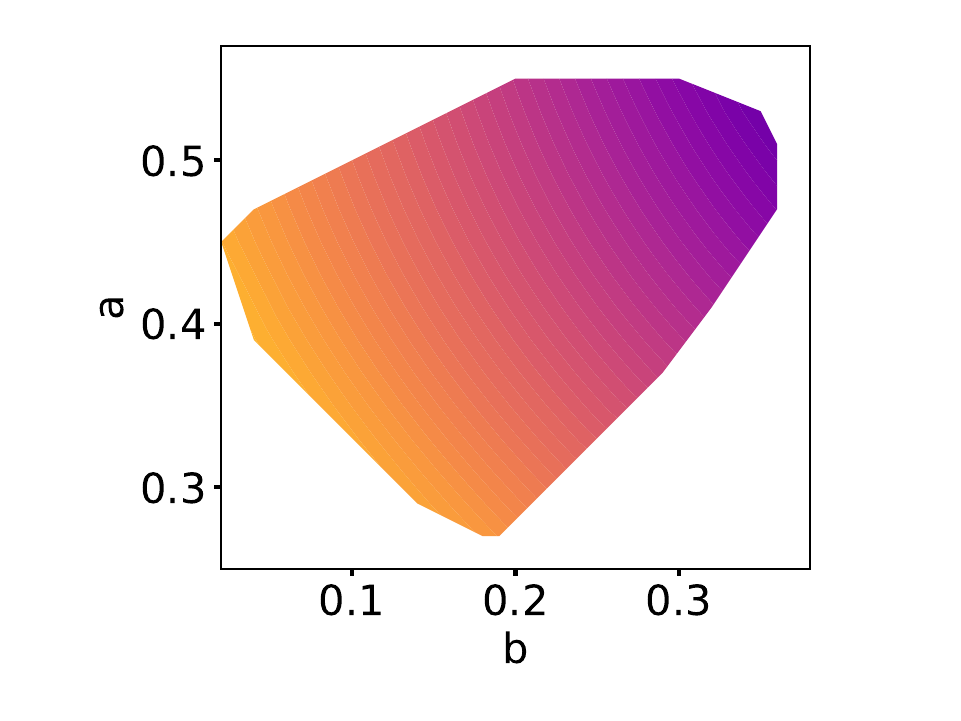}
        \caption{$q=6$; $r_E = 1.06$}
        \label{fig:Xi_6_dist_06}
    \end{subfigure}
    \begin{subfigure}[t]{0.24\textwidth}
        \includegraphics[width=1\linewidth]{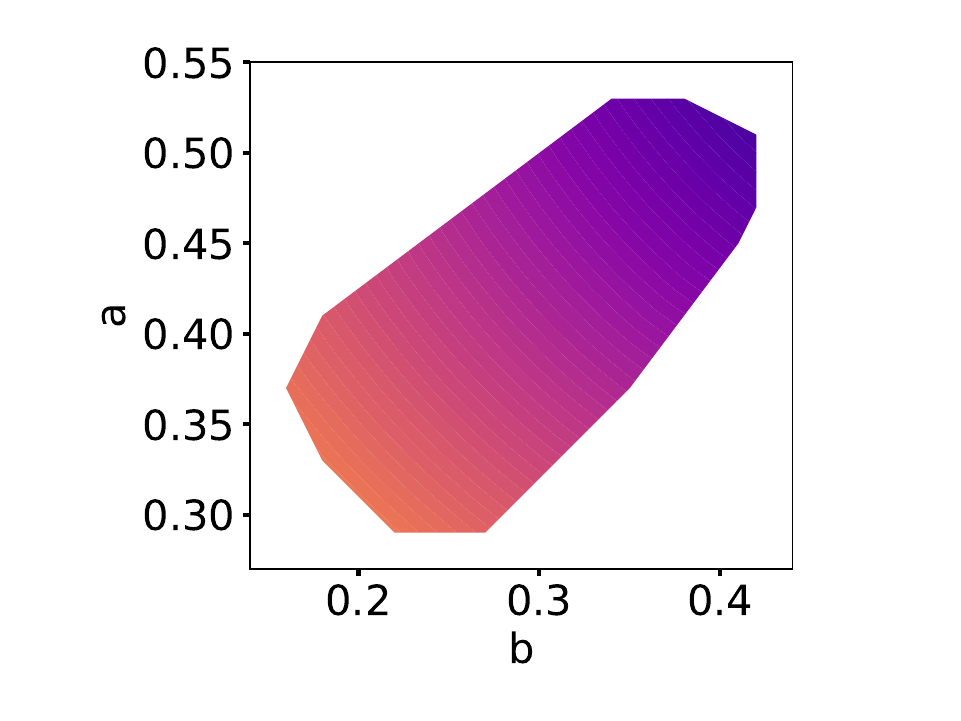}
        \caption{$q=6$; $r_E = 1.08$}
        \label{fig:Xi_6_dist_08}
    \end{subfigure}
    \begin{subfigure}[t]{0.24\textwidth}
        \includegraphics[width=1\linewidth]{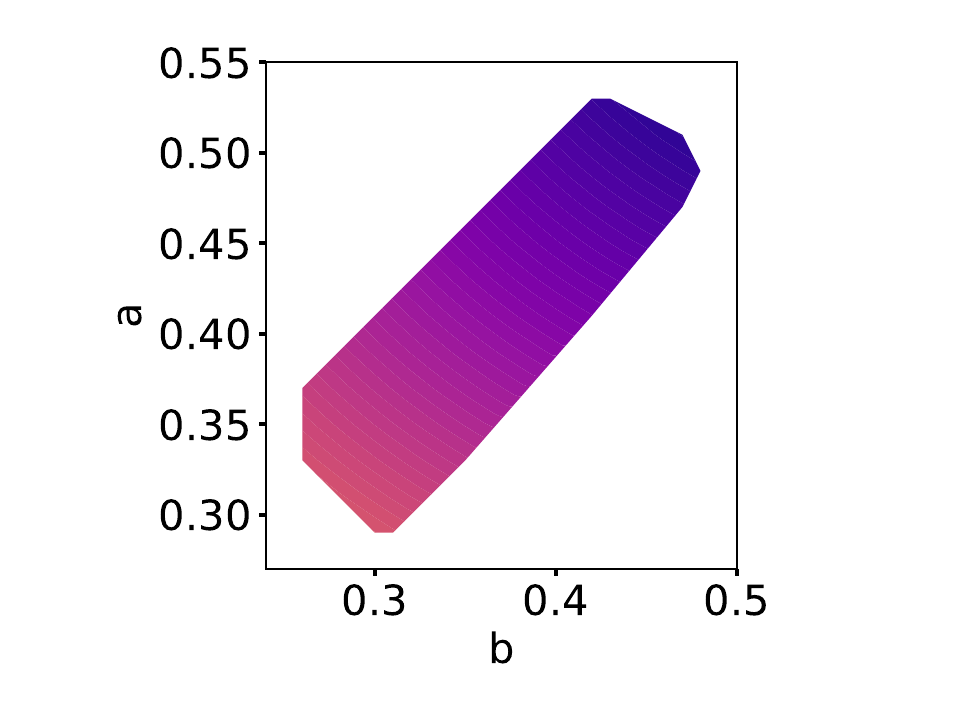}
        \caption{$q=6$; $r_E = 1.1$}
        \label{fig:Xi_6_dist_10}
    \end{subfigure}
    \begin{subfigure}[t]{0.24\textwidth}
        \includegraphics[width=1\linewidth]{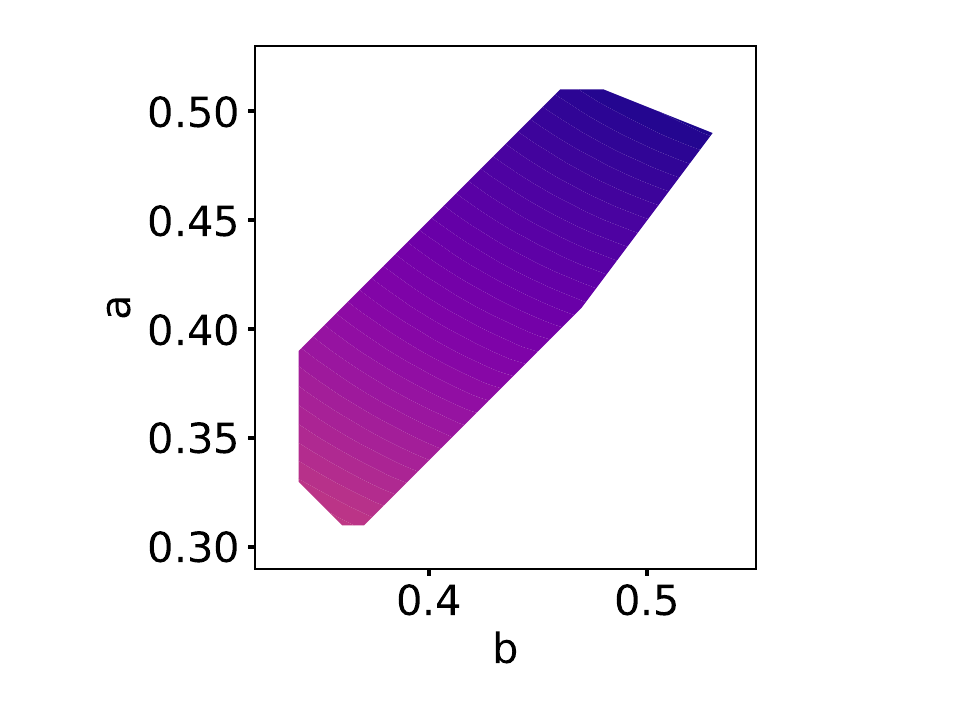}
        \caption{$q=6$; $r_E = 1.12$}
        \label{fig:Xi_6_dist_12}
    \end{subfigure}
    \begin{subfigure}[t]{0.24\textwidth}
        \includegraphics[width=1\linewidth]{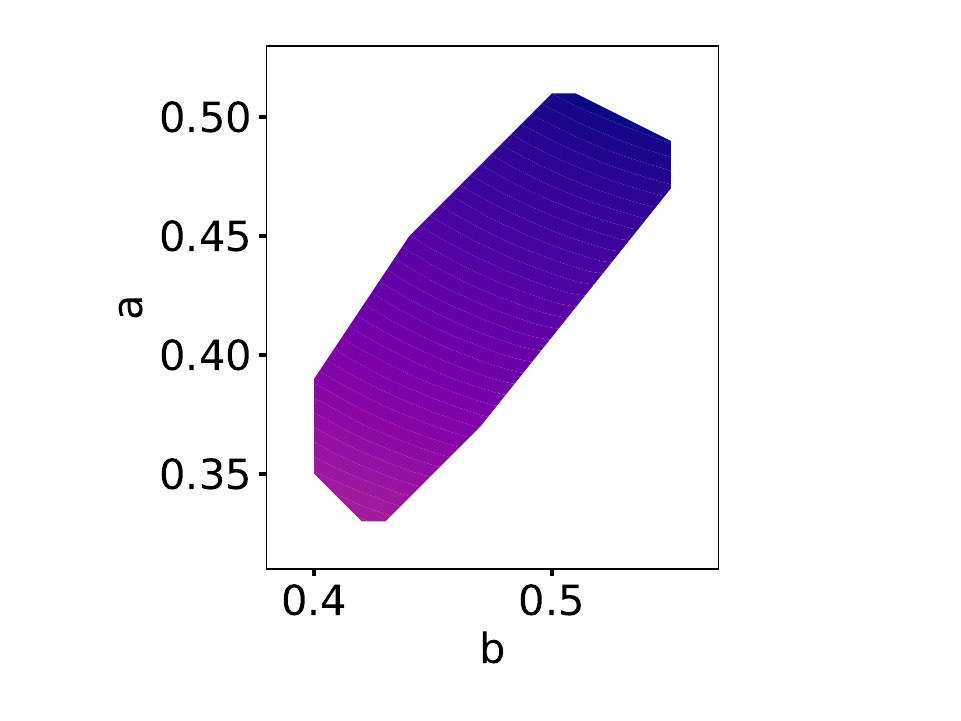}
        \caption{$q=6$; $r_E = 1.14$}
        \label{fig:Xi_6_dist_14}
    \end{subfigure}
    \begin{subfigure}[t]{0.24\textwidth}
        \centering
        \includegraphics[height=\heightof{\includegraphics[width=1\linewidth]{Xi_6_ellipses_distance_1_14.pdf}}]{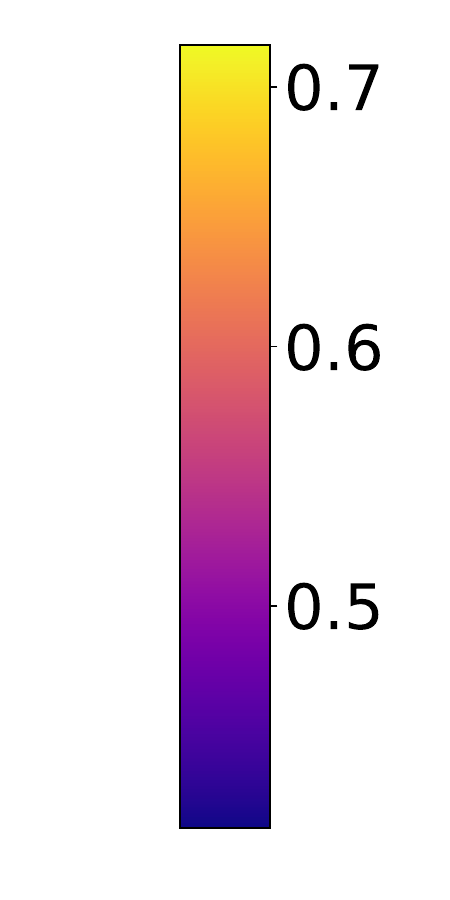}
        \caption{Colorbar}
        \label{fig:Xi_colorbar}
    \end{subfigure}
    \caption{Inner radius of the hypotrochoid based on model \eqref{ModelFitHypotrochoid} for a different amount of holes $q$ and distances of the holes to the origin $r_E$. The colorbar for these plots is visualized in (\subref{fig:Xi_colorbar})}
    \label{fig:Xi_model}
\end{figure}
Based on this model, we can calculate the maximum of the minimum pairwise distances between the nodal lines and the boundaries of the holes. For the given parameter set $r_E, a, b$, the distance between the nodal line and the boundary of the holes can be computed by 
\begin{equation}
    d_n\left({r_E},a,b\right) = {r_E} - b - {x_o}\left({r_E},a,b\right).
\end{equation}
A grid of the parameters $a$ and $b$ with a resolution of $0.001$ is generated for all $r_E$. The distance $d_n$ is then calculated for all grid points within the domains shown in the Figures \ref{fig:Parametersets_5} and \ref{fig:Parametersets_6}. The maximum distance for five holes is $d_n \approx 0.1978$ and is reached for $r_E = 1.04$, $a= 0.591$, $b=0.052$ with the eigenvalues $\lambda_1 \approx 5.2458$, $\lambda_2 \approx 13.1419$ and $\lambda_{3,4} \approx 13.1817$. The maximum distance for six holes is $0.2918$ and is reached for $r_E = 1.06$, $a= 0.45$, $b=0.021$ with the eigenvalues $\lambda_1 \approx 4.8967$, $\lambda_2 \approx 12.2005$ and $\lambda_{3,4} \approx 12.2334$. These two results demonstrate that there exist counterexamples with five and six holes, where the nodal line is sufficiently far away from the holes. This means that small variations of the parameters would not directly result in an intersection of the nodal line with the holes.

In the next step, we want to analyze, the behavior of the nodal line. Therefore, in Figure \ref{fig:Xa_model} the model of Equation \eqref{ModelFitHypotrochoid} is visualized for $x_a$. The maximum absolute error is $0.0064$. For the inner radius $x_i$, the maximum absolute error is $0.0032$, and $0.008$ for the ratio $\frac{x_i}{x_o}$. We see that for a given distance $r_E$ the smaller $a$ and $b$, the bigger is the outer radius of the hypotrochoid $x_o$. Moreover, the bigger the distance $r_E$, the smaller is in average the outer radius $x_o$ of the ellipses. 
\begin{figure}[h!]
    \centering
    \begin{subfigure}[t]{0.24\textwidth}
        \includegraphics[width=1\linewidth]{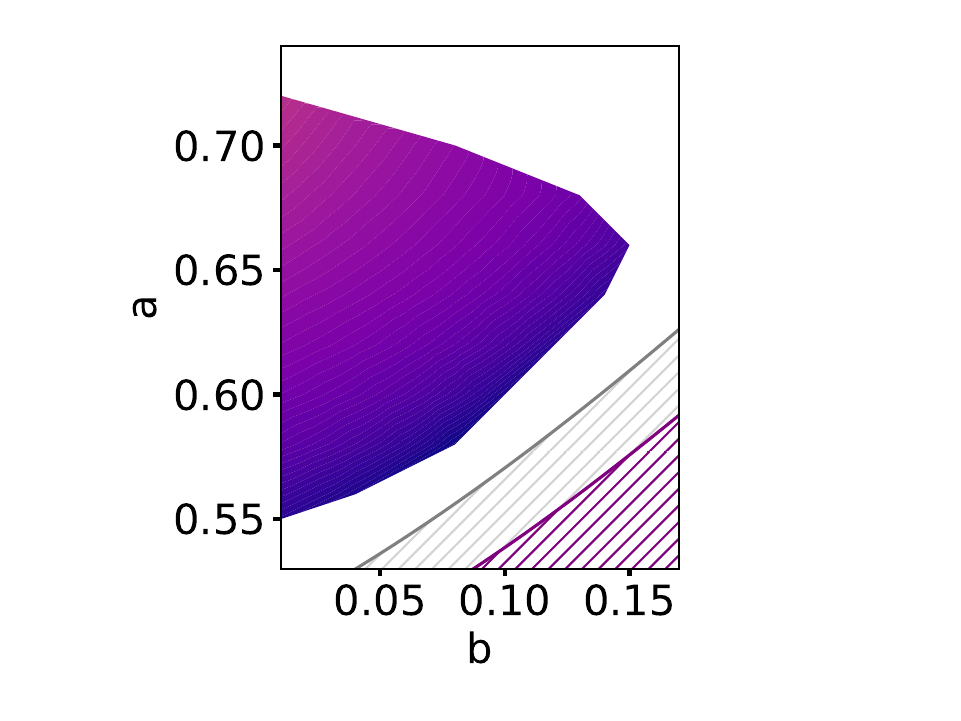}
        \caption{$q=5$; $r_E = 1.02$}
        \label{fig:Ratio_5_dist_02}
    \end{subfigure}
    \begin{subfigure}[t]{0.24\textwidth}
        \includegraphics[width=1\linewidth]{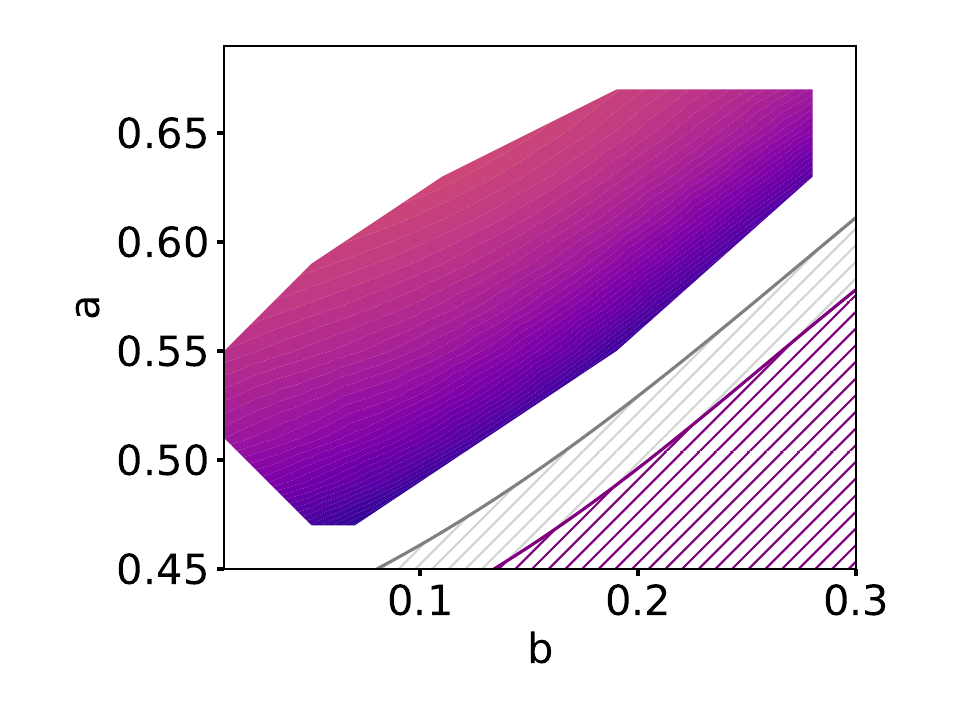}
        \caption{$q=5$; $r_E = 1.04$}
        \label{fig:Ratio_5_dist_04}
    \end{subfigure}
    \begin{subfigure}[t]{0.24\textwidth}
        \includegraphics[width=1\linewidth]{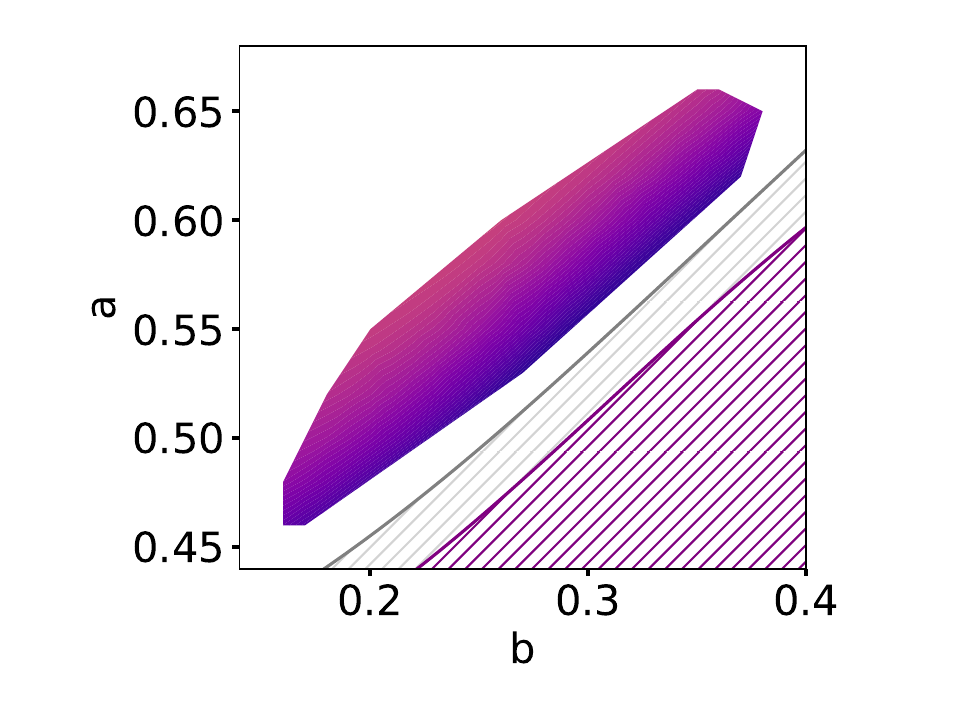}
        \caption{$q=5$; $r_E = 1.06$}
        \label{fig:Ratio_5_dist_06}
    \end{subfigure}
    \begin{subfigure}[t]{0.24\textwidth}
        \includegraphics[width=1\linewidth]{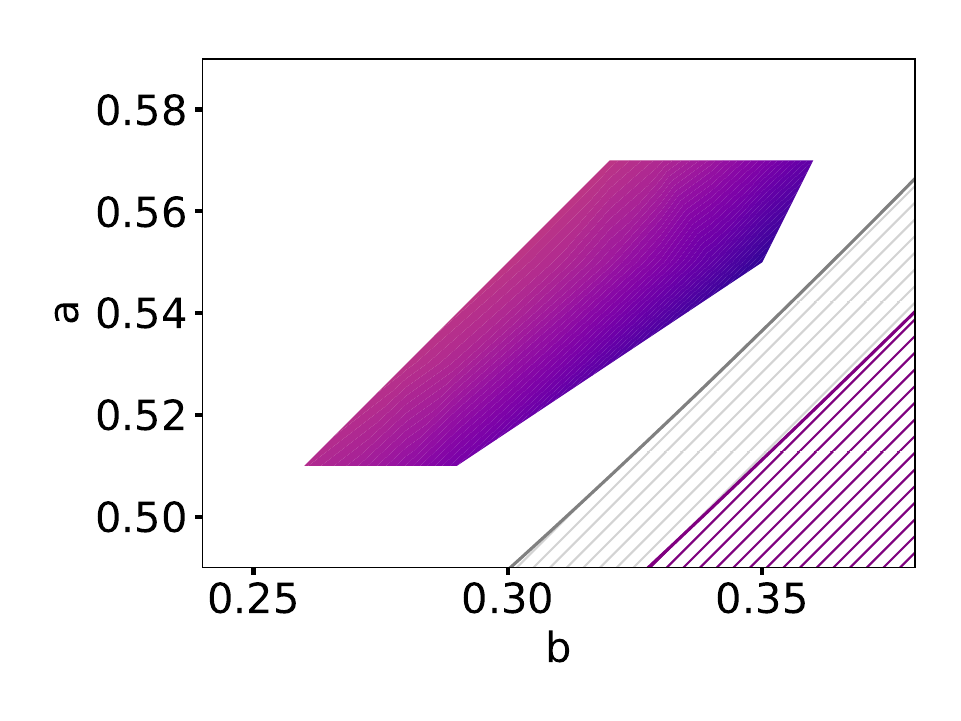}
        \caption{$q=5$; $r_E = 1.08$}
        \label{fig:Ratio_5_dist_08}
    \end{subfigure}
    \begin{subfigure}[t]{0.24\textwidth}
        \includegraphics[width=1\linewidth]{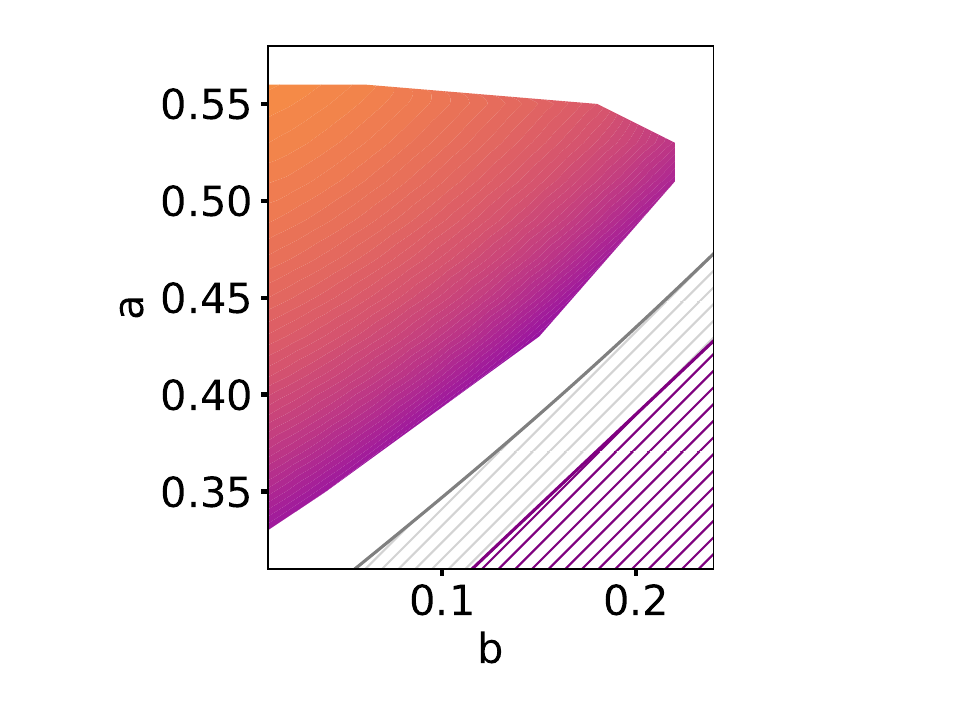}
        \caption{$q=5$; $r_E = 1.02$}
        \label{fig:Ratio_6_dist_02}
    \end{subfigure}
    \begin{subfigure}[t]{0.24\textwidth}
        \includegraphics[width=1\linewidth]{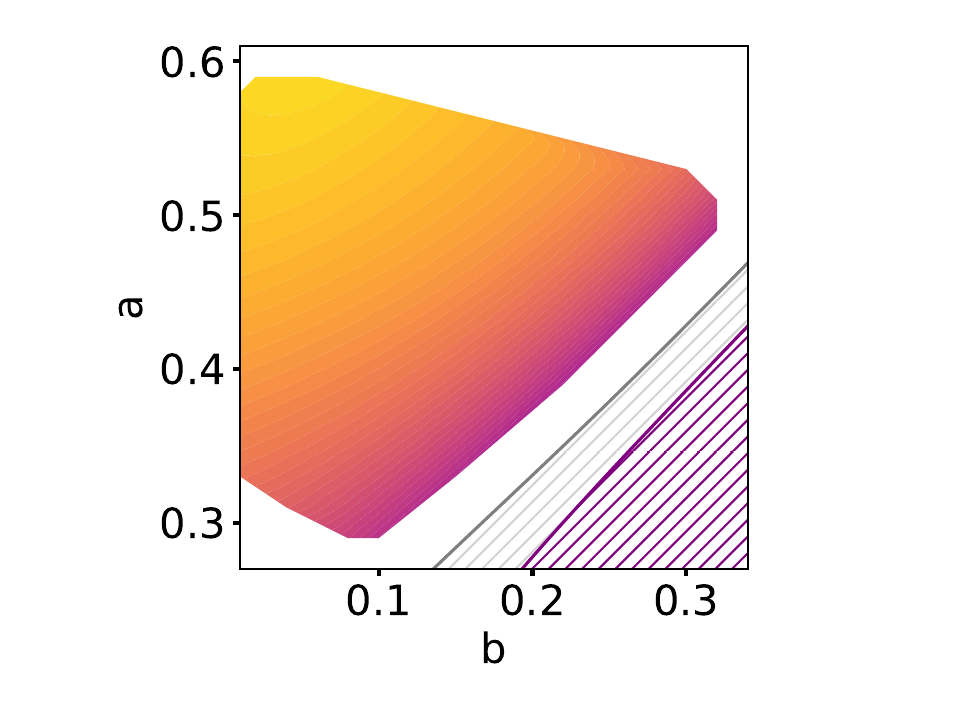}
        \caption{$q=6$; $r_E = 1.04$}
        \label{fig:Ratio_6_dist_04}
    \end{subfigure}
    \begin{subfigure}[t]{0.24\textwidth}
        \includegraphics[width=1\linewidth]{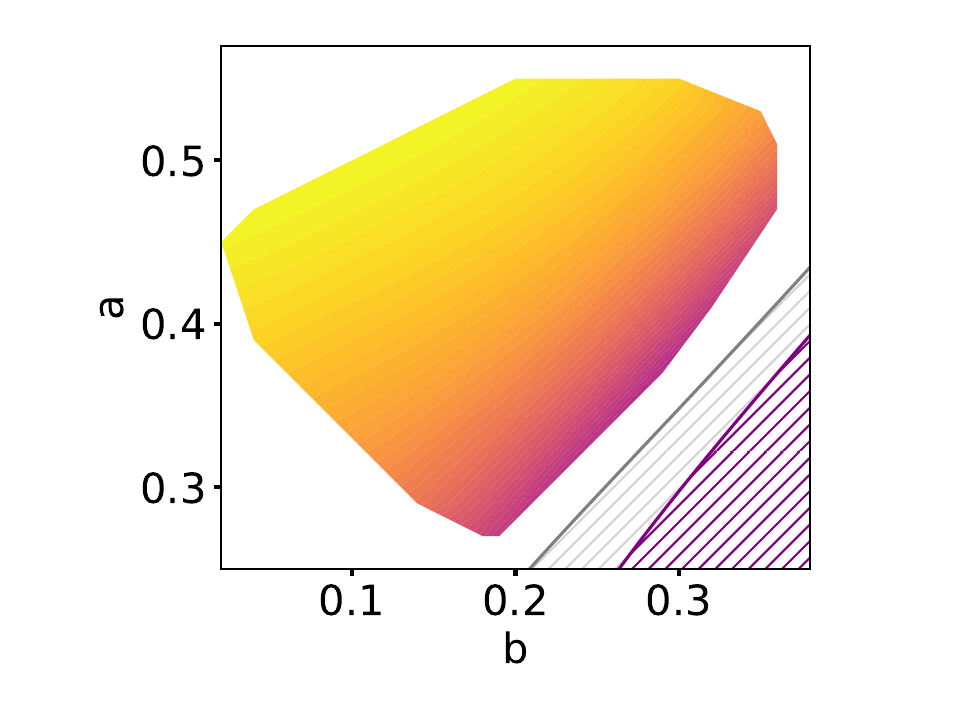}
        \caption{$q=6$; $r_E = 1.06$}
        \label{fig:Ratio_6_dist_06}
    \end{subfigure}
    \begin{subfigure}[t]{0.24\textwidth}
        \includegraphics[width=1\linewidth]{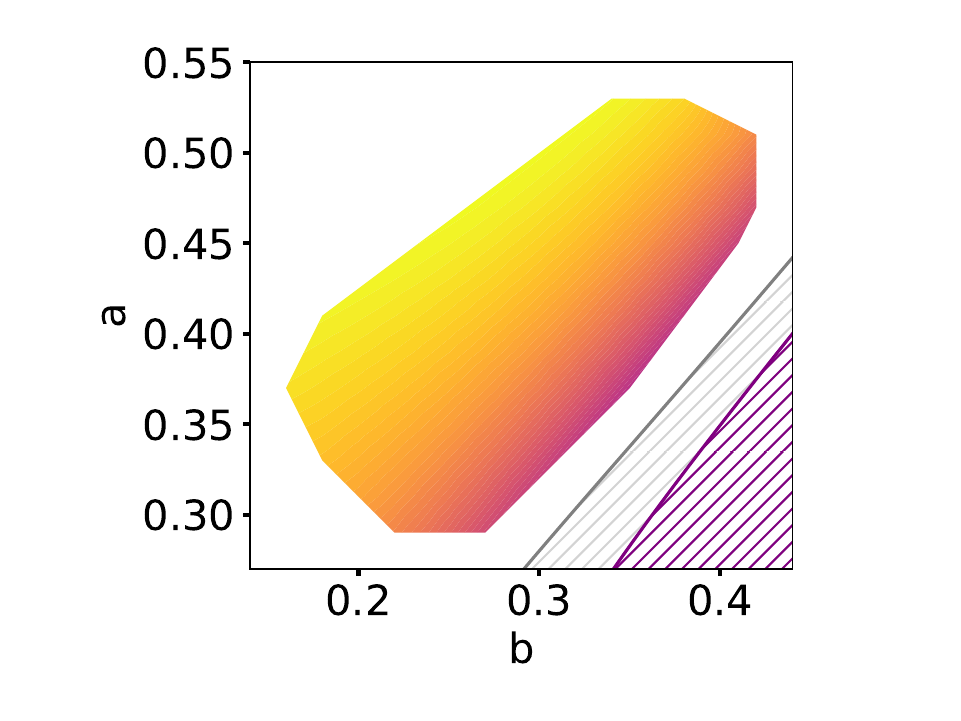}
        \caption{$q=6$; $r_E = 1.08$}
        \label{fig:Ratio_6_dist_08}
    \end{subfigure}
    \begin{subfigure}[t]{0.24\textwidth}
        \includegraphics[width=1\linewidth]{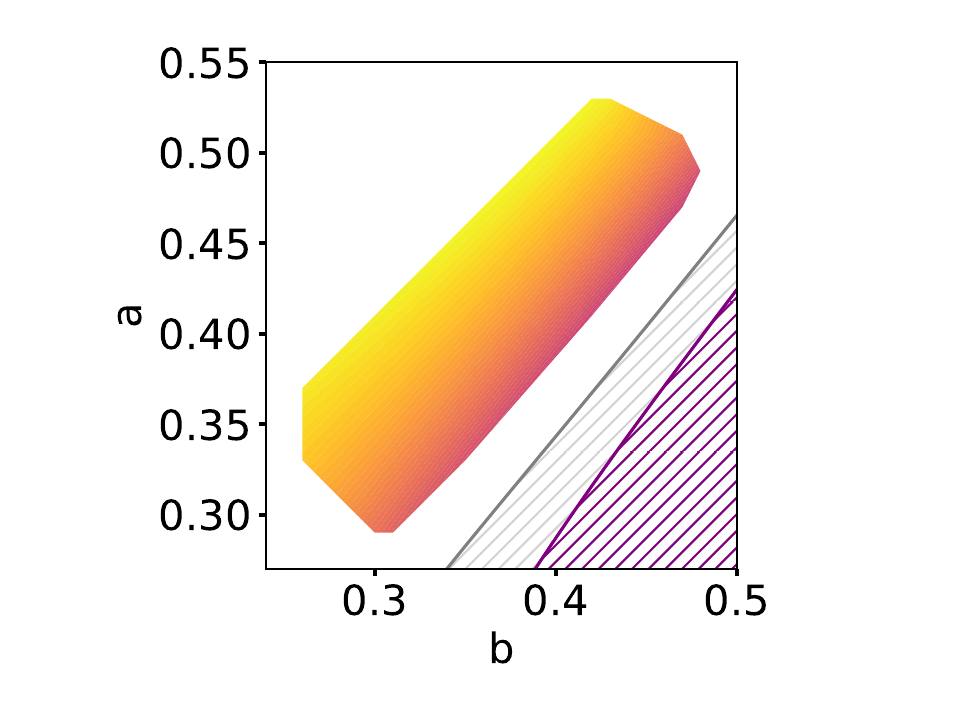}
        \caption{$q=6$; $r_E = 1.1$}
        \label{fig:Ratio_6_dist_10}
    \end{subfigure}
    \begin{subfigure}[t]{0.24\textwidth}
        \includegraphics[width=1\linewidth]{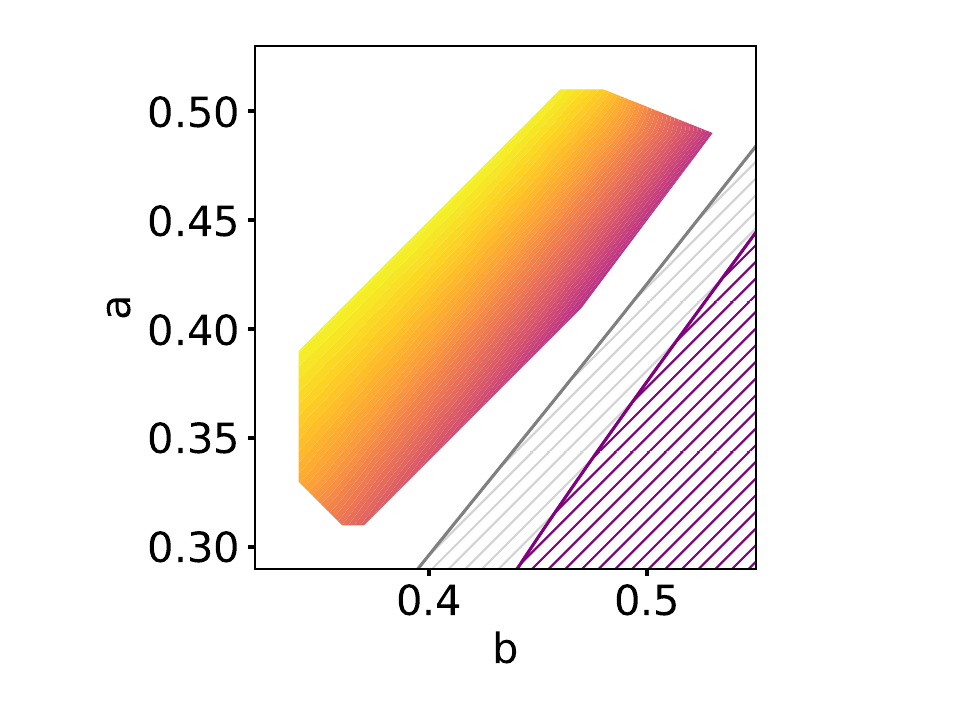}
        \caption{$q=6$; $r_E = 1.12$}
        \label{fig:Ratio_6_dist_12}
    \end{subfigure}
    \begin{subfigure}[t]{0.24\textwidth}
        \includegraphics[width=1\linewidth]{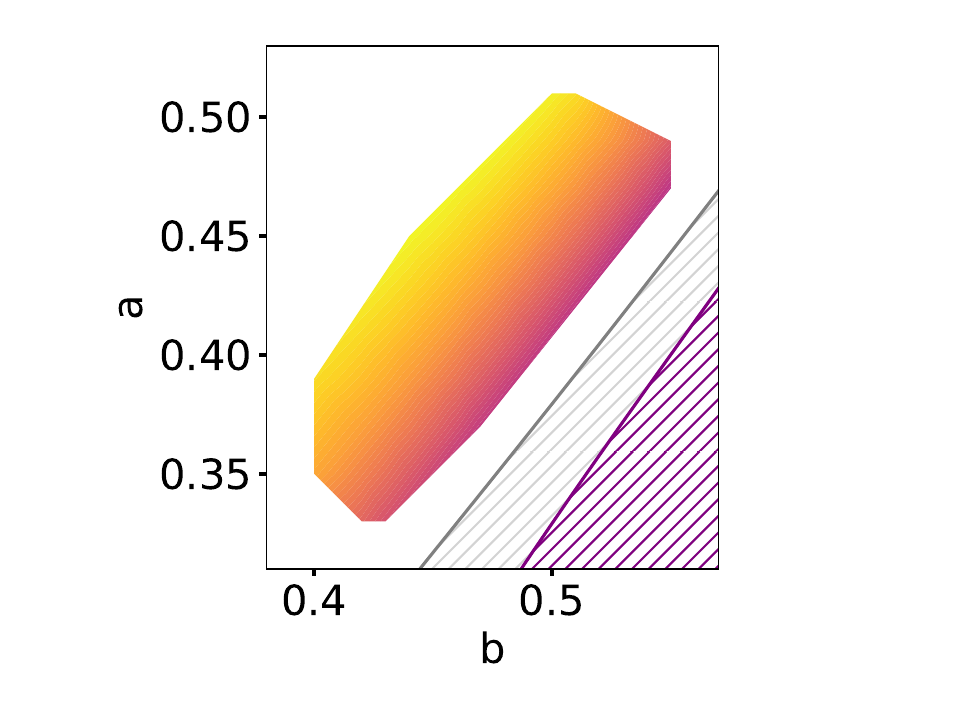}
        \caption{$q=6$; $r_E = 1.14$}
        \label{fig:Ratio_6_dist_14}
    \end{subfigure}
    \begin{subfigure}[t]{0.24\textwidth}
        \centering
        \includegraphics[height=\heightof{\includegraphics[width=1\linewidth]{Ratio_6_ellipses_distance_1_14.pdf}}]{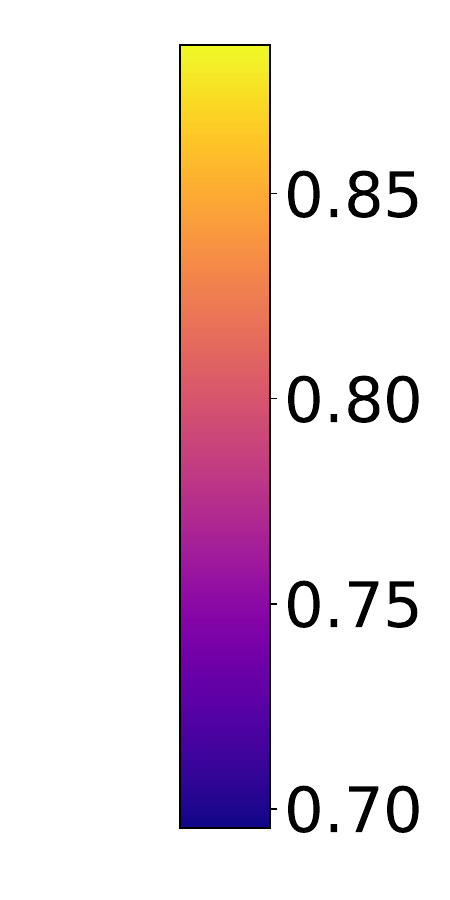}
        \caption{Colorbar}
        \label{fig:Ratio_colorbar}
    \end{subfigure}
    \caption{Ratio of the hypotrochoid $\frac{x_i}{x_o}$ based on model \eqref{ModelFitHypotrochoid} for a different amount of holes $q$ and distances of the holes to the origin $r_E$. The colorbar for these plots is visualized in (\subref{fig:Ratio_colorbar})}
    \label{fig:Ratio_model}
\end{figure}

The inner radius $x_i$ behaves differently, for a fixed distance $r_E$. The bigger $a$ and the smaller $b$, the bigger is the inner radius, see Figure \ref{fig:Xi_model}. The ratio of both radii can be seen in Figure \ref{fig:Ratio_model} and behaves similar to the inner radius. The bigger $a$ and the smaller $b$, the smaller is the ratio and the sharper pointed are the tips of the hypotrochoid. Additionally, it can be observed that the ratio decreases by reducing the amount of holes. These two observations are important, as they may explain why it is more difficult to find counter examples to Payne's conjecture with less holes. The ratio has to be $\frac{x_i}{x_o} \overset{!}{>} 1- \frac{2}{q}$, otherwise the hypotrochoid would intersect itself in the tips, see Figure \ref{HypotrochoidDifferentParameters}\subref{DomainHypotrochoid1}. A problem could occur, if the ratio would decrease faster by reducing the number of holes, than the lower level of the ratio decreases. 

Furthermore, we know $x_o \overset{!}{<} r_E - b$, since otherwise the hypotrochoid would intersect the ellipses. The parameter range that does not satisfy this condition is shown in Figure \ref{fig:Xa_model} as the hatched area. The edge of the hatched area is called $\Gamma_o$. It can be seen that the bound of the hatched area is roughly parallel to the parameter domains rejecting Payne's conjecture.
Additionally,  $x_i \overset{!}{>} 1- \frac{2}{q} x_o$ has to be valid, since otherwise the hypotrochoid would intersect itself in the tips, see Figure \ref{HypotrochoidDifferentParameters}\subref{DomainHypotrochoid1}. Invalid parameter sets of this condition are shown in Figure \ref{fig:Ratio_model} as the purple hatched area. The edge of this area is called $\Gamma_{ratio}$ This is as well parallel to the parameter domain rejecting Payne's conjecture. Additionally, we plotted the hatched line belonging to the condition $x_o \overset{!}{<} r_E - b$ in light gray. It can be seen that this condition is more strict, with the border $\Gamma_o$. Furthermore, the orthogonal distance starting from the purple line to the boundary of the parameter set rejecting Payne's conjecture, is about twice as long as the orthogonal distance of the boundaries of the two hatched domains.
This knowledge could be used to find counterexamples with less holes. First calculate the ratio and the outer radius of the hypotrochoid with the model given by \eqref{ModelFitHypotrochoid}. Then, use the two conditions to calculate the invalid parameter domains. At the end, go from the purple line $\Gamma_{ratio}$ twice the orthogonal distance of the purple line to the black line $\Gamma_o$ and a little bit further. There it may be possible to find counterexamples with less holes.
\subsection{Four ellipses}
\label{sec:4_holes}
Based on the results of the previous Section \ref{sec:BehaviorRatioNodalLine}, we aim to investigate whether it is possible to identify a region of parameter sets rejecting Payne's conjecture within a $0.01$ grid.
Therefore, we use the model defined by Equation \eqref{ModelFitHypotrochoid} to calculate the borders $\Gamma_o$ and $\Gamma_{ratio}$ described in Section \ref{sec:BehaviorRatioNodalLine} for four elliptical holes. With these borders it is possible to estimate where the searched area could be and for which distances it is not even necessary to prove whether one of the parameter sets rejects Payne's conjecture. For specific distances, the possible area outside of the hatched areas of the borders $\Gamma_0$ and $\Gamma_{ratio}$ is outside of the permitted area of Theorem \ref{TheoremOverlap}. Starting from the estimated points, we can search for the parameter sets that reject Payne's conjecture.
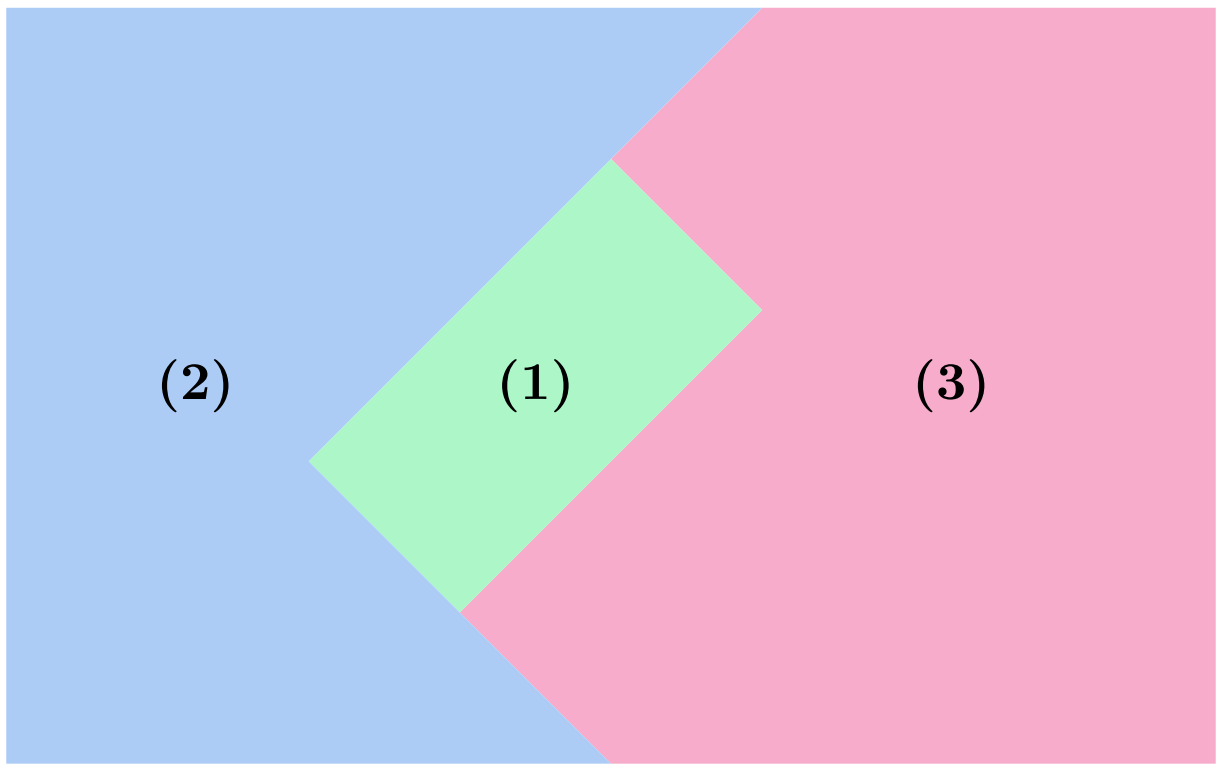
\begin{figure}[h!]
    \centering
        \tikzmath{\n = 5;
        \eps = 180/(\n+10);}
        \begin{tikzpicture}
            \definecolor{pastelgruen}{HTML}{acf6c8}
            \definecolor{rose}{HTML}{f6acca}
            \definecolor{lavenderblue}{HTML}{acccf6}
            \fill[rose] (1,-1) -- (2,-2) -- (6,-2) -- (6,3) -- (3,3) -- (2,2) -- (3,1);
            \fill[pastelgruen] (0,0) -- (1,-1) -- (3,1) -- (2,2);
            \fill[lavenderblue] (0,0) -- (3,3) -- (-2,3) -- (-2,-2) -- (2,-2);
            \draw[black, thick] (1.5,0.5) node {$\textbf{(1)}$};
            \draw[black, thick] (-0.75,0.5) node {$\textbf{(2)}$};
            \draw[black] (4.25,0.5) node {$\textbf{(3)}$};
        \end{tikzpicture}
    \caption{Types of nodal lines for a specific distance and specific shape of the domain.}
    \label{NodalLineTypes}
\end{figure}

We are using another result, which we figured out while we analyzed the nodal line for different parameter sets. In the previous Section \ref{sec:BehaviorRatioNodalLine}, we analyzed those nodal lines, that reject Payne's conjecture. Beside this, we figured out that the nodal line that does not reject Payne's conjecture has different shapes depending on the position of the parameters of the domain regarding the area rejecting Payne's conjecture. This can be seen in Figure \ref{NodalLineTypes}. The green part, labeled $(1)$, are the parameter sets for which the nodal lines does not touch the boundary. This region usually resembles a rotated rectangle or a part of it, if it is moved to the left side. The blue part, labeled $(2)$, represents the parameter sets for which the nodal line has a double eigenvalue and halves the domain, like in Figure \ref{fig:Dahne_all_results}(\subref{fig:Dahne_14_27}). This domain is obtained by extending the left boundaries of the green domain to the right side. The rest of the domain is the rose colored part, labeled $(3)$, and represents the parameter sets for which the nodal line has a simple eigenvalue and does not touch the outer boundary, like in Figure \ref{fig:Dahne_all_results}(\subref{fig:Dahne_11_27}).

As explained in the previous paragraph, we choose a parameter set for the permitted distances and consider the corresponding nodal line. Then, we search for the edge in the blue area, labeled $(2)$. If a parameter set that rejects Payne's conjecture exists, it would be right at the transition between the corner of the blue area and the rose area. Anyways, we could not find such a parameter set, which brings us to the assumption that there exist no area of parameter sets rejecting Payne's conjecture for a disc with four elliptical holes bigger than $0.01 \times 0.01 \times 0.01$.
\section{Summary, conclusion and future work}
\label{sec:Conclusion}
In this paper, a very efficient and highly parallelized algorithm to calculate the second Dirichlet eigenvalue of the negative Laplace operator was developed and presented. Therefore, first the formulation of the mathematical model using the boundary element method was discretized with the boundary element collocation method and the algorithm to solve the resulting non-linear eigenvalue problem with the Beyn algorithm was explained in detail. Based on these results, we were able to develop an algorithm that can find counterexamples to Payne's conjecture by itself for a given parameterization of the domain. This algorithm was parallelized with \texttt{MPI4py} in a well scaling way. With this, we could find counterexamples to Payne's conjecture with six and five holes by defining the domain as a disc with six and five elliptical holes. In the next step, we analyzed the parameter sets rejecting Payne's conjecture. We figured out numerically that the nodal line of such counterexamples always has the shape of an affine transformation of a hypotrochoid. The parameter sets are always simply-connected and the smaller the number of holes is without changing the parameterization of the domain, the smaller is the volume of the parameter sets. Further results can be seen in Section \ref{sec:ParametersetsRejectingPayne}. Furthermore, we achieved to find a model that describes the exterior radius of the hypotrochoid, the inner radius, and the ratio of both. Together with the restrictions that the ratio of the inner and exterior radius of a hypotrochoid has to be bigger than $\frac{x_i}{x_o} \geq 1 - \frac{2}{q}$, and the outer radius has to be smaller than $x_o < r_E - b$, we found out that if a parameter set rejecting Payne's conjecture exists for a specific distance of the holes, then the rejecting lines are parallel and close to the parameter set domain. Additionally, the boundary of the parameter domain that rejects Payne's conjecture is twice as far away from $\Gamma_{ratio}$, as the distance between the two lines, $\Gamma_o$ and $\Gamma_{ratio}$.
This information may help to find counterexamples, with less holes by starting at the line $\Gamma_{ratio}$, and moving orthogonally in the direction of the line $\Gamma_o$ twice the distance plus a little more.
We determined that it is impossible to find a counterexample to Payne's conjecture using a disc with four elliptical holes in a grid with step width $0.01$.

We conclude that it is possible to find numerically counterexamples to Payne's conjecture with six an five holes, thus reducing the minimum number of required holes. It is possible to restrict the domain of potential parameter set candidates, with the knowledge of parameter sets which reject Payne's conjecture for more holes. This makes it easier to find counterexamples with less holes. 

While this paper contributes meaningful insights, how the nodal line behaves for a different number of holes and different parameter sets of the domain in the two-dimensional case, it would be interesting to know how the nodal isosurface looks in the three-dimensional case using the code from \cite[ch. 5]{Kleefeld2011_Landweber}, \cite[ch. 3]{Kleefeld2012_ExteriorNeumann}. We suspect that the isosurface is an affine transformation of a star-shaped surface with cross-sections that resemble a hypotrochoid. This is analogous to the affine transformation of a hypotrochoid in the two-dimensional case. Moreover, an analytical proof, showing that the nodal line is always an affine transformation of a hypotrochoid in two dimensions, would verify our hypothesis, which was previously supported by numerical results.
Furthermore, we would like to know whether a counterexample to Payne's conjecture without a hole exists in the two-dimensional case. Additionally, we are looking for the minimal number of holes a three-dimensional domain needs to reject Payne's conjecture. 

\section*{Acknowledgements}
\label{sec:Acknowledgements}
We acknowledge the computing time on the supercomputer JUSUF (Jülich Support for Fenix) \cite{Vieth2021_JUSUF} at Jülich Supercomputing Centre through the project CSTMA.
\bibliography{Bibliography}

\end{document}